\documentclass[12pt,a4paper]{article}
\usepackage{graphicx}
\usepackage[latin1]{inputenc}
\usepackage[T1]{fontenc}
\usepackage[english]{babel}
\usepackage{amssymb,amsmath}

\newcommand{\thmheadercommand}[1]{\textbf{\scshape{}#1}}

\def\d{{\mathrm{d}}}
\def\N{{\mathbb{N}}}
\def\Z{{\mathbb{Z}}}

\renewcommand{\geq}{\geqslant}
\renewcommand{\leq}{\leqslant}

\def\eps{\varepsilon}
\renewcommand{\epsilon}{\varepsilon}
\renewcommand{\phi}{\varphi}

\DeclareMathOperator{\dist}{dist}

\newcommand{\abs}[1]{\left|\mskip1mu#1\right|}
\newcommand{\norm}[1]{\left\|#1\right\|}

\newcommand{\presgroup}[2]{\left\langle\,#1 \mid  #2\,\right\rangle}

\newcounter{prop}
\newcounter{defi}
\newcounter{thm}
\newcounter{lem}

\newenvironment{dem}[1][]{\noindent{\thmheadercommand{Proof#1}}\,\,--\,\,}{$\square$\medskip}

\newenvironment{enonce}[1]{\medskip\noindent{\thmheadercommand{#1}}\,\,--\,\,\begin{slshape}}{\end{slshape}\medskip}

\newenvironment{enonce2}[1]{\medskip\noindent{\thmheadercommand{#1}}\,\,--\,\,}{\medskip}

\newenvironment{defi}[1][]{\refstepcounter{prop}
\begin{enonce}{Definition \theprop{}#1}}{\end{enonce}}
\newenvironment{prop}[1][]{\refstepcounter{prop}
\begin{enonce}{Proposition \theprop{}#1}}{\end{enonce}}
\newenvironment{thm}[1][]{\refstepcounter{prop}
\begin{enonce}{Theorem \theprop{}#1}}{\end{enonce}}
\newenvironment{lem}[1][]{\refstepcounter{prop}
\begin{enonce}{Lemma \theprop{}#1}}{\end{enonce}}
\newenvironment{cor}[1][]{\refstepcounter{prop}
\begin{enonce}{Corollary \theprop{}#1}}{\end{enonce}}

\newenvironment{rem}[1][]{\refstepcounter{prop}
\begin{enonce2}{Remark \theprop{}#1}}{\end{enonce2}}

\newenvironment{defi*}[1][]{
\begin{enonce}{Definition#1}}{\end{enonce}}
\newenvironment{prop*}[1][]{
\begin{enonce}{Proposition#1}}{\end{enonce}}
\newenvironment{thm*}[1][]{
\begin{enonce}{Theorem#1}}{\end{enonce}}
\newenvironment{lem*}[1][]{
\begin{enonce}{Lemma#1}}{\end{enonce}}
\newenvironment{cor*}[1][]{
\begin{enonce}{Corollary#1}}{\end{enonce}}
\newenvironment{ex*}[1][]{
\begin{enonce}{Example#1}}{\end{enonce}}
\newenvironment{exo*}[1][]{
\begin{enonce2}{Exercise#1}}{\end{enonce2}}
\newenvironment{rem*}[1][]{
\begin{enonce2}{Remark#1}}{\end{enonce2}}

\usepackage{color}

\title{Sharp phase transition theorems for hyperbolicity of random groups}

\author{Yann Ollivier}

\begin{document}

\addtolength{\headsep}{-4mm}
\addtolength{\textheight}{4mm}

\maketitle

\begin{center}
With 31 illustrations
\end{center}

\begin{abstract}
We prove that in various natural models of a random quotient of a group,
depending on a density parameter,
for each hyperbolic group there is some critical density under which a
random quotient is still hyperbolic with high probability, whereas above
this critical value a random quotient is very probably trivial. We give
explicit characterizations of these critical densities for the various models.
\end{abstract}

\def\~{^{-1}}
\def\d{\partial}

\section*{Introduction}
\addcontentsline{toc}{section}{Introduction}

What does a generic group look like?

The study of random groups emerged from an affirmation of M.~Gromov that
``almost every group is hyperbolic'' (see~\cite{Gro1}).  More precisely,
fix $m$ and $N$ and consider the group $G$ presented by
$\presgroup{a_1,\ldots,a_m}{ r_1,\ldots,r_N}$ where the $r_i$'s are words
of length $\ell_i$ in the letters $a_i^{\pm 1}$. Then the ratio of the
number of $N$-tuples of words $r_i$ such that $G$ is hyperbolic, to the
total number of $N$-tuples of words $r_i$, tends to $1$ as $\inf
\ell_i\rightarrow \infty$. The first proof of this theorem was given by
A.Y.~Ol'shanski\u{\i} in~\cite{Ols1}, and independently by C.~Champetier
in~\cite{Ch1}, thus confirming Gromov's statement.

Later, M.~Gromov introduced (cf.~\cite{Gro2}) a finer model of random
group, in which the number $N$ of relators is allowed to be much bigger.

This model goes as follows: Choose at random $N$ cyclically reduced words
of length $\ell$ in the letters $a_i^{\pm 1}$, uniformly among the set of
all such cyclically reduced words (recall a word is called \emph{reduced}
if it does not contain a sequence of the form $a_ia_i\~$ or $a_i\~a_i$
and \emph{cyclically reduced} if moreover the last letter is not the
inverse of the first one). Let $R$ be the (random) set of these $N$
words, the random group is defined by the presentation $\presgroup
{a_1,\ldots,a_m}{R}$.

Let us explain how $N$ is taken in this model. There are
$(2m)(2m-1)^{\ell-1}\approx(2m-1)^\ell$ reduced words of length $\ell$.
We thus take $N=(2m-1)^{d\ell}$ for some number $d$ between $0$ and $1$
called \emph{density}.

The theorem stated by Gromov in this context expresses a sharp phase transition
between hyperbolicity and triviality, depending on the asymptotics of the
number of relators taken, which is controlled by the density parameter $d$.

\begin{thm}[ (M.~Gromov, \cite{Gro2})]
\label{G12}
Fix a density $d$ between $0$ and $1$. Choose a length $\ell$ and pick at
random a set $R$ of $(2m-1)^{d\ell}$ uniformly chosen cyclically reduced words of
length $\ell$ in the letters $a_1^{\pm 1},\ldots,a_m^{\pm 1}$.

If $d<1/2$ then the probability that the presentation $\presgroup
{a_1,\ldots,a_m}{R}$ defines an infinite hyperbolic group tends to
$1$ as $\ell\rightarrow\infty$.

If $d>1/2$ then the probability that the presentation $\presgroup
{a_1,\ldots,a_m}{R}$ defines the group $\{e\}$ or $\Z/2\Z$ tends to
$1$ as $\ell\rightarrow\infty$.
\end{thm}

There was a small mistake in the original proof of Gromov: the proof uses
van Kampen diagrams, and the case when some relator appears several times
in a given van Kampen diagram was forgotten (this mistake was apparently
first detected by R.~Kenyon); when no relator appears twice there is
much more independence in the probabilities and the proof is easier.
A complete proof of this theorem is included below
(section~\ref{standardcase}).

Let us discuss the intuition behind this model. What does the
density parameter $d$ mean?  Following the excellent exposition of Gromov
in~\cite{Gro2}, we assimilate $d\ell$ to a dimension. That is, $d\ell$
represents the number of ``equations'' we can impose on a random word so
that we still have a reasonable chance to find such a word in a set of
$(2m-1)^{d\ell}$ randomly chosen words (compare to the basic intersection
theory for random sets stated in section~\ref{triv}).

For example, for large $\ell$, in a set of $2^{d\ell}$ randomly chosen
words of length $\ell$ in the two letters ``a'' and ``b'', there will
probably be some word beginning with $d\ell$ letters ``a''. (This is a
simple exercise.)

As another example, in a set of $(2m-1)^{d\ell}$ randomly chosen words on
$a_i^{\pm 1}$, there will probably be two words having the same first
$2d\ell$ letters, but no more. In particular, if $d<1/12$ then the set of
words will satisfy the small cancellation property $C'(1/6)$
(see~\cite{GH} for definitions). But as soon as $d>1/12$, we are far
from small cancellation, and as $d$ approaches $1/2$ we have arbitrarily
big cancellation.

\bigskip

The purpose of this work is to give similar theorems in a more general
situation. The theorem above states that a random quotient of the free
group $F_m$ is hyperbolic. A natural question is: does a random quotient
of a hyperbolic group stay hyperbolic?

This would allow in particular to iterate the operation of taking a random
quotient. This kind of construction is at the heart of the ``wild'' group
constructed in~\cite{Gro4}.

Our theorems precisely state that for each hyperbolic group (with
``harmless'' torsion), there is a
critical density $d$ under which the quotient stays hyperbolic, and above
which it is probably trivial.
Moreover, this critical density can be characterized in terms of
well-known numerical quantities depending on the group.

We need a technical assumption of ``harmless'' torsion (see
Definition~\ref{harmlesstorsion}). Hyperbolic groups with harmless
torsion include torsion-free groups, free products of torsion-free groups
and/or finite groups (such as $\mathrm{PSL}_2(\Z)$), etc. This assumption is necessary: Appendix~\ref{appcex}
proves that Theorem~\ref{main} does not hold for some hyperbolic groups
with harmful torsion
\footnote{These results were announced in~\cite{Oll1} without this
assumption. I would like to thank Prof.~A.Yu.~Ol'shanski\u{\i} for having
pointed an error in the first version of this manuscript regarding the
treatment of torsion, which led to this assumption and
to Appendix~\ref{appcex}.}.

\bigskip

There are several ways to generalize Gromov's theorem: a good replacement
in a hyperbolic group for reduced words of length $\ell$ in a free group
could, equally likely, either be reduced words of length $\ell$ again, or
elements of norm $\ell$ in the group (the norm of an element is the
minimal length of a word equal to it). We have a theorem for each of
these two cases. We also have a theorem for random quotients by uniformly
chosen plain words (without any assumption).

In the first two versions, in order to have the number of reduced, or
geodesic, words of length $\ell$ tend to infinity with $\ell$, we have to
suppose that $G$ is not elementary. There is no problem with the case of
a quotient of an elementary group by plain random words (and the critical
density is $0$ in this case).

Let us begin with the case of reduced words, or cyclically reduced words
(the theorem is identical for these two variants).

We recall the definition and basic properties of the cogrowth $\eta$ of a
group $G$ in section~\ref{defexps} below. Basically, if $G$ is not free,
the number of reduced words of length $\ell$ which are equal to $e$ in
$G$ behaves like $(2m-1)^{\eta\ell}$. For a free group, $\eta$ is
(conventionally, by the way) equal to $1/2$. It is always at least $1/2$.

\begin{thm}[ (Random quotient by reduced words)]
\label{reducedmain}
Let $G$ be a non-elementary hyperbolic group with harmless torsion, generated by the elements $a_1,\ldots,a_m$.
Fix a density $d$ between $0$ and $1$. Choose a length $\ell$ and pick at
random a set $R$ of $(2m-1)^{d\ell}$ uniformly chosen (cyclically) reduced words of
length $\ell$ in $a_i^{\pm 1}$. Let $\langle R\rangle$ be the normal
subgroup generated by $R$.

Let $\eta$ be the cogrowth of the group $G$.

If $d<1-\eta$, then, with probability tending to $1$ as
$\ell\rightarrow\infty$, the quotient $G/\langle R\rangle$ is
non-elementary hyperbolic.

If $d>1-\eta$, then, with probability tending to $1$ as
$\ell\rightarrow\infty$, the quotient $G/\langle R\rangle$ is
either $\{e\}$ or $\Z/2\Z$.
\end{thm}

\bigskip

We go on with the case of elements on the $\ell$-sphere of the group.

In this case, for the triviality part of the theorem, some small-scale
phenomena occur, comparable to the occurrence of $\Z/2\Z$ above (think of
a random quotient of $\Z$ by any number of elements of norm $\ell$).  In
order to avoid them, we take words of norm not exactly $\ell$, but of
norm between $\ell-L$ and $\ell+L$ for some fixed $L>0$ ($L=1$ is
enough).

\begin{thm}[ (Random quotient by elements of a sphere)]
\label{geodmain}
Let $G$ be a non-elementary hyperbolic group with harmless torsion, generated by the elements $a_1,\ldots,a_m$.
Fix a density $d$ between $0$ and $1$. Choose a length $\ell$.

Let $S^\ell$ be the set of elements of $G$ which are of norm between
$\ell-L$ and $\ell+L$ with respect to $a_1^{\pm 1},\ldots,a_m^{\pm 1}$ (for some fixed
$L>0$). Let $N$ be the number of elements of $S^\ell$.

Pick at random a set $R$ of $N^d$ uniformly chosen elements of
$S^\ell$. Let $\langle R\rangle$ be the normal subgroup generated by $R$.

If $d<1/2$, then, with probability tending to $1$ as
$\ell\rightarrow\infty$, the quotient $G/\langle R\rangle$ is
non-elementary hyperbolic.

If $d>1/2$, then, with probability tending to $1$ as
$\ell\rightarrow\infty$, the quotient $G/\langle R\rangle$ is $\{e\}$.
\end{thm}

\bigskip

The two theorems above were two possible generalizations of Gromov's
theorem. One can wonder what happens if we completely relax the assumptions
on the words, and take in our set $R$ any kind of words of length $\ell$
with respect to the generating set. The same kind of theorem still
applies, with of course a smaller critical density.

The gross cogrowth $\theta$ of a group is defined in
section~\ref{defexps} below.  Basically, $1-\theta$ is the exponent (in
base $2m$) of return to $e$ of the random walk on the group. We always
have $\theta> 1/2$.

Now there are $(2m)^\ell$ candidate words of length $\ell$, so we define
density with respect to this number.

\begin{thm}[ (Random quotient by plain words)]
\label{main}
Let $G$ be a hyperbolic group with harmless torsion, generated by the elements $a_1,\ldots,a_m$.
Fix a density $d$ between $0$ and $1$. Choose a length $\ell$ and pick at
random a set $R$ of $(2m)^{d\ell}$ uniformly chosen words of
length $\ell$ in $a_i^{\pm 1}$. Let $\langle R\rangle$ be the normal
subgroup generated by $R$.

Let $\theta$ be the gross cogrowth of the group $G$.

If $d<1-\theta$, then, with probability tending to $1$ as
$\ell\rightarrow\infty$, the quotient $G/\langle R\rangle$ is
non-elementary hyperbolic.

If $d>1-\theta$, then, with probability tending to $1$ as
$\ell\rightarrow\infty$, the quotient $G/\langle R\rangle$ is
either $\{e\}$ or $\Z/2\Z$.
\end{thm}

\paragraph{Precisions on the models.} Several points in the theorems
above are left for interpretation.

There is a slight difference between choosing $N$ times a random word and
having a random set of $N$ words, since some word could be chosen
several times. But for $d<1/2$ the probability that a word is chosen
twice is very small and the difference is negligible; anyway
this does not affect our statements at all, so both interpretations are
valid.

Numbers such as $(2m)^{d\ell}$ are not necessarily integers. We can either
take the integer part, or choose two constants $C_1$ and $C_2$ and
consider taking the number of words between $C_1(2m)^{d\ell}$ and
$C_2(2m)^{d\ell}$. Once more this does not affect our statements at all.

The case $d=0$ is peculiar since nothing tends to infinity. Say that a
random set of density $0$ is a random set with a number of elements
growing subexponentially in $\ell$ (e.g.\ with a constant number of
elements).

The possible occurrence of $\Z/2\Z$ above the critical density only
reflects the fact that it may be the case that a presentation of $G$ has
no relators of odd length (as in the free group). So, when quotienting by
words of even length, at least $\Z/2\Z$ remains.

\paragraph{Discussion of the models.} Of course, the three theorems given
above are not proven separately, but are particular cases of a more
general (and more technical!) theorem. This theorem is stated in
section~\ref{thetheorem}.

Our general theorem deals with random quotients by words picked under
a given probability measure. This measure does not need to be uniform,
neither does it necessarily charge words of only one given length. It
has to satisfy some natural (once the right terminology is given...)
axioms. The axioms are stated in section~\ref{theaxioms}, and the quite
sophisticated terminology for them is given in
section~\ref{somevocabulary}.

For example, using these axioms it is easy to check that
Theorem~\ref{geodmain} still holds when quotienting by words taken in the
ball rather than in the sphere, or that taking a random quotient by
reduced words or by cyclically reduced words is (asymptotically) the
same, with the same critical density.

It is also possible to take quotients by words of different lengths, but
our method imposes that the ratio of the lengths be bounded. This is a
restriction due to the geometric nature of some parts of the argument,
which rely on the hyperbolic local-global principle, using metric
properties of the Cayley complex of the group (cf.~appendix~\ref{CHGP}).

In the case of various lengths, density has to be defined as the supremum
of the densities at each length.

The very first model of random group given in this article (the one used
by Ol'shanski\u{\i} and Champetier), with a constant number of words of
prescribed lengths, is not the case $d=0$ of our model, since in this
model the ratio of lengths can be unbounded, which completely prevents
the use of some geometric methods.  However, this model can probably be
obtained by iterating the process of taking a random quotient at $d=0$,
or by using the relative small cancellation techniques later developed by
Delzant in~\cite{D} and by Gromov in~\cite{Gro4}.

But another model encountered in the literature, which consists in
uniformly picking a fixed number of words of length between $1$ and
$\ell$, satisfies easily our axioms, as it is almost exactly our case
$d=0$. Indeed there are so much more words of length close to $\ell$ than
close to $0$, that the elements taken under this model are of length
comprised between $(1-\eps)\ell$ and $\ell$ for any $\eps$.

Whereas random plain words or random reduced words can be easily
constructed independently of the group, it could seem difficult, at first
glance, to take a quotient by random elements of a sphere. Let us simply
recall (cf.~\cite{GH}) that in a hyperbolic group, it is possible to
define for each element a normal geodesic form, and that there exists a
finite automaton which recognizes exactly the words which are normal
forms of elements of the group.

Note that all our models of random quotients depend on a generating
subset. For example, adding ``false generators'' (i.e.\ generators equal
to $e$) to our generating sets makes the cogrowth and gross cogrowth
arbitrarily close to $1$, thus the critical density for reduced words and
plain words arbitrarily small. The case of random quotients by elements
of the $\ell$-sphere seems to be more robust.

In~\cite{Z}, A.~\.Zuk proves that a random quotient of the free group by
reduced words at density greater than $1/3$ has property T. As a random
quotient of any group is the quotient of a random quotient of
the free group by the relations defining the initial group,
this means that the random quotients we consider possess property T as
well for reduced words and densities above $1/3$.

\paragraph{Other developments on generic properties of groups.} Other
properties of generic groups have been studied under one or another model
of random group. Besides hyperbolicity, this includes topics such as
small cancellation properties, torsion elements, topology of the
boundary, property T, the fact that most subgroups are free, planarity of
the Cayley graph, or the
isomorphism problem; and more
are to come.  See for example~\cite{Ch1}, \cite{AO}, \cite{A}, \cite{Z},
\cite{AC}, \cite{KS}.

Random groups have been used by M.~Gromov to construct a ``wild'' group
related to $C^\star$-algebraic conjectures, see~\cite{Gro4}.

The use of generic properties of groups also led to an announcement of an
enumeration of one-relator groups up to isomorphism, see~\cite{KSS}.

In a slightly different approach, the study of what a generic group looks
like has very interesting recent developments: genericity can also be
understood as a topological (rather than probabilistic) property in the
space of all finite type groups. See for example the work of
C.~Champetier in~\cite{Ch3}.

In all these works, properties linked to hyperbolicity are ubiquitous.

\paragraph{Acknowledgements.}
I would like to thank, in alphabetical order, Thomas Del\-zant, Étienne
Ghys, Misha Gromov, Claire Kenyon, Richard Kenyon, Pierre Pansu, Panos
Papasoglu, Frédéric Paulin and Andrzej \.Zuk for instructive talks and
comments.

Special thanks to Prof.~A.Yu.~Ol'shanksi\u{\i}, who kindly pointed out an error
in the treatment of torsion in a previous version of the manuscript (the
assumption of harmless torsion did not appear), which led to the
counterexamples of Appendix~\ref{appcex}, as well as for careful reading
and suggestions for the text.

Part of the ideas of this work emerged during my stay at the École
normale supérieure of Lyon in April 2002, at the invitation of Andrzej
\.Zuk. I would like to thank all the team of the Mathematics Department
there for their great warmth at receiving me.

\section{Definitions and notations}
\label{defs}

\subsection{Basics}

Throughout all this text, $G$ will be a discrete hyperbolic group given
by a presentation $\presgroup {a_1,\ldots,a_m}{R}$ where
$S=\{a_1,\ldots,a_m,a_1\~,\ldots a_m\~\}$ is a symmetric set of $2m$
generators, and $R$ is a finite set of words on $S$. (Every discrete
hyperbolic group is finitely presented, cf.~\cite{S}.)

We shall denote by $\delta$ a hyperbolicity constant for $G$
w.r.t.\ $S$. Let $\lambda$ be the maximal length of relations in $R$.

A hyperbolic group is
called \emph{non-elementary} if it is neither finite nor quasi-isometric
to $\Z$.

A \emph{word} will be a word made of letters in $S$. Equality of words
will always mean equality as elements of the group $G$.

A word is said to be \emph{reduced} if it does not contain a generator
$a\in S$ immediately followed by its inverse $a\~$. It is said to be
\emph{cyclically reduced} if it and all of its cyclic permutations are
reduced.

If $w$ is a word, we shall call its number of letters its \emph{length}
and denote it by $\abs{w}$. Its \emph{norm}, denoted by $\norm{w}$, will
be the smallest length of a word equal to $w$ in the group $G$.

\subsection{Growth, cogrowth, and gross cogrowth}
\label{defexps}

First, we recall the definition of the growth, cogrowth and gross
cogrowth of the group $G$ with respect to the generating set $S$.

Let $S^\ell$ be the set of all words of length $\ell$ in $a_i^{\pm 1}$.
Let $S^\ell_G$ be the set of all elements of $G$ the norm of which is
equal to $\ell$ with respect to the generating set $a_i^{\pm 1}$.
The growth $g$ controls the asymptotics of the number of elements of
$S^\ell_G$: this number is roughly equal to $(2m)^{g\ell}$. The gross
cogrowth $\theta$ controls the asymptotics of the number of
words in $S^\ell$ which are equal to the neutral element in $G$: this
number is roughly equal to $(2m)^{\theta \ell}$. The cogrowth $\eta$ is the same
with reduced words only: this number is roughly $(2m-1)^{\eta\ell}$.

These quantities have been extensively studied. Growth now belongs to the
folklore of discrete group theory (see e.g.~\cite{GdlH} or~\cite{GK} for
background and open problems).  Cogrowth has been introduced by
R.~Grigorchuk in~\cite{Gri}, and independently by J.~Cohen in~\cite{C}.
For some examples see~\cite{Ch2} or~\cite{W1}. Gross cogrowth is linked
(see below) to the spectrum of the random walk on the group, which, since
the seminal work by H.~Kesten (see~\cite{K1} and~\cite{K2}), has been
extensively studied (see for example the numerous technical results
in~\cite{W2} and the references therein).

\begin{defi}[ (Growth, cogrowth, gross cogrowth)]

The \emph{growth} of the group $G$ with respect to the generating set
$a_1,\ldots,a_m$ is defined as
\[
g=\lim_{\ell \rightarrow\infty} \frac1\ell \log_{2m} \#S^\ell_G
\]

The \emph{gross cogrowth} of the group $G$ with respect to the generating set
$a_1,\ldots,a_m$ is defined as
\[
\theta=\lim_{
\begin{subarray}{c}
\ell \rightarrow\infty
\\ \ell\text{ even}
\end{subarray}
}
\frac1\ell \log_{2m} \#\{w\in S^\ell, w=e\text{ in }G\}
\]

The \emph{cogrowth} of the group $G$ with respect to the generating set
$a_1,\ldots,a_m$ is defined as $\eta=1/2$ for a free group, and otherwise
\[
\eta=\lim_{ 
\begin{subarray}{c}
\ell \rightarrow\infty
\\ \ell\text{ even}
\end{subarray}
}
\frac1\ell \log_{2m-1} \#\{w\in S^\ell, w=e\text{ in }G, w\text{ reduced}\}
\]
\end{defi}

Let us state some properties of these quantities. All of them are proven
in~\cite{K2}, \cite{Gri} or~\cite{C}.

The limits are well-defined by a simple subadditivity (for growth) or
superadditivity (for the cogrowths) argument. We restrict
ourselves to even $\ell$ because there may be no word of odd length equal
to the trivial element, as is the case e.g.\ in a free group.

For cogrowth, the logarithm is taken in base $2m-1$ because the number of
reduced words of length $\ell$ behaves like $(2m-1)^\ell$.

Cogrowth and gross cogrowth lie between $1/2$ and $1$. Gross cogrowth
is strictly above $1/2$, as well as cogrowth except for the free group.
There exist groups with cogrowth or gross cogrowth arbitrarily close to
$1/2$.

The probability that a random walk in the group $G$ (with respect to the
same set of generators) starting at $e$, comes back to $e$ at time $\ell$
is equal to the number of words equal to $e$ in $G$, divided by the total
number of words of length $\ell$. This leads to the following
characterization of gross cogrowth, which states that the return
probability at time $t$ is roughly equal to $(2m)^{-(1-\theta)t}$. This
will be ubiquitous in our text.

\begin{enonce}{Alternate definition of gross cogrowth}
Let $P_t$ be the probability that a random walk on the group $G$ (with
respect to the generating set $a_1,\ldots,a_m$) starting at
$e$ at time $0$, comes back to $e$ at time $t$.

Then the gross cogrowth of $G$ w.r.t.\ this generating set is equal to
\[
\theta=1+\lim_{
\begin{subarray}{c}t\rightarrow\infty\\t\text{ even}\end{subarray}
}
\frac1t \log_{2m} P_t
\]
\end{enonce}

In particular, $(2m)^{\theta-1}$ is the spectral radius of the random
walk operator (denoted $\lambda$ in~\cite{K1} and $r$ in~\cite{Gri}),
which is the form under which it is studied in these papers.

A cogrowth, or gross cogrowth, of $1$ is equivalent to amenability.

It is easy to check that $g/2+\theta\geq 1$.

Gross cogrowth and cogrowth are linked by the following equation
(see~\cite{Gri}):
\[
(2m)^\theta=(2m-1)^\eta+(2m-1)^{1-\eta}
\]

The gross cogrowth of the free group $F_m$ is
$\frac12\log_{2m}\left(8m-4\right)$, and this is the only group on $m$
generators with this gross cogrowth (see~\cite{K1}). This tends to $1/2$ as
$m\rightarrow\infty$.

There are various conventions for the cogrowth of the free group. The
definition above would give $-\infty$. In~\cite{C} the cogrowth of the
free group is taken equal to $0$; in~\cite{Gri} it is not defined. Our
convention allows the formula above between cogrowth and gross cogrowth
to be valid even for the free group; it is also natural given the fact
that, for any group except the free group, the cogrowth is strictly above
$1/2$. Moreover, this leads to a single formulation for our random
quotient theorem, as with this convention, the critical density for
quotients by reduced words will be equal to $1-\eta$ in any case. So we
strongly plead for this being the right convention.

If $G$ is presented as $F_m/N$ where $N$ is a normal subgroup,
cogrowth is the growth (in base $2m-1$) of $N$. Gross cogrowth is the
same considering $N$ as a submonoid in the free monoid on $2m$
generators and in base $2m$.

Let $\Delta$ be the Laplacian on $G$ (w.r.t.\ the same generating set).
As the operator of convolution by a random walk is equal to $1-\Delta$,
we get another characterization of gross cogrowth.
The eigenvalues lie in the interval $[0;2]$. Let $\lambda_0$ be the
smallest one and $\lambda'_0$ the largest one.
Then the gross cogrowth of $G$ w.r.t.\ this generating set is equal to
\[
\theta=1+ \log_{2m} \sup (1-\lambda_0, \lambda'_0-1)
\]
(We have to consider $\lambda'_0$ due to parity problems.)

Cogrowth and gross cogrowth depend on the generating set. For
example, adding trivial generators $a_i=e$ makes them arbitrarily close
to $1$.

\subsection{Diagrams}

A \emph{filamenteous van Kampen diagram} in the group $G$ with respect to the
presentation $\presgroup{S}{R}$ will be a planar connected combinatorial $2$-complex decorated in the following way:
\begin{itemize}
\item Each $2$-cell $c$ bears some relator $r\in R$. The number of edges of
the boundary of $c$ is equal to $\abs{r}$. 
\item If $e$ is an (unoriented) edge, denote by $e_+$ and $e_-$ its two
orientations. Then $e_+$ and $e_-$ both bear some generator $a\in S$, and
these two generators are inverse.
\item Each $2$-cell $c$ has a marked vertex on its boundary, and an
orientation at this vertex.
\item The word read by going through the (oriented) edges of the boundary
of cell $c$, starting at the marked point and in the direction given by
the orientation, is the relator $r\in R$ attached to $c$.
\end{itemize}

{\footnotesize
Note on the definition of regular complexes: we do not require that each
closed $2$-cell be homeomorphic to the
standard disc. We only require the
interior of the $2$-cell to be homeomorphic to a disc, that is, the
application may be non-injective on the boundary. This makes a difference
only when the relators are not reduced words. For example, if
$abb\~c$ is a relator, then the two diagrams below are valid.
We will talk about \emph{regular diagrams} to exclude the latter.

\begin{center}
\includegraphics{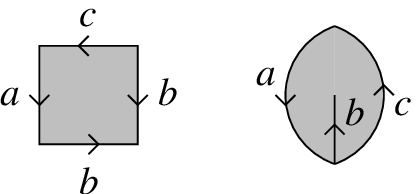}
\end{center}
}

We will use the terms \emph{$2$-cell} and \emph{face} interchangeably.

A \emph{non-filamenteous} van Kampen diagram will be a diagram in which every
$1$- or $0$-cell lies in the boundary of some $2$-cell. Unless otherwise
stated, in our text \emph{a van Kampen diagram will implicitly be
non-filamenteous}.

A \emph{$n$-hole} van Kampen diagram will be one for which the underlying
$2$-complex has $n$ holes. When the number of holes is not given, \emph{a
van Kampen diagram will be supposed to be simply connected} ($0$-hole).

\label{firstdavKd}
A \emph{decorated abstract van Kampen diagram} (davKd for short) is defined
almost the same way as a van Kampen diagram, except that no relators are
attached to the $2$-cells and no generators attached to the edges, but
instead, to each $2$-cell is attached an integer between $1$ and the
number of $2$-cells of the diagram (and yet, a starting point and
orientation to each $2$-cell).

{\footnotesize
Please note that this definition is a little bit emended in
section~\ref{newdavKds} (more decoration is added).
}

A davKd is said to be \emph{fulfillable} w.r.t.\ presentation
$\presgroup{S}{R}$ if there
exists an assignment of relators to $2$-cells and of generators to
$1$-cells, such that any two $2$-cells bearing the same number get the same
relator, and such that the resulting decorated diagram is a van Kampen diagram
with respect to presentation $\presgroup{S}{R}$.

A \emph{davKd with border $w_1,\ldots,w_n$}, where $w_1, \ldots,w_n$ are
words, will be a $(n-1)$-hole davKd with each boundary edge decorated by a
letter such that the words read on the $n$ components of the boundary are
$w_1,\ldots,w_n$. A davKd with border is said to be \emph{fulfillable} if, as a
davKd, it is fulfillable while keeping the same boundary words.

A word $w$ is equal to the neutral element $e$ in $G$ if and only if some
no-hole, maybe filamenteous, davKd with border $w$ is fulfillable
(see~\cite{LS}).

A van Kampen diagram is said to be \emph{reduced} if there is no pair of
adjacent (by an edge) $2$-cells bearing the same relator with opposite
orientations and with the common edge representing the same letter in the
relator (w.r.t.\ the starting point).  A davKd is said to be \emph{reduced}
if there is no pair of adjacent (by an edge) $2$-cells bearing the same
number, with opposite orientations and a common edge representing the
same letter in the relator.

A van Kampen diagram is said to be \emph{minimal} if it has the minimal number
of $2$-cells among those van Kampen diagrams having the same boundary word (or
boundary words if it is not simply connected). A fulfillable davKd with
border is said to be \emph{minimal} in the same circumstances.

Note that a minimal van Kampen diagram is necessarily reduced: if
there were a pair of adjacent faces with the same relator in opposite
orientations, then they could be removed to obtain a new diagram with less
faces and the same boundary (maybe adding some filaments):

\begin{center}
\includegraphics{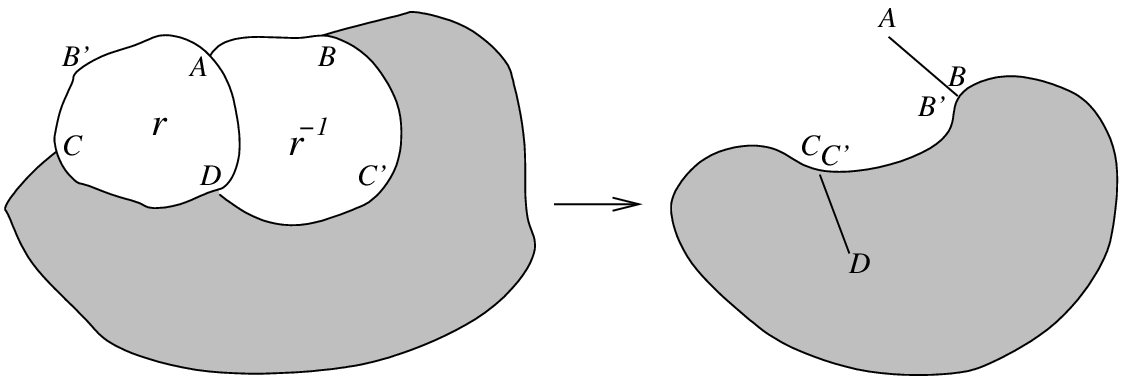}
\end{center}

\medskip

Throughout the text, we shall use the term \emph{diagram} as a short-hand
for ``van Kampen diagram or fulfillable decorated abstract van Kampen diagram''. We
will use the term \emph{minimal diagram} as a short-hand for ``minimal
van Kampen diagram or minimal fulfillable decorated abstract van Kampen diagram with
border''.

\subsection{Isoperimetry and narrowness}

There is a canonical metric on the $1$-skeleton of a van Kampen diagram (or a
davKd), which assigns length $1$ to every edge. If $D$ is a diagram, we
will denote its number of faces by $\abs{D}$ and the length of its
boundary by $\abs{\d D}$.

It is well-known (see~\cite{S}) that a discrete group is hyperbolic if and
only if there exists a constant $C>0$ such that any minimal diagram $D$
satisfies the linear isoperimetric inequality $\abs{\d D}\geq C \abs{D}$.
We show in Appendix~\ref{appiso} that in a hyperbolic group, holed diagrams
satisfy an isoperimetric inequality as well.

Throughout all the text, $C$ will be an isoperimetric constant for
$G$.

The set of $2$-cells of a diagram is also canonically equipped with a
metric: two $2$-cells sharing a common edge are defined to be at
distance $1$. The \emph{distance to the boundary} of a face will be its
distance to the exterior of the diagram considered as a face, i.e.\ a
boundary face is at distance $1$ from the boundary.

A diagram is said to be \emph{$A$-narrow} if any $2$-cell is at distance
at most $A$ from the boundary.

It is well-known, and we show in Appendix~\ref{appiso} in the form we
need, that a linear isoperimetry implies narrowness of minimal diagrams.

\section{The standard case: $F_m$}
\label{standardcase}

We proceed here to the proof of Gromov's now classical theorem
(Theorem~\ref{G12}) that a
random quotient of the free group $F_m$ is trivial in density greater
than $1/2$, and non-elementary hyperbolic in density smaller than this
value.

We include this proof here because, first, it can serve as a useful
template for understanding the general case, and,
second, it seems that no completely correct proof has been published so
far\footnote{Since the proof included here was written and diffused, a
similar but somewhat simpler proof has been published in~\cite{Z} for a
slightly different model in which relators are of length $3$ but the number
of generators $m$ tends to infinity.}.

Recall that in this case, we consider a random quotient of the free group
$F_m$ on $m$ generators by $(2m-1)^{d\ell}$ uniformly chosen
\emph{cyclically reduced} words of length $\ell$.

A random cyclically reduced word is chosen in the following way: first
choose the first letter ($2m$ possibilities), then choose the next letter
in such a way that it is not equal to the inverse of the preceding one ($2m-1$
possibilities), up to the last letter which has to be distinct both from
the penultimate letter and the first one (which lets $2m-2$ or $2m-1$
choices depending on whether the penultimate letter is the same as the
first one). The difference between $2m$ and $2m-1$ at the first position,
and between $2m-1$ and $2m-2$ at the last position is negligible (as
$\ell\rightarrow \infty$) and we will do as if we had $2m-1$ choices for
each letter exactly.

So, for the sake of simplicity of the exposition, in the following we may
assume that there are exactly $(2m-1)^\ell$ reduced words of length $\ell$, with
$2m-1$ choices for each letter. Bringing the argument to full correctness
is a straightforward exercise.

\subsection{Triviality for $d>1/2$}

The triviality of the quotient for $d>1/2$ reduces essentially to the
well-known

\begin{enonce}{Probabilistic pigeon-hole principle}
Let $\eps>0$ and put $N^{1/2+\eps}$ pigeons uniformly at random among $N$
pigeon-holes. Then there are two pigeons in the same hole with
probability tending to $1$ as $N\rightarrow\infty$ (and this happens
arbitrarily many times with growing $N$).
\end{enonce}

Now, take as your pigeon-hole the word made of the first $\ell-1$ letters
of a random word of length $\ell$. There are $(2m-1)^{\ell-1}$ pigeon-holes
and we pick up $(2m-1)^{d\ell}$ random words with $d>1/2$. Thus, with
probability arbitrarily close to $1$ with growing $\ell$, we will pick
two words of the form $wa_i$, $wa_j$ where $\abs{w}=\ell-1$ and
$a_i,a_j\in S$. Hence in the quotient group we will have $a_i=a_j$.

But as $d$ is strictly bigger than $1/2$, this will not occur only once
but arbitrarily many times as $\ell\rightarrow \infty$, with at each time
$a_i$ and $a_j$ being chosen at random from $S$. That is, for big enough
$\ell$, all couples of generators $a,b\in S$ will satisfy $a=b$ in the
quotient group. As $S$ is symmetric, in particular they will satisfy $a=a\~$.

The group presented by $\presgroup{(a_i)}{a_i=a_i\~, a_i=a_j\ \forall
i,j}$ is $\Z/2\Z$. In
case $\ell$ is even this is exactly the group we get (as there are only
relations of even length), and if $\ell$ is odd any relation of odd
length turns $\Z/2\Z$ into $\{e\}$.

This proves the second part of Theorem~\ref{G12}.

\subsection{Hyperbolicity for $d<1/2$}

We proceed as follows: We will show that the (reduced) davKd's which are
fulfillable by a random presentation necessarily satisfy some linear
isoperimetric inequality. This is stronger than proving that only minimal
diagrams satisfy an isoperimetric inequality: in fact, \emph{all} reduced
diagrams in a random group satisfy this inequality. (Of course this
cannot be true of non-reduced diagrams since one can, for example, take
any relator $r$ and arrange an arbitrarily large diagram of alternating
$r$'s and $r\~$'s like in a chessboard.)

Thus we will evaluate the probability that a given
decorated abstract van Kampen diagram can be fulfilled by a random
presentation. We show that if the davKd violates the isoperimetric
inequality, then this probability is
very small and in fact decreases exponentially with $\ell$.

Then, we apply the Cartan-Hadamard-Gromov theorem for
hyperbolic spaces, which tells us that to ensure hyperbolicity of a
group, it is not necessary to check the isoperimetric inequality for
\emph{all} diagrams but for a \emph{finite number} of them (see
section~\ref{CHGP} for details).

Say is it enough to check all diagrams with at most $K$ faces, where $K$
is some constant depending on $d$ but not on $\ell$. Assume we know that for
each of these diagrams which violates the isoperimetric inequality, the
probability that it is fulfillable decreases exponentially with $\ell$.
Let $D(K)$ be the (finite) number of davKd's with at most $K$ faces,
violating the isoperimetric inequality. The probability that at least one
of them is fulfillable is less that $D(K)$ times some quantity decreasing
exponentially with $\ell$, and taking $\ell$ large enough ensures that
with probability arbitrarily close to one, none of these davKd's is
fulfillable. The conclusion then follows by the
Cartan-Hadamard-Gromov theorem.

\bigskip

The intuitive basic picture is as follows: Consider a davKd made of two
faces of perimeter $\ell$ meeting along $L$ edges. The probability that
two given random relators $r$, $r'$ fulfill this diagram is at most
$(2m-1)^{-L}$, which is the probability that $L$ given letters of $r$ are
the inverses of $L$ given letters of $r'$. (Remember that as the relators
are taken reduced, there are only $2m-1$ choices for each letter except
for the first one. As $2m-1<2m$ we can safely treat the first letter like
the others, as doing otherwise would still sharpen our evaluation.)

\begin{center}
\includegraphics{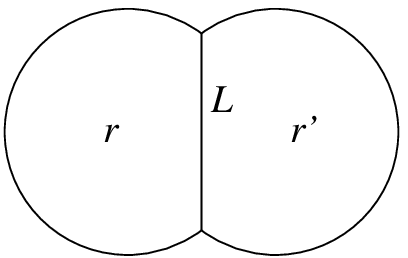}
\end{center}

Now, there are $(2m-1)^{d\ell}$ relators in the presentation. As we said,
the probability that two given relators fulfill the diagram is at most
$(2m-1)^{-L}$. Thus, the probability that there exist two relators in the
presentation fulfilling the diagram is at most
$(2m-1)^{2d\ell}\,(2m-1)^{-L}$, with the new factor accounting for the
choice of the two relators.

This evaluation becomes non-trivial for $L>2d\ell$. Observe that the
boundary length of the diagram is $2\ell-2L=2(1-2d)\ell-2(L-2d\ell)$.
That is, if $L\leq 2d\ell$ then the boundary is longer than
$2(1-2d)\ell$, and if $L>2d\ell$ then the probability that the diagram
can be fulfilled is exponentially small with $\ell$.

To go on with our intuitive reasoning, consider a graph with $n$ relators
instead of two. The number of ``conditions'' imposed by the graph is
equal to the total length $L$ of its internal edges, that is, the
probability that a random assignment of relators satisfy them is
$(2m-1)^{-L}$, whereas the number of choices for the relators is
$(2m-1)^{nd\ell}$ by definition. So if $L>nd\ell$ the probability is too
small. But if $L\leq nd\ell$, then the boundary length, which is equal to
$n\ell-2L$, is bigger than $(1-2d)n\ell$ which is the isoperimetric
inequality we were looking for.

This is the picture we will elaborate on. In fact, what was false in the
last paragraph is that if the same relator is to appear several times in the
diagram, then we cannot simply multiply probabilities as we did, as these
probabilities are no more independent.

\bigskip

Thus, let $D$ be a reduced davKd. We will evaluate the probability that it
can be fulfilled by relators of a random presentation. Namely

\begin{prop}
Let $D$ be a reduced davKd. The probability that $D$ can be fulfilled by
relators of a random presentation is at most $(2m-1)^{(\abs{\d
D}-\ell\abs{D}(1-2d))/2\abs{D}}$.
\end{prop}

\begin{dem}
Each face of $D$ bears a number between $1$ and $\abs{D}$. Let $n$ be the
number of distinct numbers the faces bear in $D$. Of course, $n\leq
\abs{D}$. (The original proof by Gromov was valid only when $n=\abs{D}$, so
that all relators are chosen independently, which simplifies the proof.
If $n<\abs{D}$ then we cannot simply multiply probabilities as in the basic
picture.) Suppose, for simplicity, that these $n$ distinct numbers are
$1,2,\ldots,n$.

To fulfill $D$ is to give $n$ relators $r_1,\ldots, r_n$ satisfying the
relations imposed by the diagram.

We will construct an auxiliary graph $\Gamma$ summarizing all letter
relations imposed by the diagram $D$. Vertices of $\Gamma$ will represent
the letters of $r_1,\ldots,r_n$, and edges of $\Gamma$ will represent
inverseness (or equality, depending on orientation) of letters imposed by
shared edges between faces of $D$.

Thus, take $n\ell$ vertices for $\Gamma$, arranged in $n$ parts of
$\ell$ vertices. Call the vertices corresponding to the faces of $D$ bearing
number $i$ the $i$-th part of the graph. Each part is made of $\ell$
vertices.

We now explain what to take as edges of $\Gamma$.

In the diagram, every face is marked with a point on its boundary, and an
orientation. Label the edges of each face $1,2,\ldots,\ell$ starting at
the marked point, following the given orientation.

If, in the davKd $D$, the $k$-th edge of a face bearing number $i$ is
equal to the $k'$-th edge of an adjacent face bearing number $j$, then
put an edge in $\Gamma$ between the $k$-th vertex of the $i$-th part and
the $k'$-th vertex of the $j$-th part. Decorate the newly
added edge with $-1$ if the two faces' orientations agree, or with $+1$
if they disagree.

Thus, a $-1$ edge between the $k$-th vertex of the $i$-th part and the
$k'$-th vertex of the $j$-th part means that the $k$-th letter of
relator $r_i$ has to be the inverse of the $k'$-th letter of relator
$r_j$.

Successively add an edge to $\Gamma$ in this way for each interior edge
of the davKd $D$, so that the total number of edges of $\Gamma$ is equal
to the number of interior edges of $D$.

As $D$ is reduced, the graph $\Gamma$ can contain no loop. It may well
have multiple edges, if, in the davKd, several pairs of adjacent faces
bear the same numbers and have common edges at the same position.

Note that this graph only depends on the davKd $D$ and in no way on the
random presentation.

The graph $\Gamma$ for the basic picture above is:

\begin{center}
\includegraphics{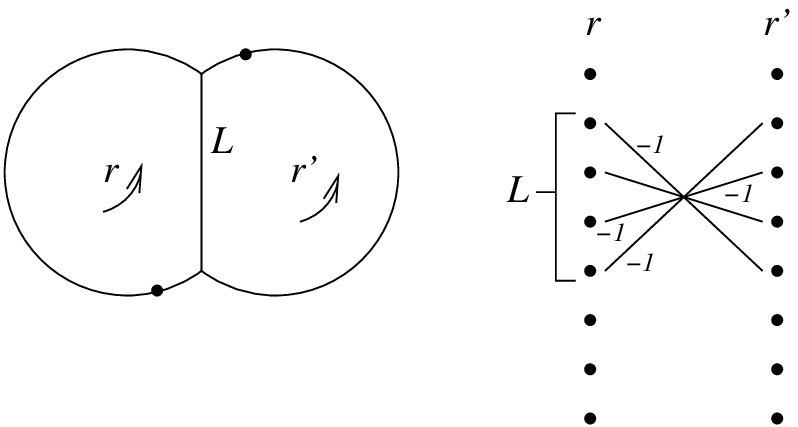}
\end{center}

\bigskip

Now let us evaluate the probability that $D$ is fulfillable. To fulfill
$D$ is to assign a generator to each vertex of $\Gamma$ and see if the
relations imposed by the edges are satisfied.

Remark that if the generator of any vertex of the graph is assigned, then this
fixes the generators of its whole connected component. (And, maybe,
depending on the signs of the edges of $\Gamma$, there is no correct
assignation at all.) Thus, the number of degrees of freedom is at most
equal to the number of connected components of $\Gamma$.

Thus (up to our approximation on the number of cyclically reduced words),
the number of random assignments of cyclically reduced words to the
vertices of $\Gamma$ is $(2m-1)^{n\ell}$, whereas the number of those
assignments satisfying the constraints of the edges is at most $(2m-1)^C$
where $C$ is the number of connected components. Hence, the probability
that a given assignment of $n$ random words to the vertices of $\Gamma$
satifies the edges relations is at most $(2m-1)^{C-n\ell}$.

This is the probability that $n$ \emph{given} relators of a random
presentation fulfill the diagram. Now there are $(2m-1)^{d\ell}$ relators
in a random presentation, so the probability that we can find $n$ of them
fulfilling the diagram is at most $(2m-1)^{nd\ell}\,(2m-1)^{C-n\ell}$.

Now let $\Gamma_i$ be the subgraph of $\Gamma$ made of those vertices
corresponding to a face of $D$ bearing a number $\leq i$. Thus
$\Gamma_1\subset\Gamma_2\subset\ldots\subset \Gamma_n=\Gamma$. Of course,
the probability that $\Gamma$ is fulfillable is less than any of the
probabilities that $\Gamma_i$ is fulfillable for $i\leq n$.

The above argument on the number of connected components can be repeated
for $\Gamma_i$: the probability that $\Gamma_i$ is fulfillable is at most
$(2m-1)^{id\ell+C_i-i\ell}$ where $C_i$ is the number of connected
components of $\Gamma_i$.

This leads to setting 
\[
d_i=id\ell+C_i-i\ell
\]
and following Gromov we
interpret this number as the dimension of $\Gamma_i$, or, better, the
dimension of the set of random presentations for which there exist $i$
relators satisfying the conditions imposed by $\Gamma_i$. Thus:
\[
\Pr(D\text{ is fulfillable})\leq (2m-1)^{d_i}\qquad \forall i
\]

Before concluding we need a further purely combinatorial lemma.

\begin{lem}
\[
\abs{\d D} \geq 
\ell \abs{D} (1-2d) + 2 \sum d_i(m_i-m_{i+1})
\]
where $m_i, 1\leq i\leq n$ is the number of faces of $D$ bearing relator number $i$.
\end{lem}

Before proving the lemma let us end the proof of the proposition.
We are free to choose the order of the construction, and we may suppose that
the $m_i$'s are non-increasing, i.e.\ that we began with the relator
appearing the biggest number of times in $D$, etc., so that $m_i-m_{i+1}$
is non-negative.

If all $d_i$'s are non-negative, then we have the isoperimetric
inequality $\abs{\d D}\geq \ell \abs{D} (1-2d)$ and the proposition is
true since the probability at play is at most $1$.

If some $d_i$ is negative, we use the fact established above that the
probability that the diagram is fulfillable is less than $(2m-1)^{\inf
d_i}$. As $\sum m_i=\abs{D}$, we have $\sum d_i(m_i-m_{i+1})\geq
\abs{D}\inf d_i$. Thus 
$\inf d_i\leq \left(\abs{\d D}-\ell\abs{D}(1-2d)\right)/2\!\abs{D}$ hence
the proposition.
\end{dem}

\begin{dem}[ of the lemma]
A vertex in the $i$-th part of $\Gamma$ is thus of
multiplicity at most $m_i$. Let $A$ be the number of edges in $\Gamma$.
We have
\[
\abs{\d D}\geq\abs{D}\ell-2A=\ell\sum m_i-2A
\]
(where the equality $\abs{\d D}= \abs{D}\ell-2A$ holds when $D$ has no
filaments).

Thus we want to show that either the number of edges is small, or the
fulfillability probability is small. The latter grows with the number of
connected components of $\Gamma$, so this looks reasonable.

Let $A_i$ be the number of edges in $\Gamma_i$. We now show that
\[
A_{i+1}-A_i+m_{i+1}(d_{i+1}-d_i)\leq m_{i+1}d\ell
\]
or equivalently that
\label{edgedim}
\[
A_{i+1}-A_i+m_{i+1}(C_{i+1}-(C_i+\ell))\leq 0
\]

Depart from $\Gamma_i$ and add the new vertices and edges of
$\Gamma_{i+1}$. When adding the $\ell$ vertices, the number of connected
components increases by $\ell$. So we only have to show that when adding the
edges, the number of connected components decreases at least by $1/m_{i+1}$
times the number of edges added.

Call \emph{external point} a point of $\Gamma_{i+1}\setminus \Gamma_i$
which shares an edge with a point of $\Gamma_i$. Call \emph{internal
point} a point of $\Gamma_{i+1}\setminus \Gamma_i$ which is not external.
Call \emph{external edge} an edge between an external point and a point
of $\Gamma_i$, \emph{internal edge} an edge between two internal points,
and \emph{external-internal edge} an edge between an external and
internal point. Call \emph{true internal point} a point which has at
least one internal edge.

While adding the external edges, each external point is connected to a
connected component inside $\Gamma_i$, and thus the number of connected
components decreases by $1$ for each external point.

Now add the internal edges (but not yet the external-internal ones): If
there are $N$ true internal points, these make at most $N/2$
connected components after adding the internal edges, so the number of
connected components has decreased by at least $N/2$.

After adding the external-internal edges the number of connected
components still decreases. Thus it has decreased by at least the number
of external points plus half the number of true internal points.

Now as each external point is of degree at most $m_{i+1}$, the number of
external plus external-internal edges is at most $m_{i+1}$ times the
number of external points. If there are $N$ true internal points, the
number of internal edges is at most $Nm_{i+1}/2$ (each edge is counted $2$
times). So the total number of edges is at most $m_{i+1}$ times the
number of external points plus half the number of true internal points,
which had to be shown.

Thus we have proved that $A_{i+1}-A_i+m_{i+1}(d_{i+1}-d_i)\leq
m_{i+1}d\ell$.
Summing over $i$ yields \[A+\sum m_i(d_i-d_{i-1})\leq d\ell \sum m_i\]

Thus,
\begin{eqnarray*}
\abs{\d D}&\geq&\ell\sum m_i-2A
\\&\geq &
\ell\sum m_i-2d\ell \sum m_i +2\sum m_i(d_i-d_{i-1})
\\ &= &
\ell \abs{D} (1-2d) + 2 \sum d_i(m_i-m_{i+1})
\end{eqnarray*}
as was needed.
\end{dem}

\begin{cor}
Let $D$ be a davKD. Then, either $D$ satisfies the isoperimetric
inequality
\[
\abs{\d D}\geq \ell\abs{D} (1/2-d)
\]
or the probability that is can be fulfilled by relators of a random
presentation is at most $(2m-1)^{-\ell(1/2-d)/2}$.
\end{cor}

Hence the interest of taking $d<1/2$...

This was for a given davKd $D$. In order to show that the group is
hyperbolic, we have to show that the probability that there exists a
davKd violating the isoperimetric inequality tends to $0$ when
$\ell\rightarrow \infty$. But here we use the local-global principle for
hyperbolic grometry (or Cartan-Hadamard-Gromov theorem, see
Appendix~\ref{CHGP}), which can be stated as:

\begin{prop*}
For each $\alpha>0$, there exist an integer $K(\alpha)\geq 1$ and an
$\alpha'>0$ such that,
if a group is given by relations of length $\ell$ for some $\ell$ and if
any reduced van Kampen diagram with at most $K$ faces satifies
\[
\abs{\d D}\geq \alpha\ell \abs{D}
\]
then any reduced van Kampen diagram $D$ satisfies
\[
\abs{\d D}\geq \alpha' \ell \abs{D}
\]
(hence the group is hyperbolic).
\end{prop*}

Now take $\alpha=1/2-d$ and the $K$ given by the proposition. If
$N(K,\ell)$ is the number of davKd's with at most $K$ faces and each face
has $\ell$ edges, then the probability that one of them is fulfillable
and violates the isoperimetric inequality is at most
$N(K,\ell)\,(2m-1)^{-\ell(1/2-d)/2}$.

\begin{prop}
For fixed $K$, the number $N(K,\ell)$ grows polynomially with $\ell$.
Hence, the probability $N(K,\ell)\,(2m-1)^{-\ell(1/2-d)/2}$ tends
exponentially to $0$ as $\ell\rightarrow \infty$.
\end{prop}

\begin{dem}
Let us evaluate $N(K,\ell)$. As the relators in the presentation are
taken to be cyclically reduced, we only have to consider regular diagrams
(see section~\ref{defs}). A regular davKd is only a planar graph with some
decoration on the edges, namely, a planar graph with on each edge a
length indicating the number of edges of the davKd it represents, and with
vertices of degree at least $3$ (and, as in a davKd, every face is
decorated with a starting point, an orientation, and a number between $1$
and $K$). Let $G(K)$ be the number of planar graphs with vertex degree at
least $3$. In such a graph there are (by Euler's formula) at most $3K$
edges, so there are at most $\ell^{3K}$ choices of edge lengths, and we
have $(2\ell K)^K$ choices for the decoration of each face (orientation,
starting point and number between $1$ and $K$).

So $N(K,\ell)\leq G(K)(2K)^K \ell^{4K}$.
\end{dem}

This proves that the quotient is hyperbolic; we now show that it is
infinite. We can of course use the general argument of
section~\ref{infiniteness} but there is a shorter proof in this case.
First, as any reduced diagram
satisfies $\abs{\d D}\geq \alpha' \ell \abs{D}\geq \alpha'\ell$, the
ball of radius $\alpha'\ell/2$ injects into the quotient, hence the
quotient contains at least one non-trivial element and cannot be $\{e\}$.

Second, we prove that the presentation is aspherical. With our
conventions on van Kampen diagrams, our asphericity implies asphericity
of the Cayley complex and thus cohomological dimension at most $2$
(indeed, thanks to the marking of each face by a starting point and a
relator number, two faces are reducible in a diagram only if they really
are the same face in the Cayley complex, so that diagram reduction is a
homotopy in the Cayley complex).  This will end the proof:
indeed, cohomological dimension at most $2$ implies torsion-freeness
(see~\cite{B}, p.~187), hence the quotient cannot be a non-trivial finite
group.

Indeed, the isoperimetric inequality above is not only valid for minimal
diagrams, but for \emph{any} reduced diagram. Now suppose that there is
some reduced spherical diagram. It will have zero boundary length and
thus will violate any isoperimetric inequality, hence a contradiction.
Thus the presentation is aspherical.

This proves Theorem~\ref{G12}.

\section{Outline of the argument}
\label{outline}

Here we explain some of the ideas of the proof of
Theorems~\ref{reducedmain}, \ref{geodmain} and~\ref{main}.

We will give a general theorem for hyperbolicity of random quotients by
words taken from some probability measures on the set of all words. We
will need somewhat technical axioms on the measures (for example, that
they weight only long words). Here we give a heuristic justification of
why these axioms are needed.

We proceed by showing that van Kampen diagrams of the quotient $G/\langle R\rangle$
satisfy a linear isoperimetric inequality.

If $D$ is a van Kampen diagram of the quotient, let $D'$ be the subcomplex of
$D$ made of relators of the presentation of $G$ (``old relators'') and
$D''$ the subcomplex made of relators in $R$ (``new relators'').

Say the new relators have length of order $\ell$ where $\ell$ is much
bigger than the hyperbolicity constant of $G$. (This will be Axiom~1.)

The main point will be that $D'$ is a diagram in the hyperbolic group
$G$, and, as such, is narrow (see Appendix~\ref{appiso}). We show below
that its narrowness is of order $\log\ell$. Hence, if $\ell$ is big
enough, the diagram $D$ can be viewed as big faces representing the new
relators, separated by a thin layer of ``glue'' representing the old
relators. The ``glue'' itself may contain invaginations in the new
relators and narrow excrescences on the boundary.

\begin{center}
\includegraphics{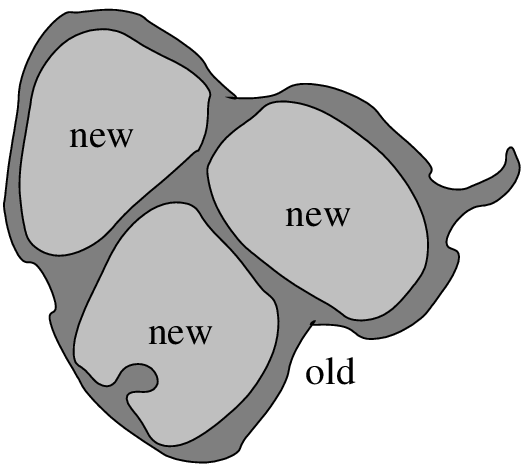}
\end{center}

\subsection{A basic picture}

As an example, let us study a basic picture consisting of two new
relators separated by some old stuff. Say that two random new relators
$r, r'$ are ``glued'' along subwords of length $L, L'$ (we may have
$L\neq L'$). Let $w$ be the word bordering the part of the diagram made
of old relators, we have $\abs{w}=L+L'+o(\ell)$. By construction, $w$ is
a word representing the trivial element in $G$. Write $w=xux'v$ where $x$
is a subword of $r$ of length $L$, $x'$ is a subword of $r'$ of length
$L'$, and $u$ and $v$ are short words.

\begin{center}
\includegraphics{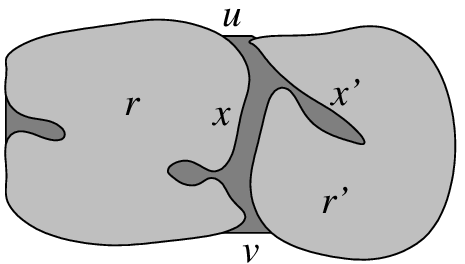}
\end{center}

Let us evaluate the probability that such a diagram exists. Take two
given random relators $r, r'$ in $R$. The probability that they can be
glued along subwords $x,x'$ of lengths $L, L'$ by narrow glue in $G$ is
the probability that there exist short words $u,v$ such that $xux'v=e$ in
$G$.

If, as in the standard case, there were no glue (no old relators) and $r$
and $r'$ were uniformly chosen random reduced words, the probability that
$r$ and $r'$ could be glued along subwords $x$, $x'$ of length $L$ (we
would have $L=L'$ in this case) would be $(2m-1)^{-L}$. But we now have
to consider the case when then $x$ and $x'$ are equal, not as words, but
as elements of $G$ (and up to small words $u$ and $v$, which we will
neglect).

If, for example, the relators are uniformly chosen random words, then $x$
and $x'$ are independent subwords, and the probability that $x$ and $x'$
are (almost) equal in $G$ is the probability that $x{x'}\~=e$; but
$x{x'}\~$ is a uniformly chosen random word of length $L+L'$, and by
definition the probability that it is equal to $e$ is controlled by the
gross cogrowth of $G$: this is roughly $(2m)^{-(1-\theta)(L+L')}$ (recall
the alternate definition of gross cogrowth in section~\ref{defexps}).

In order to deal not only with uniformly chosen random words but with
other situations such as random geodesic words, we will need a control on
the probability that two relators can be glued (modulo $G$) along
subwords of length $L$ and $L'$. This will be our Axiom~3: we will ask
this probability to decrease like $(2m)^{-\beta(L+L')}$ for some exponent
$\beta$ (equal to $1-\theta$ for plain random words).

Now in the simple situation with two relators depicted above, the length
of the boundary of the diagram is not exactly $2\ell-L-L'$, since there
can be invaginations of the relators, i.e.\ long part of the relators
which are equal to short elements in $G$ (as in the left part of the
picture above). In the case of uniformly chosen random relators, by
definition the probability that a part of length $L$ of a relator is
(nearly) equal to $e$ in $G$ is roughly $(2m)^{-(1-\theta)L}$. So,
again inspired by this case, we will ask for an axiom controlling the
length of subwords of our relators. This will be our Axiom~2.

Axiom~4 will deal with the special case when $r={r'}\~$, so that the words $x$
and $x'$ above are equal, and not at all chosen independently as we
implicitly assumed above. In this case, the size of centralizers of
torsion elements in the group will matter.

\bigskip

This was for given $r$ and $r'$. But there are $(2m)^{d\ell}$ relators in
$R$, so we have $(2m)^{2d\ell}$ choices for $r,r'$. Thus, the probability
that in $R$, there are two new relators that glue along subwords of
length $L, L'$ is less than $(2m)^{2d\ell}(2m)^{-\beta(L+L')}$.

Now, just observe that the length of the boundary of the diagram is (up
to the small words $u$ and $v$) $2\ell-L-L'$. On the other hand, when
$d<\beta$, the exponent $2d\ell-\beta(L+L')$ of the above probability
will be negative as soon as $L+L'$ is bigger than $2\ell$. This is
exactly what we want to prove: either the boundary is big, or the
probability of existence of the diagram is small.

This is comparable to the former situation with random quotients of the
free group: in the free group, imposing two random relators to glue along
subwords of lengths $L$ and $L'=L$ results in $L$ ``equations'' on the
letters. Similarly, in the case of plain random words, in a group of
gross cogrowth $\theta$, imposing two random words to glue along subwords
of lengths $L, L'$ results in $\beta(L+L')$ ``equations'' on these random
words, with $\beta=1-\theta$.

Now for diagrams having more than two new relators, essentially the
number of ``equations'' imposed by the gluings is $\beta$ times the total
internal length of the relators. The boundary is the external length. If
there are $n$ new relators and the total internal length is $A$, then the
boundary is roughly $n\ell-A$. But the probability of existence of such a
diagram is $(2m)^{-\beta A}(2m)^{nd\ell}$ where the last factor accounts
for the choice of the $n$ relators among the $(2m)^{d\ell}$ relators of
$R$. So if $d<\beta$, as soon as $A>n\ell d/\beta$, the probability decreases
exponentially with $\ell$; otherwise, the boundary is longer than
$n\ell(1-d/\beta)$.

\subsection{Foretaste of the Axioms}

As suggested by the above basic picture, we will demand four axioms: one
saying that our random relators are of length roughly $\ell$, another
saying that subwords of our relators are not too short, another one
controlling the probability that two relators glue along long subwords
(that is, the probability that these subwords are nearly equal in $G$),
and a last one controlling the probability that a relator glues along
its own inverse.

As all our estimates are asymptotic in the length of the words
considered, we will be allowed to apply them only to sufficiently long
subwords of our relators (and not to one individual letter, for example),
that is, to words of length at least $\eps\ell$ for some $\eps$.

Note that in order to be allowed to apply these axioms to any subword of
the relators at play, whatever happens elsewhere, we will need to ask
that different subwords of our relators behave quite independently from
each other; in our axioms this will result in demanding that the
probability estimates hold for a subword of a relator conditionnally to
whatever the rest of the relator is.

This is a strong independence condition, but, surprisingly enough, is it
valid not only for uniformly chosen random words (where by definition
everything is independent, in any group), but also for randomly chosen
geodesic words. This is a specific property of hyperbolic groups.

\bigskip

Several exponents will appear in the axioms. As we saw in the basic
picture, the maximal density up to which the quotient is non-trivial is
exactly the minimum of these exponents. Back to the intuition behind the
density model of a random quotient (see the introduction), the exponents
in our axioms indicate how many equations it takes in $G$ to have
certain gluings in our relators, whereas the density of the random
quotient is a measure of how many equations we can reasonably impose so
that it is still possible to find a relator satisfying them among our
randomly chosen relators. So this intuition gets a very precise numerical
meaning.

\section{Axioms on random words implying hyperbolicity of a random
quotient, and statement of the main theorem}
\label{axioms}

We want to study random quotients of a (non-elementary) hyperbolic group
$G$ by randomly chosen elements. Let $\mu_\ell$ be the law, indexed by
some parameter $\ell$ to tend to infinity, of the random elements
considered.

We will always assume that $\mu_\ell$ is a symmetric measure, i.e.\ for
any $x\in G$, we have $\mu_\ell(x)=\mu_\ell(x\~)$.

We will show that if the measure satisfies some simple axioms, then the
random quotient by elements picked under the measure is hyperbolic.

For each of the elements of $G$ weighted bu $\mu_\ell$, fix once and for
all a representation of it as a word (and choose inverse words for
inverse elements), so that $\mu_\ell$ can be considered as a measure on
words.  Satisfaction of our axioms may depend on such a choice.

Let $\mu_\ell^L$ be the law $\mu_\ell$ restricted (and rescaled) to words
of length $L$ (or $0$ if there are no such words in the support of
$\mu$). In most applications, $\mu_\ell$ will weight only words of length
$\ell$, but we will occasionally use laws $\mu_\ell$ weighting words of
length comprised between, say, $A\ell$ and $B\ell$.

To pick a random set $R$ of density at most $d$ is to pick, for each
length $L$, independently, at most $(2m)^{dL}$ random words of length $L$
according to law $\mu_\ell^L$.  That is, for each length, the density is
at most $d$.

(We say ``at most'' because we do not require that exactly $(2m)^{dL}$
words of length $L$ are taken for each $L$. Taking smaller $R$ will result
in a hyperbolic quotient as well.)

We want to show that if $d$ is less than some quantity depending on
$\mu_\ell$ (and $G$, since $\mu_\ell$ takes value in $G$), then the
random quotient $G/\langle R \rangle$ is very probably non-elementary
hyperbolic.

\subsection{Asymptotic notations}

By the notation $f(\ell)\approx g(\ell)$ we shall mean that
\[
\lim_{
\ell\rightarrow \infty} \frac1\ell \log f(\ell)
=
\lim_{
\ell\rightarrow \infty} \frac1\ell \log g(\ell)
\]

We define the notation $f(\ell)\lesssim g(\ell)$
similarly. We will say, respectively, that $f$
is roughly equal or roughly less than $g$.

Accordingly, we will say that $f(\ell,L)\approx g(\ell,L)$ uniformly for
all $L\leq \ell$ if whatever the sequence $L(\ell)\leq \ell$ is, we have
\[
\lim_{
\ell\rightarrow \infty} \frac1\ell \log f(\ell,L(\ell)) = \lim_{
\ell\rightarrow \infty} \frac1\ell \log g(\ell,L(\ell))
\]and if this
limit is uniform in the sequence $L(\ell)$.

\subsection{Some vocabulary}
\label{somevocabulary}

Here we give technical definitions designed in such a manner that
the axioms can be stated in a natural way. We recommend to look at the
axioms first.

Let $x$ be a word. For each $a,b$ in $[0;1]$ such that $a+b\leq 1$,
we denote by $x_{a;b}$ the subword of $x$
going from the $(a\abs{x})$-th letter (taking integer part, and
inclusively) to the $((a+b)\abs{x})$-th letter (taking integer part, and
exclusively), so that $a$ indicates the position of the subword, and $b$
its length. If $a+b>1$ we cycle around $x$.

\begin{defi}
\label{subword}
Let $P_\ell$ be a family of properties of words, indexed by the integer
$\ell$. We say that
\[
\text{for any subword }x\text{ under }\mu_\ell,\; \Pr(P_\ell(x))\lesssim p(\ell)
\]
if for any $a,b\in [0;1]$, $b>0$,
whenever we pick a word $x$ according to
$\mu_\ell$ we have
\[
\Pr\left(P_\ell(x_{a;b}) \;|\; \abs{x},\, x_{0;a}\right)\lesssim p(\ell)\quad\text{if
}a+b\leq 1
\]
or
\[
\Pr\left(P_\ell(x_{a;b}) \;|\; \abs{x},\, x_{a+b-1;a}\right)\lesssim p(\ell)\quad\text{if
}a+b> 1
\]
and if moreover the constants implied in $\lesssim$ are uniform in $a$,
and, for each $\eps>0$, uniform when $b$ ranges in the interval $[\eps;1]$.
\end{defi}

That is, we pick a subword of a given length and ask the probability to
be bounded independently of whatever happened in the word up to this
subword (if the subword cycles around the end of the word, we condition
by everything not in the subword).

We also have to condition w.r.t.\ the length of the word since in the
definition of a random set of density $d$ under $\mu_\ell$ above, we made
a sampling for each length separately.

It would not be reasonable to ask that the constants be independent of
$b$ for arbitrarily small $b$. For example, if $\mu_\ell$ consists in
choosing uniformly a word of length $\ell$, then taking $b=1/\ell$
amounts to considering subwords of length $1$, which we are unable to say
anything interesting about.

We give a similar definition for properties depending on two words, but
we have to beware the case when they are subwords of the same word.

\begin{defi}
Let $P_\ell$ be a family of properties depending on two words, indexed by
the integer $\ell$. We say that
\[
\text{for any two disjoint subwords }x, y\text{ under }\mu_\ell,\;
\Pr(P_\ell(x,y))\lesssim p(\ell)
\]
if for any $a,b,a',b'\in [0;1]$
such that $b>0, b'>0, a+b\leq 1, a'+b'\leq 1$, whenever we
pick two independent words $x$, $x'$ according to $\mu_\ell$ we have
\[
\Pr\left(P_\ell(x_{a;b},\,x'_{a';b'}) \;|\; \abs{x},\, \abs{x'},\, x_{0;a},\, x'_{0;a'}\right)\lesssim p(\ell)
\]
and if for any $a,b,a',b'\in [0;1]$
such that $a< a+b\leq a'< a'+b'\leq 1$,
whenever we pick a word $x$ according to $\mu_\ell$, we have
\[
\Pr\left(P_\ell(x_{a,b},\, x_{a';b'}) \; |\; \abs{x},\, \abs{x'},\, x_{0;a},\, x_{a+b;a'}\right) \lesssim
p(\ell)
\]

We give similar definitions when $a+b>1$ or $a'+b'>1$, conditioning by
every subword not in $x_{a;b}$ or $x'_{a';b'}$.

Furthermore, we demand that the constants implied in $\lesssim$ be
uniform in $a$, $a'$, and, for each $\eps>0$, uniform when $b, b'$
range in the interval $[\eps;1]$.
\end{defi}

We are now ready to express the axioms we need on our random words.

\subsection{The Axioms}
\label{theaxioms}

Our first axiom states that $\mu_\ell$ consists of words of length
roughly $\ell$ up to some constant factor. This is crucial for the
hyperbolic local-global principle (Appendix~\ref{CHGP}).

\begin{enonce}{Axiom 1}
There is a constant $\kappa_1$ such that $\mu_\ell$ weights only words of
length between $\ell/\kappa_1$ and $\kappa_1\ell$.
\end{enonce}

Note this axiom applies to words picked under $\mu_\ell$, and not
especially subwords, so it does not rely on our definitions above. But of
course, if $\abs{x}\leq \kappa_1\ell$, then $\abs{x_{a;b}}\leq b \kappa_1\ell$.

Our second axiom states that subwords do not probably represent short
elements of the group.

\begin{enonce}{Axiom 2}
There are constants $\kappa_2,\beta_2$ such that for any subword $x$ under
$\mu_\ell$, for any $t\leq 1$, we have
\[
\Pr\left(\norm{x}\leq \kappa_2\abs{x}(1-t)\right)\lesssim (2m)^{-\beta_2t\abs{x}}
\]
uniformly in $t$.
\end{enonce}


Our next axiom controls the probability that two subwords are
almost inverse in the group. We will generally apply it with
$n(\ell)=O(\log\ell)$.

\begin{enonce}{Axiom 3}
There are constants $\beta_3$ and $\gamma_3$ such that for any function
$n=n(\ell)$,
for any two disjoint subwords $x,y$ under $\mu_\ell$,
the probability that there exist words $u$ and $v$ of
length at most $n$, such that $xuyv=e$ in $G$, is
roughly less than $(2m)^{\gamma_3 n}(2m)^{-\beta_3(\abs{x}+\abs{y})}$.
\end{enonce}

Our last axiom deals with algebraic properties of commutation with short
words.

\begin{enonce}{Axiom 4}
There exist constants $\beta_4$ and $\gamma_4$ such that,
for any function $n=n(\ell)$,
for any subword $x$ under $\mu_\ell$, the probability that there exist words $u$ and $v$ of length at most $n$, such that $ux=xv$ and $u\neq e$, $v\neq e$, is
roughly less than $(2m)^{\gamma n}(2m)^{-\beta_4 \abs{x}}$
\end{enonce}

If $G$ has big centralizers, this axiom will probably fail to be true.
We will see below (section~\ref{ax4}) that, in a hyperbolic group with
``strongly harmless'' torsion, the algebraic Axiom~4 is a consequence of
Axioms~1 and~3 combined with a more geometric axiom which we state now.

\begin{enonce}{Axiom 4'}
There are constants $\beta_{4'}$ and $\gamma_{4'}$ such that, for any $C>0$,
for any function $n=n(\ell)$,
for any subword $x$ under $\mu_\ell$, the probability that there exists a word $u$ of
length at most $n$ such that some cyclic permutation $x'$ of $xu$
satisfies $\norm{x'}\leq C\log\ell$, is roughly less than
$(2m)^{\gamma_{4'} n}
(2m)^{-\beta_{4'} \abs{x}}$.
\end{enonce}

\begin{rem}
\label{smallermeasure}
Let $\mu'_\ell$ be a family of measures such that $\mu'_\ell\lesssim
\mu_\ell$. As our axioms consist only in rough upper bounds, if the family
$\mu_\ell$ satisfy them, then so does the family $\mu'_\ell$.
\end{rem}

Note that as we condition every subword by whatever happened before
(i.e.\ by what the rest of the word is up to the position of the
subword), our axioms imply that subwords at different places are
essentially independent. This is of course true of plain random words,
but also of geodesic words and reduced words as we will see below.

\bigskip

In~\cite{Gro4}, p.~139--141, M.~Gromov uses similar-looking properties.
His $\mathbf{pr_1}$ is similar to our Axiom~2, and his $\mathbf{pr_3}$
controls the same kind of event as our Axiom~4. We no not use any
analogue of his $\mathbf{pr_2}$, and analogues of our Axioms~1 and~3 are
indeed present in~\cite{Gro4} but in a more ``diffuse'' way in the paper.
Also note that in~\cite{Gro4} emphasis is put on very small densities, so
that the properties considered therein are of the form ``such event is
realized with probability exponentially close to $1$'', whereas since we
work in large densities we have to get a precise control of the tails of
the distributions, and so our axioms take the form ``the probability of
a deviation of size $L$ from such event is at most $\exp(-\beta L)$'', with
a tight value of $\beta$ needed.  So our axioms (which have been found
independently of~\cite{Gro4}) are more precise quantitatively.

\subsection{The Theorem}
\label{thetheorem}

Our main tool is the following

\begin{thm}
\label{techmain}
Let $G$ be a non-elementary hyperbolic group with trivial virtual centre.
Let $\mu_\ell$ be a family of symmetric measures indexed by $\ell$, satisfying
Axioms~1, 2, 3 and 4. Let $R$ be a set of random words of density at most $d$ picked under $\mu_\ell$.

If $d<\min(\beta_2,\beta_3,\beta_4)$, then with probability
exponentially close to $1$ as $\ell\rightarrow\infty$, the random
quotient $G/\langle R \rangle$ is non-elementary hyperbolic, as well as
all the intermediate quotients $G/\langle R' \rangle$ with $R'\subset R$.
\end{thm}

Section~\ref{mainsection} is devoted to the proof.

\begin{rem}
\label{smallerquotient}
Remark~\ref{smallermeasure} tells that if the theorem applies to some family of
measures $\mu_\ell$, it applies as well to any family of measures
$\mu'_\ell\lesssim \mu_\ell$.
\end{rem}

\subsection{On torsion and Axiom 4}
\label{ax4}

We show here that in a hyperbolic group with ``harmless'' torsion,
Axioms~1, 3 and~4' imply Axiom~4. The proof makes the algebraic nature of
this axiom clear: in a hyperbolic group, it means
that subwords under $\mu_\ell$ are probably not torsion elements, neither
elements commuting with torsion elements, nor close to powers of
short elements.

Recall that the virtual centre of a hyperbolic group is the set of
elements whose action on the boundary at infinity is trivial. For basic
properties see~\cite{Ols2}.

\begin{defi}[ (Harmless torsion)]
\label{harmlesstorsion}

A torsion element in a hyperbolic group is said to be \emph{strongly
harmless} if its centralizer is
either finite or virtually $\Z$.

A torsion element is said to be \emph{harmless} if it is either strongly
harmless or lying in the virtual centre.

A hyperbolic group is said to be \emph{with (strongly) harmless
torsion} if each non-trivial torsion element is (strongly) harmless.
\end{defi}

Harmfulness is defined as the opposite of harmlessness.

For example, torsion-free groups are with harmless torsion, as well as
free products of free groups and finite groups. Strongly harmless torsion
is stable by free product, but harmless torsion is not.

\bigskip

Let $\mu_\ell$ be a measure satisfying Axioms~1, 3 and~4'.

\begin{prop}
\label{torsion}
The probability that, for a subword $x$ under $\mu_\ell$, there exists a
word $u$ of length at most $n=n(\ell)$ such that $xu$ is a torsion
element, is roughly less than $(2m)^{\gamma_{4'} n}(2m)^{-\beta_{4'}\abs{x}}$.
\end{prop}

\begin{dem}
In a hyperbolic group, there are only finitely many conjugacy classes of
torsion elements (see~\cite{GH}, p.~73). Let $L$ be the maximal length of
a shortest element of a conjugacy class of torsion elements, we have
$L<\infty$. Now every torsion element is conjugated to an element of
length at most $L$.

Suppose $xu$ is a torsion element. It follows from
Corollary~\ref{conj} that some cyclic permutation of it is conjugate to
an element of length at most $L$ by some word of length at most $\delta
\log_2\abs{xu}+C'_c+1$ where $C'_c$ is a constant depending on the group. In
particular, this cyclic conjugate has norm at most $L+2(\delta
\log_2\abs{xu}+C'_c+1)$.

Suppose, by Axiom~1, that $\abs{x}\leq \kappa_1\ell$.

There are $\abs{xu}\leq \kappa_1\ell+n$ cyclic conjugates of $xu$. The choice of
the cyclic conjugate therefore only introduces a polynomial factor in
$\ell$. Let $x'$ denote the cyclic conjugate of $xu$ at play.

Thus we have to evaluate the probability that $\norm{x'}\leq L+2(\delta
\log_2\abs{x'}+C'_c+1)$. As $L$ and $C'_c$ are mere constants, 
Axiom~4' precisely says that this probability is roughly less than
$(2m)^{\gamma_{4'} n}(2m)^{-\beta_{4'}\abs{x}}$.
\end{dem}

\begin{prop}
\label{preq}
Let $w\in G$. For any subword $x$ under $\mu_\ell$, the probability that
$x=w$ in $G$ is roughly less than $(2m)^{-\beta_3\abs{x}}$ (uniformly in
$w$).
\end{prop}

\begin{dem}
Suppose that the probability that a subword $x$ under $\mu_\ell$ is equal to
$w$ is equal to $p$. Then, by symmetry, the probability that an
independent disjoint subword $y$ with $\abs{y}=\abs{x}$ is equal to $w\~$
is equal to $p$ as well. So the probability that two disjoint subwords
$x$ and $y$ are inverse is at least $p^2$. But Axiom~3 tells (taking
$u=v=e$) that this probability is roughly at most
$(2m)^{-\beta_3(\abs{x}+\abs{y})}=(2m)^{-2\beta_3\abs{x}}$, hence
$p\lesssim (2m)^{-\beta_3\abs{x}}$.
\end{dem}

\begin{prop}
\label{commhyp}
Suppose $G$ has strongly harmless torsion, and that Axioms~1, 3 and~4'
are satisfied. Set $\beta=\min(\beta_3,\beta_{4'})$.

There is a constant $\gamma$ such that for any
subword $x$ under $\mu_\ell$, the
probability that there exist words $u, v$ of length at most
$n=n(\ell)$, such that $ux=xv$ in $G$, with $u, v$ not
equal to $e$, is roughly less than $(2m)^{\gamma n-\beta\abs{x}}$.

So Axiom~4 is satisfied with $\beta_4=\min(\beta_3,\beta_{4'})$.
\end{prop}

\begin{dem}
Denote by $x$ again a geodesic word equal to $x$ in $G$.

The words $u$ and $v$ are conjugate (by $x$), and are of length at most
$n$. After Corollary~\ref{conj} they are conjugate by a word $w$ of length at
most $Cn$ where $C$ is a constant depending only on $G$.

Let us draw the hyperbolic quadrilateral $xwuw\~x\~u\~$. This is a
commutation diagram between $xw$ and $u$.

\begin{center}
\includegraphics{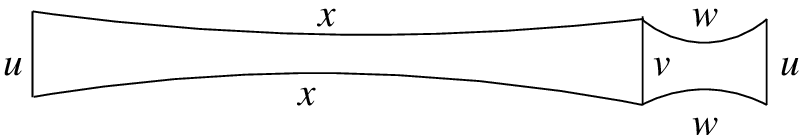}
\end{center}

The word $xw$ may or may not be a torsion element. The probability that
there exists a word $w$ of length at most $Cn$, such that $xw$ is a
torsion element, is roughly less than $(2m)^{\gamma_{4'} Cn-\beta\abs{x}}$ by
Proposition~\ref{torsion}. In this case we conclude.

Now suppose that $xw$ is not a torsion element. Then we can glue the
above diagram to copies of itself along their $u$-sides. This way we get two
quasi-geodesics labelled by $((xw)^n)_{n\in\Z}$ that stay at finite
distance from each other. The element $u$ acting on the first
quasi-geodesic gives the second one.

These two quasi-geodesics define an element $\tilde{x}$ in the boundary
of $G$. This element is of course stabilized by $xw$, but it is
stabilized by $u$ as well. This means that either $u$ is a hyperbolic
element, or (by strong harmlessness) that $u$ is a torsion element with
virtually cyclic centralizer.

The idea is that in this situation, $xw$ will lie close to some geodesic
$\Delta$ depending only on the short element $u$. As there are not many
such $\Delta$'s (and as the probability for a random word to be close to
a given geodesic behaves roughly like the probability to be close to
the origin), this will be unlikely.

First, suppose that $u$ is hyperbolic.
Let us use the same trick as above with the roles of $xw$ and
$u$ exchanged: glue the diagram above to copies of itself by the
$(xw)$-side. This defines two quasi-geodesics labelled by
$(u^n)_{n\in\Z}$, one of which goes to the other when acted upon by $xw$.

Namely, let $\Delta$ be a geodesic equivalent to $(u^n)$, and set
$\Delta'=xw\Delta$. As $xw$ stabilizes the
limit of $\Delta$, $\Delta'$ is equivalent to $\Delta$. But two
equivalent geodesics in a hyperbolic group stay at Hausdorff distance at
most $R_1$ where $R_1$ is a constant depending only on the group
(see~\cite{GH}, p.~119).

The distance from $xw$ to $\Delta'$ is equal to the distance from $e$ to
$\Delta$. By Proposition~\ref{cyclsubgeod} applied to $u^0=e$, this
distance is at most $\abs{u}+R_2$ where $R_2$ is a constant depending
only on $G$. Hence the distance from $xw$ to $\Delta$ is at most
$\abs{u}+R$ with $R=R_1+R_2$.  Let $y$ be a point on $\Delta$ realizing
this distance.  As $\abs{xw}\leq \abs{x}+\abs{w}$, we have $\abs{y}\leq
\abs{x}+\abs{w}+\abs{u}+R$.  There are at most
$2\abs{x}+2\abs{w}+2\abs{u}+2R+1$ such possible points on $\Delta$ (since
$\Delta$ is a geodesic). For each of these points, the probability that
$x$ falls within distance $\abs{u}+R+\abs{w}$ of it is roughly less than
$(2m)^{\abs{u}+R+\abs{w}} (2m)^{-\beta\abs{x}}$ by Proposition~\ref{preq}
applied to all of these points. So the probability that $x$ falls within
distance less than $\abs{u}+R+\abs{w}$ of any one of the possible $y$'s
on a given geodesic $\Delta$ is roughly less than
$(2\abs{x}+2\abs{w}+2\abs{u}+2R+1)(2m)^{\abs{u}+R+\abs{w}}
(2m)^{-\beta\abs{x}}$ which in turn is roughly less than
$(2m)^{Cn-\beta\abs{x}}$ as $\abs{w}\leq Cn$ and $R$ is a constant.

This was for one fixed $u$. But each different $u$ defines a different
$\Delta$. There are at most $(2m)^{\abs{u}}\leq(2m)^{n}$ possibilities
for $u$. Finally, the probability that $x$ falls within distance
$R+\abs{w}$ of any one of the geodesics defined by these $u$'s is less
than $(2m)^{n+Cn-\beta\abs{x}}$ as was to be shown. Thus we can conclude
when $u$ is hyperbolic.

Second, if $u$ is a torsion element with virtually cyclic centralizer $Z$, we
use a similar argument.
Let $L$ as above be the maximal length
of a shortest element of a conjucacy class of a torsion element. By
Proposition~\ref{gconj}, $u$ is conjugate to some torsion element $u'$ of
length at most $L$ by a conjugating word $v$ with $\abs{v}\leq
\abs{u}/2+R_1$ where $R_1$ is a constant. The centralizer of $u'$ is
$Z'=vZv\~$. We know that $xw\in Z$.

There are two subcases: either $Z$ is finite or $Z$ is virtually $\Z$.

Let us begin with the former. If $Z$ is finite, let $\norm{Z}$ be the
maximal norm of an element in $Z$. We have $\norm{Z}\leq
2\abs{v}+\norm{Z'}$. Let $R_2=\max \norm{Z'}$ when $u'$ runs through all
torsion elements of norm at most $L$. As $xw$ lies
in $Z$ we have $\norm{x}\leq \abs{w}+\norm{Z}\leq \abs{w}+2\abs{v}+R_2\leq
\abs{w}+\abs{u}+2R_1+R_2$. So by Proposition~\ref{preq} the probability
of this event is roughly less than
$(2m)^{\abs{w}+\abs{u}+2R_1+R_2}\lesssim (2m)^{Cn+n}$ as $\abs{w}\leq Cn$
and as $R_1,R_2$ are mere constants.

Now if $Z$ is virtually $\Z$, let $\Delta$ be a geodesic joining the two
limit points of $Z$. The element $u'$ defined above
stabilizes the endpoints of the geodesic $v\Delta$, and so does
$vxwv\~$.

By Corollary~\ref{diststab}, $vxwv\~$ lies at distance at most
$R(v\Delta)$ from $v\Delta$. As there are only a finite number of
torsion elements $u'$ with $\norm{u'}\leq L$, the supremum $R$ of the
associated $R(v\Delta)$ is finite, and so, independently of $u$, the
distance between $vxwv\~$ and $v\Delta$ is at most $R$.

Now $\dist(xw, \Delta)\leq \abs{v}+\dist(xwv\~,\Delta)=
\abs{v}+\dist(vxwv\~,v\Delta)\leq \abs{v}+R$ and we conclude exactly as in the
case when $u$ was hyperbolic,
using that $\abs{v}\leq \abs{u}/2+R_1$. This ends the proof in case $u$
is a torsion element with virtually cyclic centralizer.
\end{dem}

\section{Applications of the main theorem}

We now show how Theorem~\ref{techmain} leads, with some more work, to the
theorems on random quotients by plain words, reduced words and geodesic
words given in the introduction.

We have three things to prove:
\begin{itemize}
\item first, that these three models of a random
quotient satisfy our axioms with the right critical densities;
\item second, as
Theorem~\ref{techmain} only applies to hyperbolic groups with strongly
harmless torsion (instead of harmless torsion), we have to find a way to get rid of the virtual centre;
\item third, we have to
prove triviality for densities above the critical one.
\end{itemize}

Once this is done, Theorems~\ref{reducedmain}, \ref{geodmain} and
\ref{main} will be proven.

\bigskip

We will have to work differently if we consider quotients by plain random
words, by random reduced words or by random geodesic words.

For instance, satisfaction of the axioms is very different for plain
words and for geodesic words, because in plain random words, two given
subwords fo the same word are chosen independently, which is not the case
at all \emph{a priori} for a geodesic word.

Furthermore, proving triviality of a quotient involves small scale
phenomena, which are very different in our three models of random words
(think of a random quotient of $\Z$ by random words of $\ell$ letters
$\pm 1$ or by elements of size exactly $\ell$).

These are the reasons why the next three sections are divided in cases,
and why we did not include these properties in a general and technical
theorem such as Theorem~\ref{techmain}.

\bigskip

Note that it is natural to express the critical densities in terms of the
$\ell$-th root of the total number of words of the kind considered, that
is, in base $2m$ for plain words, $2m-1$ for reduced words and $(2m)^g$
for geodesic words.

\subsection{Satisfaction of the axioms}

\subsubsection{The case of plain random words}

We now take as our measure for random words the uniform measure on all
words of length $\ell$. Axiom~1 is satisfied by definition.

In this section, we denote by $B_\ell$ (as ``Brownian'') a random word of
length $\ell$ uniformly chosen among all $(2m)^\ell$ possible words.

Recall $\theta$ is the gross cogrowth of the group, that is, the number
of words of length $\ell$ which are equal to $e$ in the group is roughly
$(2m)^{\theta\ell}$ for even $\ell$.

Recall the alternate definition of gross cogrowth given in the
introduction: the exponent of return to $e$ of the random walk in $G$ is
$1-\theta$. This is at the heart of what follows.

We will show that

\begin{prop}
Axioms~1, 2, 3, 4' are satisfied by plain random uniformly chosen words,
with exponent $1-\theta$ (in base $2m$).
\end{prop}

By definition, disjoint subwords of a uniformly taken random word are
independent. So we do not have to care at all with the conditional
probabilities of the axioms (contrary to the case of geodesic words
below). Conditionnally to anything else, every subword $x$ follows the
law of $B_{\abs{x}}$.

The definition of gross cogrowth only applies to even lengths. If $\ell$
is odd, either there are some relations of odd length in the presentation
of the group, and then the limits holds, or there are no such relations,
and the number of words of length $\ell$ equal to $e$ is zero. In any
case, this number is $\lesssim (2m)^{\theta \ell}$.

This is a delicate (but irrelevant) technical point: We should care with
parity of the length of words. If there are some relations of odd length
in our group, then the limit in the definition of gross cogrowth is valid
regardless of parity of $\ell$, but in general this is not the case (as
is examplified by the free group). In order to get valid results for any
length, we therefore often have to replace a $\approx$ sign with a
$\lesssim$ one. In many cases, our statements of the form
``$\Pr(\ldots)\lesssim f(\ell)$'' could in fact be replaced by
``$\Pr(\ldots)\approx f(\ell)$ if $\ell$ is even or if there are
relations of odd length, and $\Pr(\ldots)=0$ otherwise''. Here is the
first example of such a situation.

\begin{prop}
The probability that $B_\ell$ is equal to $e$ is roughly less than
$(2m)^{-(1-\theta)\ell}$.
\end{prop}

\begin{dem}
Alternate definition.
\end{dem}

%

\begin{prop}
\label{randomnorm}
\[
\Pr(\norm{B_\ell}\leq \ell') \lesssim
(2m)^{-(1-\theta)\left(\ell-\frac{\theta}{1-\theta}\ell'\right)}
\]
uniformly in $\ell'\leq \ell$.

In particular, the escaping speed is at least $\frac{1-\theta}{\theta}$.
So Axiom~2 is satisfied with $\kappa_2=\frac{1-\theta}{\theta}$ and
$\beta_2=1-\theta$.
\end{prop}

\begin{dem}
For any $L$ between $0$ and $\ell'$, we have that
\[\Pr(B_{\ell+L}=e)\geq (2m)^{-L} \Pr(\norm{B_\ell}=L)\]

But $\Pr(B_{\ell+L}=e)\lesssim(2m)^{-(1-\theta)(\ell+L)}$ (and
this is uniform in $L\leq \ell$ since in any case,
$\ell+L$ is at least equal to $\ell$), hence the evaluation for a
given $L$.

Now, summing over $L$ between $0$ and $\ell'$ introduces only a
subexponential factor in $\ell$.
\end{dem}

\begin{prop}
\label{qeind}
The probability that, for two independently chosen words
$B_\ell$ and $B'_{\ell'}$, there exist words $u$ and $v$ of length at
most $n=n(\ell)$, such that $B_\ell u B'_{\ell'} v=e$ in $G$, is roughly
less than $(2m)^{(2+2\theta)n}(2m)^{-(1-\theta)(\ell+\ell')}$.

That is, Axiom~3 is satisfied with exponent $1-\theta$.
\end{prop}

\begin{dem}
For any word $u$, we have $\Pr(B_{\abs{u}}=u)\geq (2m)^{-\abs{u}}$.

So let $u$ and $v$ be any two fixed words of length at most $n$. We have
\[\Pr(B_{\ell+\abs{u}+\ell'+\abs{v}}=e)\geq
(2m)^{-\abs{u}-\abs{v}}\Pr(B_\ell u B'_{\ell'}v=e)\]

We know that $\Pr(B_{\ell+\abs{u}+\ell'+\abs{v}}=e) \lesssim
(2m)^{-(1-\theta)\left(\ell+\abs{u}+\ell'+\abs{v}\right)}$.

So
$\Pr(B_\ell u B'_{\ell'}v=e) \lesssim (2m)^{\theta(\abs{u}+\abs{v})}
(2m)^{-(1-\theta)(\ell+\ell')}$.

Now there are $(2m)^{\abs{u}+\abs{v}}$ choices for $u$ and $v$.
\end{dem}

\begin{prop}
The probability that there exists a word $u$ of length at most
$n=n(\ell)$, such that some cyclic conjugate
of $B_\ell u$ is of norm less than $C\log\ell$, is roughly less than
$(2m)^{(1+\theta)n} (2m)^{-(1-\theta)\ell}$.

So Axiom~4' is satisfied with exponent $1-\theta$.
\end{prop}

\begin{dem}
As above, for any word $u$, we have $\Pr(B_{\abs{u}}=u)\geq
(2m)^{-\abs{u}}$. So any property of $B_\ell u$ occurring with some
probability will occur for $B_{\ell+\abs{u}}$ with at least
$(2m)^{-\abs{u}}$ times this probability. We now work with
$B_{\ell+\abs{u}}$.

Any cyclic conjugate of a uniformly chosen random word is itself a
uniformly chosen random word, so we can assume that the cyclic conjugate
at play is $B_{\ell+\abs{u}}$ itself. There are $\ell+\abs{u}$ cyclic
conjugates, so the choice of the cyclic conjugate only introduces a
subexponential factor in $\ell$ and $\abs{u}$.

But we just saw above in Proposition~\ref{randomnorm}
that the probability that
$\norm{B_{\ell+\abs{u}}}\leq L$ is roughly less than
$(2m)^{-(1-\theta)\left(\abs{u}+\ell-\frac{\theta}{1-\theta}L\right)}$.

Summing over the $(2m)^{\abs{u}}$ choices for $u$ yields the desired
result, taking $L=C\log \ell$.
\end{dem}

So plain random words satisfy our axioms.

\subsubsection{The case of random geodesic words}

The case of geodesic words is a little bit more clever, as subwords of a
geodesic word are not \emph{a priori} independent.

For each element $x\in G$ such that $\norm{x}=\ell$, fix once and for all
a representation of $x$ by a word of length $\ell$. We are going to prove
that when $\mu_\ell$ is the uniform law on the sphere of radius $\ell$ in
$G$, Axioms~1-4' are satisfied.

Recall that $g$ is the growth of the group: by definition, the number of
elements of length $\ell$ in $G$ is roughly $(2m)^{g\ell}$. As $G$ is
non-elementary we have $g>0$ (otherwise there is nothing to prove).

\begin{prop}
Axioms~1, 2, 3, 4' are satisfied by random uniformly chosen elements of
norm $\ell$, with exponent $1/2$ (in base $(2m)^g$).
\end{prop}

Our proofs also work if $\mu_\ell$ is the uniform measure on the spheres
of radius between $\ell-L$ and $\ell+L$ for any fixed $L$. We will use
this property later.

Note that Axioms~1 and~2 are trivially satisfied for geodesic words, with
$\kappa_1=\kappa_2=1$ and $\beta_2=\infty$.

\bigskip

The main obstacle is that two given subwords of a geodesic word are not
independent. We are going to replace the model of randomly chosen
elements of length $\ell$ by another model with more independence, and
prove that these two models are roughly equivalent.

Let $X_\ell$ denote a random uniformly chosen element on the sphere of
radius $\ell$ in $G$.  For any $x$ on this sphere, we have
$\Pr(X_\ell=x)\approx (2m)^{-g\ell}$.

Note that for any $\eps>0$, for any $\eps\ell\leq L\leq \ell$ the rough
evaluation of the number of points of length $L$ by $(2m)^{gL}$ can by
taken uniform for $L$ in this interval (take $\ell$ so that $\eps \ell$
is big enough).

\bigskip

First, we will change a little bit the model of random geodesic words.
The axioms above use a strong independence property of subwords of the
words taken. This independence is not immediately satisfied for subwords
of a given random geodesic word (for example, in the hyperbolic group
$F_2\times\Z/2\Z$, the occurrence of a generator of order $2$ somewhere
prevents it from occurring anywhere else in a geodesic word). So we will
cheat and consider an alternate model of random geodesic words.

For a given integer $N$, let $X^N_\ell$ be the product of $N$ random
uniformly chosen geodesic words of length $\ell/N$. We will compare the
law of $X_\ell$ to the law of $X^N_\ell$.

Let $x\in G$ such that $\norm{x}=\ell$. We have $\Pr(X_\ell=x)\approx
(2m)^{-g\ell}$. Let $x=x_1x_2\ldots x_N$ where each $x_i$ is of length
$\ell/N$. The probability that the $i$-th segment of $X^N_\ell$ is equal
to $x_i$ is roughly $(2m)^{-g\ell/N}$. Multiplying, we get
$\Pr(X^N_\ell=x)\approx (2m)^{-g\ell}$.

Thus, if $P$ is a property of words, we have for any given $N$ that
\[
\Pr(P(X_\ell))\lesssim \Pr(P(X^N_\ell))
\]

(The converse inequality is false as the range of values of $X^N_\ell$ is
not contained in that of $X_\ell$.)

Of course, the constants implied in $\lesssim$ depend on $N$. We are
stating that for any fixed $N$, when $\ell$ tends to infinity the law of
the product of $N$ words of length $\ell/N$ encompasses the law of
$X_\ell$, and \emph{not} that for a given $\ell$, when $N$ tends to
infinity the law of $N$ words of length $\ell$ is close to the law of a
word of length $N\ell$, which is false.

We are going to prove the axioms for $X^N_\ell$ instead of $X_\ell$. As the
axioms all state that the probability of some property is roughly less
than something, these evaluations will be valid for $X_\ell$.

The $N$ to use will depend on the length of the subword at play in the
axioms. With notations as above, if $x_{a;b}$ is a subword of length
$b\ell$ of $X_\ell$, we will choose an $N$ such that $\ell/N$ is small
compared to $b\ell$, so that $x_{a;b}$ can be considered the product of
a large number of independently randomly chosen smaller geodesic words.
This is fine as our axioms precisely \emph{do not} require the evaluations to be
uniform when the relative length $b$ tends to $0$.

\bigskip

First, we need to study multiplication by a random geodesic word.

Let $(x|y)$ denote the Gromov product of two elements $x,y\in G$.
That is, $(x|y)=\frac12 \left(\norm{x}+\norm{y}-\norm{x\~y}\right)$.

\begin{prop}
\label{ravi}
Let $x\in G$ and $L\leq \ell$. We have
\[
\Pr\left((x|X_\ell)\geq L\right)\lesssim (2m)^{-gL}
\]
uniformly in $x$ and $L\leq \ell$.
\end{prop}

\begin{dem}
Let $y$ be the point at distance $L$ on a geodesic joining $e$ to
$x$.
By the triangle-tripod transformation in $exX_\ell$, the inequality
$(x|X_\ell)\geq L$ means that $X_\ell$ is at distance at most
$\ell-L+4\delta$ from $y$. There are roughly at most
$(2m)^{g(\ell-L+4\delta)}$ such points. Thus, the probability that
$X_\ell$ is equal to one of them is roughly less than
$(2m)^{g(\ell-L+4\delta)-g\ell}\approx(2m)^{-gL}$.

Let us show that this evaluation can be taken uniform in $L\leq \ell$.
The problem comes from the evaluation of the number of points  at
distance at most $\ell-L+4\delta$ from $y$ by
$(2m)^{g(\ell-L+4\delta)}$: when $\ell-L+4\delta$ is not large
enough, this cannot be taken uniform.  So take some $\eps>0$ and first
suppose that $L\leq(1-\eps)\ell$, so that $\ell-L+4\delta\geq
\eps'\ell$ for some $\eps'>0$. The evaluation of the number of points at
distance at most $\ell-L+4\delta$ from $y$ by
$(2m)^{g(\ell-L+4\delta)}$ can thus be taken uniform in $L$ in this
interval.

Second, let us suppose that $L\geq (1-\eps)\ell$. Apply the trivial
estimate that the number of points at distance $\ell-L+4\delta\leq
\eps\ell+4\delta$ from $y$ is less than
$(2m)^{\eps\ell+4\delta}$. The probability that $X_\ell$ is equal to one
of them is roughly less than $(2m)^{\eps\ell-g\ell}\leq
(2m)^{-(g-\eps)L}$ uniformly for these values of $L$.

So for any $\eps$, we can show that for any $L\leq \ell$, the
probability at play is uniformly roughly less than $(2m)^{-(g-\eps)L}$.
Writing out the definition shows that this exacly says that our
probability is less than $(2m)^{-gL}$ uniformly in $L$.
\end{dem}

\begin{cor}
\label{multgeod}
Let $x\in G$ and $L\leq 2\ell$. Then
\[
\Pr\left(\norm{xX_\ell}\leq \norm{x}+\ell-L\right)
\lesssim (2m)^{-gL/2}
\]
and
\[
\Pr\left(\norm{X_\ell x}\leq \norm{x}+\ell-L\right)
\lesssim (2m)^{-gL/2}
\]
uniformly in $x$ and $L$.
\end{cor}

\begin{dem}
Note that the second case follows from the first one applied to $x\~$ and
$X_\ell\~$, and symmetry of the law of $X_\ell$.

For the first case, apply Proposition~\ref{ravi} to $X_\ell$ and $x\~$ and write out the definition
of the Gromov product.
\end{dem}

\begin{prop}
\label{multgeods}
For any fixed $N$, uniformly for any $x\in G$ and any $L\leq 2\ell$ we have
\[
\Pr\left(\norm{xX^N_\ell}\leq \norm{x}+\ell-L\right)
\lesssim (2m)^{-gL/2}
\]
and
\[
\Pr\left(\norm{X^N_\ell x}\leq \norm{x}+\ell-L\right)
\lesssim (2m)^{-gL/2}
\]
\end{prop}

\begin{dem}
Again, note that the second inequality follows from the first one by
taking inverses and using symmetry of the law of $X^N_\ell$.

Suppose $\norm{xX^N_\ell}\leq\norm{x}+\ell-L$.
Let $x_1,x_2,\ldots,x_N$ be $N$ random uniformly chosen geodesic words of
length $\ell/N$. Let $L_i\leq 2\ell/N$ such that $\norm{xx_1\ldots
x_i}=\norm{xx_1\ldots x_{i-1}}+\ell/N-L_i$. By $N$ applications of
Corollary~\ref{multgeod}, the probability of such an event is roughly
less than $(2m)^{-g\eps\sum L_i/2}$. But·we have $\sum L_i\geq L$. Now the
number of choices for the $L_i$'s is at most $(2\ell)^N$, which
is polynomial in $\ell$, hence the proposition.
\end{dem}

Of course, this is not uniform in $N$.

\bigskip

We now turn to satisfaction of Axioms~3 and~4' (1 and 2 being trivially
satisfied). We work under the model of $X^N_\ell$. Let $x$ be a subword
of $X^N_\ell$. By taking $N$ large enough (depending on $\abs{x}/\ell$), we
can suppose that $x$ begins and ends on a multiple of $\ell/N$. If not,
throw away an initial and final subword of $x$ of length at most
$\ell/N$. In the estimates, this will change $\norm{x}$ in
$\norm{x}-2\ell/N$ and, if the estimate to prove is of the form
$(2m)^{-\beta\norm{x}}$, for each $\eps>0$ we can find an $N$ such that
we can prove the estimate $(2m)^{-\beta(1-\eps)\norm{x}}$. Now if
something is roughly less than $(2m)^{-\beta(1-\eps)\norm{x}}$ for every
$\eps>0$, it is by definition roughly less than $(2m)^{-\beta\norm{x}}$.

Note that taking $N$ depending on the relative length $\abs{x}/\ell$ of
the subword is correct since we did not ask the estimates to be uniform
in this ratio.

The main advantage of this model is that now, the law of a subword is
independent of the law of the rest of the word, so we do not have to care
about the conditional probabilities in the axioms.

\begin{prop}
Axiom~3 is satisfied for random geodesic words, with exponent $g/2$.
\end{prop}

\begin{dem}
Let $x$ and $y$ be subwords. The word $x$ is a product of $N\abs{x}/\ell$
geodesic words of length $\ell/N$, and the same holds for $y$. Now take
two fixed words $u$, $v$, and let us evaluate the probability that
$xuyv=e$.

Fix some $L\leq \ell$, and suppose $\norm{x}=L$. By
Proposition~\ref{multgeods} starting at $e$, this occurs with probability
$(2m)^{-g(\abs{x}-L)/2}$. Now we have $\norm{xu}\geq L-\norm{u}$, but
$\norm{xuy}=\norm{v\~}$. By Proposition~\ref{multgeods} starting at $xu$
this occurs with probability $(2m)^{-g(L-\norm{u}+\abs{y}-\norm{v})/2}$.

So the total probability is at most the number of choices for $u$ times the
number of choices for $L$ times $(2m)^{-g(\abs{x}-L)/2}$ times
$(2m)^{-g(L-\norm{u}+\abs{y}-\norm{v})/2}$. Hence the proposition.
\end{dem}

\begin{prop}
Axiom 4' is satisfied for random geodesic words, with exponent $g/2$.
\end{prop}

\begin{dem}
Taking notations as in the definitions, let $x$ be a subword of
$X^N_\ell$ of length $b\ell$ with $b\leq 1$. The law of $x$ is
$X^{bN}_{b\ell}$.

Note that applying Proposition~\ref{multgeods} starting with the neutral
element $e$ shows that $\Pr(\norm{x}\leq L)\lesssim
(2m)^{-g(\abs{x}-L)/2}$.

Fix a $u$ of length at most $n$ and consider a cyclic conjugate $y$ of $xu$.

First, suppose that the cutting made in $xu$ to get the cyclic conjugate
$y$ was made in $u$, so that $y=u''xu'$ with $u=u'u''$. In this case, we
have $\norm{y}\geq \norm{x}-\norm{u''}-\norm{u}\geq \norm{x}-\abs{u}$,
and so we have $\Pr(\norm{y}\leq C\log\ell)\leq \Pr(\norm{x}\leq
C\log\ell+\norm{u})\lesssim (2m)^{-g(\abs{x}-C\log\ell-\abs{u})/2}\approx
(2m)^{g\abs{u}/2-g\abs{x}/2}$.

Second, suppose that the cutting was made in $x$, so that $y=x''ux'$ with $x=x'x''$.

Up to small words of length at most $\ell/N$ at the beginning and end of
$x$, the words $x'$ and $x''$ are products of randomly chosen geodesic
words of length $\ell/N$.

Apply Proposition~\ref{multgeods} starting with the element $u$,
multiplying on the right by $x'$, then on the left by $x''$. This
shows that $\Pr(\norm{y}\leq \norm{u}+\abs{x'}+\abs{x''}-L)\lesssim
(2m)^{-gL/2}$, hence the evaluation, taking
$L=\abs{x'}+\abs{x''}+\norm{u}-C\log\ell$.

To conclude, observe that there are at most $(2m)^{\abs{u}}$ choices for
$u$ and at most $\abs{x}+\abs{u}$ choices for the cyclic conjugate, hence
an exponential factor in $\abs{u}$.
\end{dem}

\subsubsection{The case of random reduced words}

Recall $\eta$ is the cogrowth of the group $G$, i.e.\ the number of
reduced words of length $\ell$ which are equal to $e$ is roughly
$(2m-1)^{\eta \ell}$.

Here we have to suppose $m>1$. (A random quotient of $\Z$ by reduced words
of length $\ell$ is $\Z/\ell\Z$.)

\begin{prop}
Axioms~1, 2, 3, 4' are satisfied by random uniformly chosen reduced words,
or random uniformly chosen cyclically reduced words, with exponent
$1-\eta$ (in base $2m-1$).
\end{prop}

The proof follows essentially the same lines as that for plain random
words. We do not include it explicitly here.

Nevertheless, there are two changes encountered.

The first problem is that we do not have as much independence
for reduced words as for plain words. Namely, the occurrence of a
generator at position $i$ prevents the occurrence of its inverse at
position $i+1$.

We solve this problem by noting that, though the $(i+1)$-th letter depends
on what happened before, the $(i+2)$-th letter does not depend too much
(if $m>1$).

Indeed, say the $i$-th letter is $x_j$. Now it is immediate to check that the
$(i+2)$-th letter is $x_j$ with probability $1/(2m-1)$, and is each other
letter with probability $(2m-2)/(2m-1)^2$. This is close to a uniform
distribution up to a factor of $(2m-2)/(2m-1)$.

This means that, conditioned by the word up to the $i$-th letter, the law
of the word read after the $(i+2)$-th letter is, up to a constant
factor, an independently chosen random reduced word.

This is enough to allow to prove satisfaction of the axioms for random
reduced words by following the same lines as for plain random words.

The second point to note is that a reduced word is not necessarily
cyclically reduced. The end of a reduced word may collapse with the
beginning. Collapsing along $L$ letters has probability precisely
$(2m-1)^{-L}$, and the induced length loss is $2L$. So this
introduces an exponent $1/2$, but the cogrowth $\eta$ is bigger
than $1/2$ anyway.

In particular, everything works equally fine with reduced and cyclically
reduced words (the difference being non-local), with the same critical
density $1-\eta$.

\subsection{Triviality of the quotient in large density}

Recall $G$ is a hyperbolic group generated by $S=a_1^{\pm
1},\ldots,a_m^{\pm 1}$. Let $R$ be a set of $(2m)^{d\ell}$ randomly
chosen words of length $\ell$. We study $G/\langle R \rangle$.

As was said before, because triviality of the quotient involves
small-scale phenomena, we have to work separately on plain random words,
reduced random words or random geodesic words.

\label{triv}

Generally speaking, the triviality of the quotient reduces essentially to
the following fact, which is analogue to the fact that two (say generic
projective complex algebraic) submanifolds whose sum of dimensions is
bigger than the ambient dimension do intersect (cf.~our discussion of the
density model of random groups in the introduction).

\begin{enonce}{Basic intersection theory for random sets}
\label{bitrs}
Let $S$ be a set of $N$ elements. Let $\alpha, \beta$ be two numbers in $[0;1]$ such that $\alpha+\beta>1$.
Let $A$ be a given part of $S$ of cardinal
$N^\alpha$. Let $B$ be a set of $N^\beta$ randomly uniformly chosen
elements of $S$. Then $A\cap B\neq \varnothing$ with probability tending
to $1$ as $N\rightarrow \infty$ (and the intersection is arbitrarily
large with growing $N$).
\end{enonce}

This is of course a variation on the probabilistic pigeon-hole
principle where $A=B$.
 
\begin{rem*}
Nothing in what follows is specific to quotients of hyperbolic groups:
for the triviality of a random quotient by too many relators, any group
(with $m>1$ in the reduced word model and $g>0$ in the geodesic word
model) would do.
\end{rem*}

\subsubsection{The case of plain random words}

We suppose that $d>1-\theta$.

Recall that $\theta$ is the gross cogrowth of the group, i.e.\ that
\[
\theta=\lim_{
\ell\rightarrow \infty, \ell\text{ even}} \frac1\ell \log_{2m} \#\{w\in
B^\ell,
w=e\text{ in }G\}
\]

\bigskip

We want to show that the random quotient $G/\langle R \rangle$ is either $\{1\}$ or
$\Z/2\Z$. Of course the case $\Z/2\Z$ occurs when $\ell$ is even and when
the presentation of $G$ does not contain any odd-length relation.

To use gross cogrowth, we have to distinguish according to parity of $\ell$. We
will treat only the least simple case when $\ell$ is even. The other case is
even simpler.

Rely on the intersection theory for random sets stated above. Take for
$A$ the set of all words of length $\ell-2$ which are equal to $e$ in
$G$. There are roughly $(2m)^{\theta(\ell-2)}\approx (2m)^{\theta\ell}$ of them.
Take for $B$ the set made of the random words of $R$ with the last two letters
removed, and recall that $R$ consists of $(2m)^{d\ell}$ randomly chosen
words with $d>1-\theta$.

Apply the intersection principle: very probably, these sets will
intersect. This means that in $R$, there will probably be a word of the
form $wab$ such that $w$ is trivial in $G$ and $a,b$ are letters in
$S$ or $S^{-1}$.

This means that in the quotient $G/\langle R \rangle$, we have $ab=e$.

Now as $d+\theta>1$ this situation occurs arbitrarily many times as
$\ell\rightarrow \infty$. Due to our uniform choice of random words, the
$a$ and $b$ above will exhaust all pairs of generators of $S$ and $S\~$.

Thus, in the quotient, the product of any two generators $a,b\in S\cup
S\~$ is equal to $e$. Hence the quotient is either trivial or $\Z/2\Z$
(and is it trivial as soon as $\ell$ is odd or the presentation of $G$
contains odd-length relators).

This proves the second part of Theorem~\ref{main}.

\subsubsection{The case of random geodesic words}

When taking a random quotient by geodesic words of the same length, some
local phenomena may occur. For example, the quotient of $\Z$ by any
number of randomly chosen elements of norm $\ell$ will be $\Z/\ell\Z$.
Think of the occurrence of either $\{e\}$ or $\Z/2\Z$ in a quotient by
randomly chosen non-geodesic words.

In order to avoid this phenomenon, we consider a random quotient by
randomly chosen elements of norm comprised between $\ell-L$ and
$\ell+L$ for some fixed small $L$. Actually we will take $L=1$. 

Recall $g$ is the growth of the group, that is, the number of elements of
norm $\ell$ is roughly $(2m)^{g\ell}$, with $g>0$ as $G$ is
non-elementary.

We now prove that a random quotient of any group $G$ by $(2m)^{d\ell}$
randomly chosen elements of norm $\ell-1$, $\ell$ and $\ell+1$, with
$d>g/2$, is trivial with probability tending to $1$ as
$\ell\rightarrow\infty$. 

(By taking $(2m)^{d\ell}$ elements of norm $\ell$, $\ell+1$ or
$\ell-1$
we mean either taking $(2m)^{d\ell}$ elements of each of these norms,
or taking $1/3$ at each length, or deciding for each element with a given
positive probability what its norm will be, or any other roughly
equivalent scheme.)

\bigskip

Let $a$ be any of the generators of the group. Let $x$ be any element of
norm $\ell$. The product $xa$ is either of norm $\ell$, $\ell+1$ or
$\ell-1$.

Let $S$ be the sphere of radius $\ell$, we have $\abs{S}\approx
(2m)^{g\ell}$.

Let $R$ be the set of random words taken. Taking $d>g/2$ precisely
amounts to taking more than ${\abs{S}}^{1/2}$ elements of $S$.

Let $R'$ be the image of $R$ by $x\mapsto xa$. By an easy variation on
the probabilistic pigeon-hole principle applied to $R$, there will very
probably be one element of $R$ lying in $R'$. This means that $R$ will
contain elements $x$ and $y$ such that $xa=y$. Hence, $a=e$ in the
quotient by $R$.

As this will occur for any generator, the quotient is trivial. This
proves the second part of Theorem~\ref{geodmain}.

\subsubsection{The case of random reduced words}

For a quotient by random reduced words in density $d>1-\eta$ (where
$\eta$ is the cogrowth of the group), the proof of triviality is nearly
identical to the case of a quotient by plain random words, except that
in order to have the number of words taken go to infinity, we have to
suppose that $m\geq 2$.

\subsection{Elimination of the virtual centre}
\label{virt}

Theorem~\ref{techmain} only applies to random quotients of hyperbolic
groups with strongly harmless torsion. We have to show that the presence
of a virtual centre does not change random quotients. The way to do this
is simply to quotient by the virtual centre; but, for example, geodesic
words in the quotient are not geodesic words in the original group, and
moreover, the growth, cogrowth and gross cogrowth may be different. Thus
something should be said.

Recall the virtual centre of a hyperbolic group is the set of elements
whose action on the boundary at infinity is trivial. It is a normal
subgroup (as it is defined as the kernel of some action). It is finite,
as any element of the virtual centre has force $1$ at each point of the
boundary, and in a (non-elementary) hyperbolic group, the number of
elements having force less than a given constant at some point is finite
(cf.~\cite{GH}, p.~155). See~\cite{Ols2} or~\cite{Ch3} for an exposition
of basic properties and to get an idea of the kind of problems arising
because of the virtual centre.

Let $H$ be the virtual centre of $G$ and set $G'=G/H$. The quotient $G'$
has no virtual centre.

\subsubsection{The case of plain or reduced random words}

Note that the set $R$ is the same, since the notion of plain random word
or random reduced word is defined independently of $G$ or $G'$.

As $(G/H)/\langle R \rangle=(G/\langle R \rangle)/H$, and as a quotient
by a finite normal subgroup is a quasi-isometry, $G/\langle R \rangle$
will be infinite hyperbolic if and only if $G'/\langle R \rangle$ is.

So in order to prove that we can assume a trivial virtual centre, it is
enough to check that $G$ and $G/H$ have the same cogrowth and gross
cogrowth, so that the notion of a random quotient is really the same.

We prove it for plain random words, as the case of reduced words is
identical with $\theta$ replaced with $\eta$ and $2m$ replaced with
$2m-1$.

\begin{prop}
\label{prneighbor}
Let $H$ be a subset of $G$, and $n$ an integer. Then
\[
\Pr(\exists u\in G, \abs{u}= n, B_\ell u\in H) \leq (2m)^n
\Pr(B_{\ell+n}\in H)
\]
\end{prop}

\begin{dem}
Let $H_n$ be the $n$-neighborhood of $H$ in $G$. We have that
$\Pr(B_{\ell+n}\in H)\geq (2m)^{-n} \Pr(B_\ell\in H_n)$.
\end{dem}

\begin{cor}
\label{fqcogrowth}
A quotient of a group by a finite normal subgroup has the same gross cogrowth.
\end{cor}

\begin{dem}
Let $H$ be a finite subgroup of $G$ and let $n=\max\{\norm{h}, h\in H\}$
so that $H$ is included in the $n$-neighborhood of $e$.
Then $\Pr(B_\ell =_{G/H} e)=\Pr(B_\ell\in H)\leq \sum_{k\leq n}(2m)^k
\Pr(B_{\ell+k}=e)\lesssim (2m)^{-(1-\theta)\ell}$.
\end{dem}

\begin{rem*}
Gross cogrowth is the same only if defined with respect to the same set of
generators. For example, $F_2\times \Z/2\Z$ presented by $a,b,c$ with
$ac=ca$, $bc=cb$ and $c^2=e$ has the same gross cogrowth as $F_2$ presented by
$a,b,c$ with $c=e$.
\end{rem*}

So in this case, we can safely assume that the virtual centre of $G$ is
trivial.

\subsubsection{The case of random geodesic words}

A quotient by a finite normal subgroup preserves growth, so $G$ and $G'$ have
the same growth.

But now a problem arises, as the notion of a random element of norm
$\ell$ differs in $G$ and $G'$. So our random set $R$ is not defined the
same way for $G$ and $G'$.

Let us study the image of the uniform measure on the $\ell$-sphere of $G$
into $G'$. Let $L$ be the maximal norm of an element in $H$. The image
of this sphere is contained in the spheres of radius between $\ell-L$ and
$\ell+L$.

The map $G\rightarrow G'$ is of index $\abs{H}$. This proves that the
image of the uniform probability measure $\mu_\ell$ on the sphere of
radius $\ell$ in $G$ is, as a measure, at most $\abs{H}$ times the sum of
the uniform probability measures on the spheres of $G'$ of radius between
$\ell-L$ and $\ell+L$.  In other words, it is roughly less than the
uniform probability measure $\nu_\ell$ on these spheres.

The uniform measure $\nu_\ell$ on the spheres of radius between $\ell-L$
and $\ell+L$ (for a fixed $L$) satisfies our axioms.
So we can apply Theorem~\ref{techmain} to the quotient of $G'$ by a set
$R'$ of random words chosen using measure $\nu_\ell$. This random quotient will
be non-elementary hyperbolic for $d<g/2$.

By Remark~\ref{smallerquotient}, for a random set $R$ picked from measure
$\mu_\ell$ (the one we are interested in), the quotient $G'/\langle
R\rangle$ will be non-elementary hyperbolic as well.

But $G'/\langle R \rangle=G/H/\langle R\rangle=G/\langle R \rangle/H$,
and quotienting $G/\langle R \rangle$ by the finite normal subgroup
$H$ is a quasi-isometry, so $G/\langle R \rangle$ is non-elementary
hyperbolic if and only if $G'/\langle R\rangle$ is.

\section{Proof of the main theorem}
\label{mainsection}

We now proceed to the proof of Theorem~\ref{techmain}.

$G$ is a hyperbolic group without virtual centre generated by $S=a_1^{\pm
1},\ldots$, $a_m^{\pm 1}$. Say that $G$ has presentation
$\presgroup{S}{Q}$.  Let $R$ be a set of random words of density at most
$d$ picked under the measure $\mu_\ell$. We will study $G/\langle R
\rangle$.

Let $\beta=\min(\beta_2,\beta_3,\beta_4)$
where $\beta_2,\beta_3,\beta_4$ are given by the axioms. We assume
that $d<\beta$.

\bigskip

We will study van Kampen diagrams in the group $G/\langle R \rangle$.  If
$G$ is presented by $\presgroup{S}{Q}$, call \emph{old relator} an
element of $Q$ and \emph{new relator} an element of $R$.

We want to show that van Kampen diagrams of $G/\langle R \rangle$ satisfy a linear
isoperimetric inequality. Let $D$ be such a diagram. $D$ is made of old
and new relators. Denote by $D'$ the subdiagram of $D$ made of old
relators and by $D''$ the subdiagram of $D$ made of new relators.

If $\beta=0$ there is nothing to prove. Hence we suppose that $\beta>0$.
In the examples we consider, this is equivalent to $G$ being non-elementary.

\subsection{On the lengths of the relators}

In order not to make the already complex notations even heavier, \emph{we will
suppose that all the words taken from $\mu_\ell$ are of length
$\ell$}. So $R$ is made of $(2m)^{d\ell}$ words of length $\ell$.
This is the case in all the applications given in this text.

For the general case, there are only three ways in which the length of the
elements matters for the proof:
\begin{enumerate}
\item As we are to apply asymptotic estimates, the length of the elements
must tend to infinity.
\item The hyperbolic local-global theorem of Appendix~\ref{CHGP}
crucially needs that the ratio of the lengths of relators be bounded independently of $\ell$.
\item In order not to perturb our probability estimates, the number of
distinct lengths of the relators in $R$ must be subexponential in $\ell$.
\end{enumerate}

All these properties are guaranteed by Axiom~1.

\subsection{Combinatorics of van Kampen diagrams of the quotient}
\label{combdiag}

We now proceed to the application of the program outlined in
section~\ref{outline}. The reader may want to refer to this section while
reading the sequel of this text.

We consider a van Kampen diagram $D$ of $G/\langle R\rangle$. Let $D'$ be the
part of $D$ made of old relators of the presentation of $G$, and $D''$
the part made of new relators in $R$.

Redefine $D'$ by adding to it all edges of $D''$: this amounts to adding
some filaments to $D'$. This way, we ensure that faces of $D''$ are
isolated and that $D'$ is connected; and that if a face of $D''$ lies on
the boundary of $D$, we have a filament in $D'$, such that $D''$ does
not intersect the boundary of $D$; and last, that if the diagram $D''$ is
not regular (see section~\ref{defs} for definition), we have a
corresponding filament in $D'$.

\begin{center}
\includegraphics{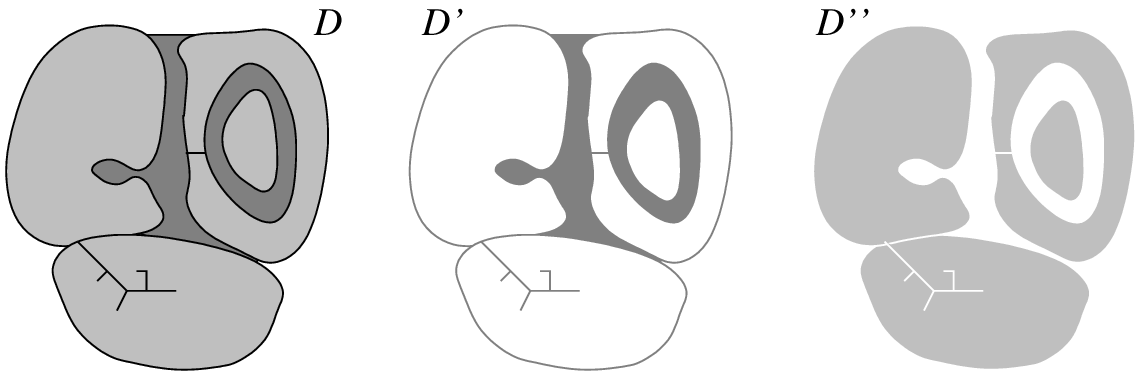}
\end{center}

After this manipulation, we consider that each edge of $D''$ is in
contact only with an edge of $D'$, so that we never have to deal with
equalities between subwords of two new relators (we will treat them as two
equalities to the same word).

\bigskip

We want to show that if $D$ is minimal, then it satisfies some
isoperimetric inequality. In fact, as in the case of random quotients of
a free group, we do not really need that $D$ is minimal. We need that
$D$ is reduced in a slightly stronger sense than previously, which we
define now.

\begin{defi}
A van Kampen diagram $D=D'\cup D''$ on $G/\langle R\rangle$ (with $D'$ and
$D''$ as above) is said to be \emph{strongly
reduced} with respect to $G$ if there is no pair of faces of $D''$ bearing the
same relator with opposite orientations, such that their marked starting
points are joined in $D'$ by a simple path representing the trivial element in
$G$.
\end{defi}

In particular, a strongly reduced diagram is reduced.

\begin{prop}
\label{strongreduction}
Every van Kampen diagram has a strong reduction, that is, there exists a
strongly reduced diagram with the same boundary.
\end{prop}

In particular, to ensure hyperbolicity of a group it is enough to prove
the isoperimetric inequality for all strongly reduced diagrams.

\begin{dem}
Suppose that some new relator
$r$ of $D''$ is joined to some $r\~$ by a path $w$ in $D'$ representing
the trivial element in $G$. Then incise the diagram along $w$ and apply
surgery to cancel $r$ with $r\~$. This leaves a new diagram with two
holes $w, w\~$. Simply fill up these two holes with diagrams in $G$
bordered by $w$ (this is possible precisely since $w$ is the trivial
element of $G$).

\begin{center}
\includegraphics{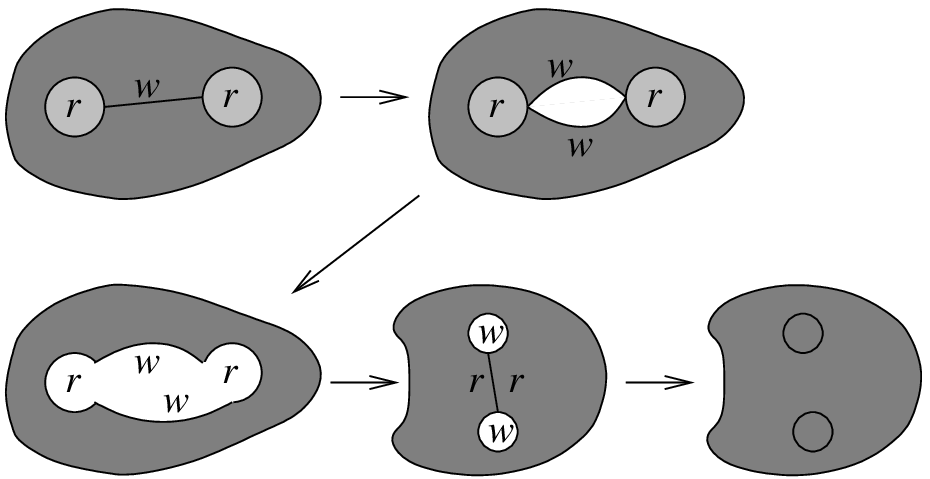}
\end{center}

Note that this way we introduce only old relators and no new ones in the
diagram. Iterate the process to get rid of all annoying pairs of new relators.
\end{dem}

Often in geometric group theory, problems arise when two relators are
conjugate (or when a conjugate of a relator is very close to another
relator), and such cases are typically excluded by reinforcing the
definition of ``strongly reduced''. In the case of random presentations,
however, below the critical density it never occurs that two relators are
conjugated. So we do not have to care about these problems: these cases
are automatically wiped off by our axioms.

\bigskip

We will show that any strongly reduced van Kampen diagram $D$ such that
$D'$ is minimal very probably satisfies some linear isoperimetric
inequality. By the local-global principle for hyperbolic spaces
(Cartan-Hadamard-Gromov theorem, cf.\ Appendix~\ref{CHGP}), it is enough
to show it for diagrams having less than some fixed number of faces. This
crucial point considerably simplifies the geometric and probabilistic
treatment. More precisely, we will show the following.

\begin{prop}
\label{mainprop}
There exist constants $\alpha, \alpha'>0$ (depending on $G$, $d$ and the
random model but not
on $\ell$) such that, for any integer $K$, with probability exponentially close to $1$ as
$\ell\rightarrow \infty$ the set of relators $R$ satisfies the
following:

For any van Kampen diagram $D=D'\cup D''$
satisfying the three conditions:
\begin{itemize}
\item The number of faces of $D''$ is at most $K$;
\item $D'$ is minimal among van Kampen diagrams in $G$ with the same boundary;
\item $D$ is strongly reduced with respect to $G$;
\end{itemize}
then $D$ satisfies the isoperimetric inequality
\[
\abs{\d D}\geq \alpha\ell\abs{D''}+\alpha' \abs{D'}
\]
\end{prop}

(Of course, the constant implied in ``exponentially close'' depends on
$K$.)

\bigskip

Before proceeding to the proof of this proposition, let us see how it
implies hyperbolicity of the group $G/\langle R\rangle$, as well as that
of all intermediate quotients. This step uses the local-global
hyperbolic principle (Appendix~\ref{CHGP}), which essentially states that
it is enough to check the isoperimetric inequality for a finite number of
diagrams.

\begin{prop}
\label{mainprop2}
There exists an integer $K$ (depending on $G$ and $d$ but not on $\ell$) such that if the set of
relators $R$ happens to satisfy the conclusions of
Proposition~\ref{mainprop}, with $\ell$ large enough, then $G/\langle R\rangle$ is hyperbolic.
Better, then there exist constants $\alpha_1, \alpha_2>0$ such
that for any strongly reduced diagram $D$ such that $D'$ is minimal, we have
\[
\abs{\d D}\geq \alpha_1 \ell \abs{D''}+\alpha_2\abs{D'}
\]
\end{prop}

\begin{rem}
Proposition~\ref{mainprop2} implies that a quotient of $G$ by a smaller
set $R'\subset R$ is hyperbolic as well. Indeed, any strongly reduced
diagram on $R'$ is, in particular, a strongly reduced diagram on $R$.
\end{rem}

\begin{dem}
By our strongly reduction process, for any van Kampen diagram there exists
another van Kampen diagram $D$ with the same boundary, such that $D'$ is
minimal (otherwise replace it by a minimal diagram with the same boundary)
and $D$ is strongly reduced. Thus, it is enough to show the isoperimetric
inequality for strongly reduced diagrams to ensure hyperbolicity.

We want to apply Proposition~\ref{localglobal}. Take for property $P$
in this proposition ``to be strongly reduced''. Recall the notations of
the appendix: $L_c(D)=\abs{\d D}$ is the boundary length of $D$, and
$A_c(D)$ is the area of $D$ in the sense that a relator of length $L$ has
area $L^2$. Note that $\ell\abs{D''}+\abs{D'}\geq A_c(D)/\ell$.

Take a van Kampen diagram $D$
such that $k^2/4 \leq A_d(D) \leq
480 k^2$ for some $k^2=K\ell^2$ where $K$ is some constant independent of
$\ell$ to be chosen later. As $A_d(D)\leq K\ell^2$, we have
$\abs{D''}\leq K$. Proposition~\ref{mainprop} for this $K$ tells us that $L_c(D)=\abs{\d D}\geq
\alpha\ell\abs{D''}+\alpha'\abs{D'}\geq \min(\alpha,\alpha') A_c(D)/\ell$.
Thus
\[
L_c(D)^2\geq \min(\alpha,\alpha')^2 A_c(D)^2/\ell^2
\geq \min(\alpha,\alpha')^2 A_c(D) K/4
\]
as $A_c(D)\geq k^2/4$,
so taking $K=10^{15}/\min(\alpha,\alpha')^2$ is enough to ensure that the
conditions of Proposition~\ref{localglobal} are fulfilled by $K\ell^2$. (The
important point is that this $K$ is independent of $\ell$.)

The conclusion is that any strongly reduced van Kampen diagram $D$ satisfies
the linear isoperimetric inequality
\[
L_c(D)\geq A_c(D)\min(\alpha,\alpha')/10^{12}\ell
\]
and, fiddling with the constants and using the isoperimetry from $D$, we
can even put it in the form
\[
\abs{\d D}\geq \alpha_1 \ell \abs{D''}+\alpha_2\abs{D'}
\]
if it pleases,
where $\alpha_{1,2}$ depend on $G$ and $d$ but not on $\ell$.

So the proposition above, combined with the
local-global hyperbolicity principle of Appendix~\ref{CHGP}, is sufficient
to ensure hyperbolicity.
\end{dem}

\bigskip

A glance through the proof can even show that if $\ell$ is taken large
enough, the constant $\alpha_2$ in the inequality
\[
\abs{\d D}\geq \alpha_1 \ell \abs{D''}+\alpha_2\abs{D'}
\]
is arbitrarily close to the original
isoperimetry constant in $G$.

This suggests, in the spirit of~\cite{Gro4}, to iterate the operation of
taking a random quotient, at different lengths $\ell_1$, then $\ell_2$,
etc., with fast growing $\ell_i$. The limit group will not be hyperbolic
(it will be infinitely presented), but it will satisfy an isoperimetric
inequality like
\[
\abs{\d D}\geq \alpha \sum_{f\text{ face of }D} \ell(f)
\]
where $\ell(f)$ denotes the length of a face. This property could be
taken as a definition of a kind of loose hyperbolicity, which should be
related in some way to the notion of ``fractal hyperbolicity'' proposed
in~\cite{Gro4}.

\bigskip

Now for the proof of Proposition~\ref{mainprop}.

We have to assume that $D'$ is minimal, otherwise we know nothing about
its isoperimetry in $G$. But as in the case of a random quotient of $F_m$
(section~\ref{standardcase}), the isoperimetric inequality will not only
be valid for minimal diagrams but for all (strongly reduced)
configurations of the random relators.

If $D''=\varnothing$ then $D=D'$ is a van Kampen diagram of $G$ and as $D'$ is
minimal, it satifies the inequality $\abs{\d D}\geq C\abs{D}$ as this is
the isoperimetric inequality in $G$. So we can
take $\alpha'=C$ and any $\alpha$ in this case. Similarly, if the old
relators are much more numerous that the new ones, then isoperimetry of
$G$ is enough. Namely:

\begin{lem}
\label{lemsmallold}
Proposition~\ref{mainprop} holds for
diagrams satisfying $\abs{D'}\geq 4\abs{D''}\ell/C$.
\end{lem}

\begin{dem}[ of the lemma]
Suppose that the old relators are much more numerous than the new ones,
more precisely that $\abs{D'}\geq 4\abs{D''}\ell/C$. In this case as
well, isoperimetry in $G$ is enough to ensure isoperimetry of $D$. Note
that $D'$ is a diagram with at most $\abs{D''}$ holes. We have of course
that $\abs{\d D}\geq \abs{\d D'}-\abs{\d D''}\geq \abs{\d
D'}-\abs{D''}\ell$.

By Proposition~\ref{isohole} for diagrams with holes in $G$, we have that
$\abs{\d D'}\geq C\abs{D'}-\abs{D''}\lambda (2+4\alpha \log \abs{D'})$. So, for
$\ell$ big enough,
\begin{eqnarray*}
\abs{\d D} &\geq & \abs{\d D'}-\abs{D''}\ell
\\&\geq&
C\abs{D'}-\abs{D''}\ell-\abs{D''}\lambda (2+4\alpha \log \abs{D'})
\\&\geq&
C\abs{D'}/3+\left(C\abs{D'}/3-\abs{D''}\ell\right)
\\& &+\left(C\abs{D'}/3-\abs{D''}\lambda (2+4\alpha \log \abs{D'})\right)
\\&\geq&
C\abs{D'}/3+(4\abs{D''}\ell/3-\abs{D''}\ell)
\\& &+\left(4\abs{D''}\ell/3-\abs{D''}\lambda (2+4\alpha \log 4\abs{D''}\ell/C)\right)
\\&\geq&
C\abs{D'}/3+\ell\abs{D''}/3 
\end{eqnarray*}
as for $\ell$ big enough, the third term is positive. So in this case we
can take $\alpha=1/3$ and $\alpha'=C/3$.
\end{dem}

So \emph{we now suppose that $1\leq\abs{D''}\leq K$ and that $\abs{D'}\leq
4\abs{D''}\ell/C$}. In particular, the boundary length of $D$ is at most
$\abs{D''}\ell+\abs{D'}\lambda\leq \ell\abs{D''}(1+4\lambda/C)$.

\subsection{New decorated abstract van Kampen
diagrams}
\label{newdavKds}

We now redefine decorated abstract van Kampen diagrams so that they
better fit our needs (the definition given in the introduction fits the
case of free groups only).  The idea is that since $D'$ is very narrow
(at the scale of $\ell$), at scale $\ell$ $D$ looks like a van Kampen
diagram with respect to the new relators, with some narrow ``glue'' (that
is, old relators) between faces. This intuition will be formalized
using Proposition~\ref{coarsening} in Appendix~\ref{appiso}, which will
help tell
which parts of the boundary words of the new relators are facing
which.

The diagram $D'$ has at most $K$ holes.  First, after
Corollary~\ref{narrowhole}, we can suppose that $D'$ is
$E_1\log\ell$-narrow for some constant $E_1$ depending on $G$ and $K$ but
not on $\ell$ (here we used $\abs{D'}\leq 4K\ell/C$ to get logarithmic
dependence on $\ell$).

Besides, we can apply Proposition~\ref{coarsening} to $D'$. This defines
a $(8K,E_2\log\ell)$-matching $X$ (see Definition~\ref{defimatching})
between at most $8K$ subwords of the boundary words of $D'$, for some
constant $E_2$ depending on $G$ and $K$ but not on $\ell$ (here again we
used $\abs{D'}\leq 4K\ell/C$ to get logarithmic dependence on $\ell$).
Set $E=\max(E_1,E_2)$.

The boundary words of $D'$ are precisely the new relators on one side,
and the boundary word of $D$ on the other side.

\begin{center}
\includegraphics{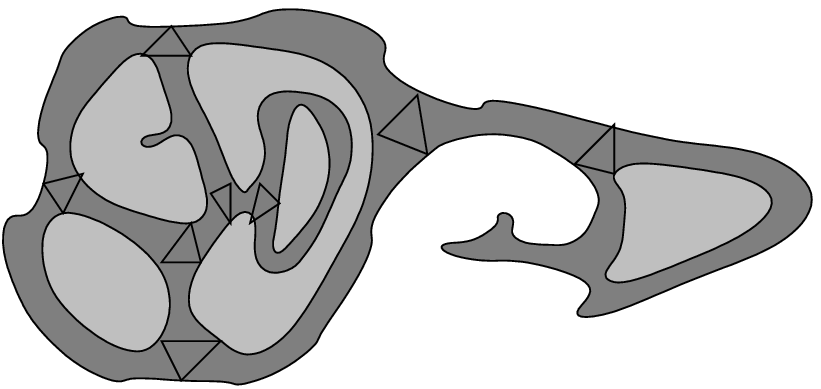}
\end{center}

Each match in $X$ is a pair of two subwords $w$, $w'$ of one of the new
relators (or of the boundary word), together with two short words $u$, $v$
of length at most $E\log\ell$, such that $w=uw'v$ in $G$.

\begin{center}
\includegraphics{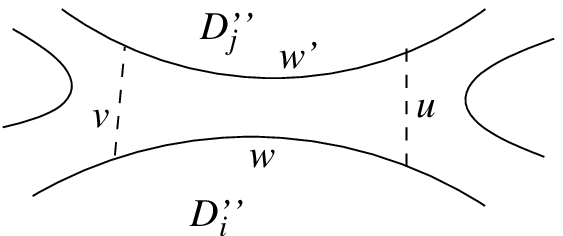}
\end{center}

As there can be ``invaginations'' of $D'$ into $D''$, the
lengths of $w$ and $w'$ may not be equal at all. It may even be the case
that one of these two words is of length $0$, as in the following
picture. This is not overmuch disturbing but should be kept in mind.

\begin{center}
\includegraphics{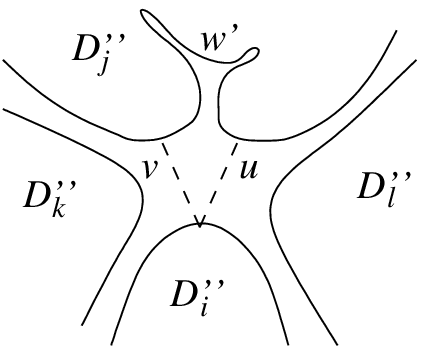}
\end{center}

\bigskip

Intuitively, we can reconstruct $D$ ``at scale $\ell$'' if we know this
matching $X$: simply take the new relators and glue them along the
matches in $X$. This leads to redefining what a decorated abstract van Kampen
diagram is. Knowing $D$, the associated abstract diagram $\mathcal{D}$
will keep the combinatorial and geometric information but will forget what are
the precise values of the new relators. Namely, given $D$ we only keep
the following information: How many new faces there are (that is,
$\abs{D''}$, which is at most $K$); Which new faces bear the same new
relator or not (this can be done by attributing a number between $1$ and
$\abs{D''}$ to each face, two faces getting the same number if and only
if they bear the same new relator); Which subword is matched to which one
(that is, where the cuttings of the subwords were done and what the
pairing is). This leads to the following (compare the definition of a
davKd given in the introduction page~\pageref{firstdavKd}, together with
Proposition~\ref{coarsening}).

\begin{defi}
A \emph{decorated abstract van Kampen diagram} (davKd for short)
$\mathcal{D}$ is the
following data:
\begin{itemize}
\item An integer $k$ (the number of faces), also denoted
$\abs{\mathcal{D}}$; any number between $1$ and
$k$ will be called a \emph{face} of $\mathcal{D}$.
\item An integer $\abs{\d\mathcal{D}}$ between $0$ and
$\abs{\mathcal{D}}\ell(1+4\lambda/C)$, called the \emph{boundary length}
of $\mathcal{D}$.
\item A set of $k$ integers between $1$ and $k$ (which faces bear the
same relator).
\item For each face, a number between $1$ and $\ell$ (a starting point
for the relator) and an orientation $\pm 1$.
\item A partition of the set $\{1,\ldots,\ell\}\times
\{1,\ldots,k\}\cup \{1,\ldots,\abs{\d\mathcal{D}}\}\times \{k+1\}$ into $8k$ subsets (some of which may be empty) of the form
$\{i,i+1,\ldots,i+j\}\times\{p\}$. The subsets
$\{1,\ldots,\ell\}\times\{p\}$ will be called \emph{words} in
$\mathcal{D}$, with $\{1,\ldots,\abs{\d\mathcal{D}}\}\times\{k+1\}$ being the
\emph{boundary word} and the others \emph{internal words}. The
elements of the partition will be called \emph{subwords} in
$\mathcal{D}$, and the \emph{length} of a subword
$\{i,i+1,\ldots,i+j\}\times\{p\}$ will be $j+1$. 
\item A partition of the set of subwords into two parts and a bijection
between these parts (which subword is matched to which). A pair of two bijected
subwords will be called a \emph{match} in $\mathcal{D}$.
\end{itemize}
\end{defi}

A very important fact is the following one.

\begin{prop}
\label{numberofdavKds}
For a fixed $K$, the number of different davKd's with at most $K$ faces
is less than some polynomial in $\ell$.
\end{prop}

\begin{dem}
This is at most $K.K\ell(1+4\lambda/C).K^K.\ell^K.2^K.((K+1)\ell)^{8K}.(8K)^{8K}$.
\end{dem}

We will still add some decoration below in section~\ref{themainargument}.
This further decoration will again be polynomial in $\ell$.

We just saw that to our van Kampen diagram $D$ we can associate a
decorated abstract van Kampen diagram $\mathcal{D}$, coming from
Proposition~\ref{coarsening}. This we will call the \emph{davKd
associated to $D$}.  This davKd sums up all quasi-equalities in $G$
imposed by the diagram on the new relators.

%
%
%

\begin{defi}
\label{isodavKd}
Let $\mathcal{D}$ be a davKd. We say that a van Kampen diagram $D=D'\cup
D''$ of
$G/\langle R\rangle$ \emph{fulfills} $\mathcal{D}$ if
$\mathcal{D}$ is the davKd associated to $D$ and if $D$ 
satisfies the
assumptions of Proposition~\ref{mainprop} and Lemma~\ref{lemsmallold}
that is:
\begin{itemize}
\item The number of faces of $D''$ is at most $K$;
\item $D'$ is minimal among van Kampen diagrams in $G$ with the same
boundary;
\item $D$ is strongly reduced with respect to $G$;
\item $\abs{D'}\leq 4\abs{D''}\ell/C$.
\end{itemize}

A davKd $\mathcal{D}$ is said to be \emph{fulfillable} if some van Kampen
diagram fulfills it.

A davKd $\mathcal{D}$ is said to satisfy an
$\alpha$-isoperimetric inequality if
\[
\abs{\d \mathcal{D}}\geq \alpha \ell \abs{\mathcal{D}}
\]
\end{defi}

\begin{prop}
\label{davKdtovKd}
If some davKd $\mathcal{D}$ satisfies an $\alpha$-isoperimetric
inequality, then any van Kampen diagram $D=D'\cup D''$ fulfilling
$\mathcal{D}$
satisfies the isoperimetric inequality
\[
\abs{\d D}\geq \alpha\ell\abs{D''}/2+C\alpha \abs{D'}/8
\]
\end{prop}

\begin{dem}
Indeed, since $\abs{D'}\leq 4\abs{D''}\ell/C$ we have
$\alpha\ell\abs{D''}/2+C\alpha\abs{D'}/8 \leq\alpha\ell\abs{D''}$.
\end{dem}

Thus, to prove Proposition~\ref{mainprop} we have to show that, with high
probability, any fulfillable davKd satisfies some linear isoperimetric
inequality (with some isoperimetric constants depending on $G$, the
density $d$ and the random model but not on $K$ or $\ell$).

\medskip

In the matching $X$ of $D$, there may be matches between subwords of
the boundary, as in the following figure. Such parts of the
diagram always improve isoperimetry (up to $2E\log\ell$). So \emph{in the
following we consider that all matches in $X$ match
either two subwords of the new relators or a subword of a new relator and
a subword of the boundary}.

\begin{center}
\includegraphics{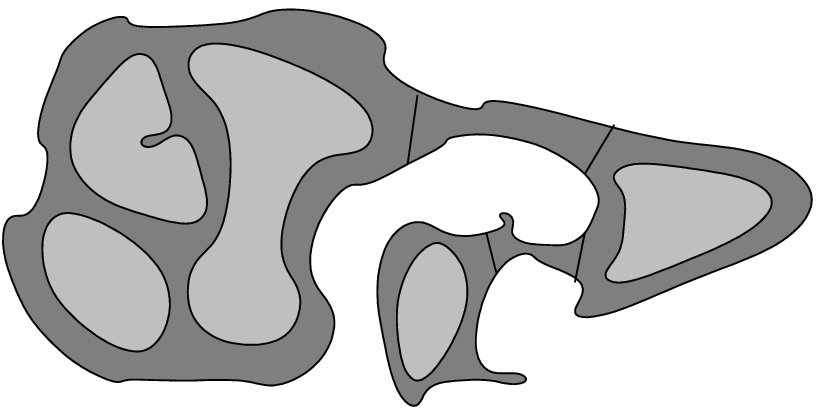}
\end{center}

\subsection{Graph associated to a decorated abstract van Kampen diagram}

As in the case of random quotients of the free group, we will construct
an auxiliary graph $\Gamma$ summarizing all conditions imposed by a davKd
on the random relators of $R$. But instead of imposing equality between
letters of these relators, the conditions will rather be interpreted as
equality modulo $G$.

Let now $D$ be a davKd. We will evaluate the probability that it is
fulfillable by the relators of $R$.

Each face of $D$ bears a number between $1$ and $\abs{D}$. Let $n$ be the
number of such distinct numbers, we have $n\leq \abs{D}$. Suppose for the
sake of simplicity that these $n$ distinct numbers are $1,2,\ldots,n$.

To fulfill the diagram is to give $n$ relators $r_1,\ldots,r_n$
satisfying the conditions that if we put these relators in the
corresponding faces, then by gluing ``up to small words'' the faces along
the subwords desribed in the davKd, we get a (strongly reduced) van
Kampen diagram of $G/\langle R \rangle$.

We now construct the auxiliary graph $\Gamma$.

Take $n\ell$ points as vertices of $\Gamma$, arranged in $n$ parts of
$\ell$ vertices called the \emph{parts} of $\Gamma$.  Interpret the
$k$-th vertex of the $i$-th part as the $k$-th letter of relator $r_i$ in
$R$. Internal subwords in $D$ are identified with successive vertices
of $\Gamma$ (with a reversal if the face in $D$ to which the
subword belongs is negatively oriented).

We now explain what to take as edges of $\Gamma$.

Let $f$ be a match in $D$. First, suppose that this is a match between
two internal subwords in $D$. Say it is a match between subwords of faces
of $D$ bearing numbers $i$ and $i'$.  These two subwords in $D$
correspond to two sets of successive vertices in the $i$-th part and the
$i'$-th part of $\Gamma$.

Add to $\Gamma$ a special vertex $w$ called an \emph{internal
translator}. Add edges between $w$ and each of the vertices of the $i$-th
part of $\Gamma$ represented by the first subword of $f$; symmetrically,
add edges between $w$ and each of the vertices of the $i'$-th part of
$\Gamma$ belonging to the other subword of $f$.

(This may result in double edges if $i=i'$. We will deal with this
problem later, but for the moment we keep the double edges.)

Follow this process for all matches between internal subwords of $D$.
Each translator so obtained is connected with two (or maybe one if
$i=i'$) parts of $\Gamma$.

As several faces of $D$ may bear the same number (the same relator of $R$), a
vertex of $\Gamma$ is not necessarily of multiplicity one. The
multiplicity of a vertex of the $i$-th part is at most the number of
times relator $i$ appears on a face of $D$.

For each match in $D$ involving a boundary subword and an internal
subword adjacent to, say, face $i$,
add a special vertex $b$ to $\Gamma$, called a \emph{boundary
translator}. Add edges between $b$ and the vertices of the $i$-th
part of $\Gamma$ corresponding to the internal subword of the match at play.

At the end of the process, the number of edges in $\Gamma$ is equal to
the cumulated length of all internal subwords of $D$, which is
$\ell\abs{D}$.

Here is an example of a simple van Kampen diagram on $G/\langle
R\rangle$, its associated davKd (represented graphically by a diagram
``at scale $\ell$''), and the associated graph $\Gamma$.

\label{gammaexample}
\begin{center}
\includegraphics{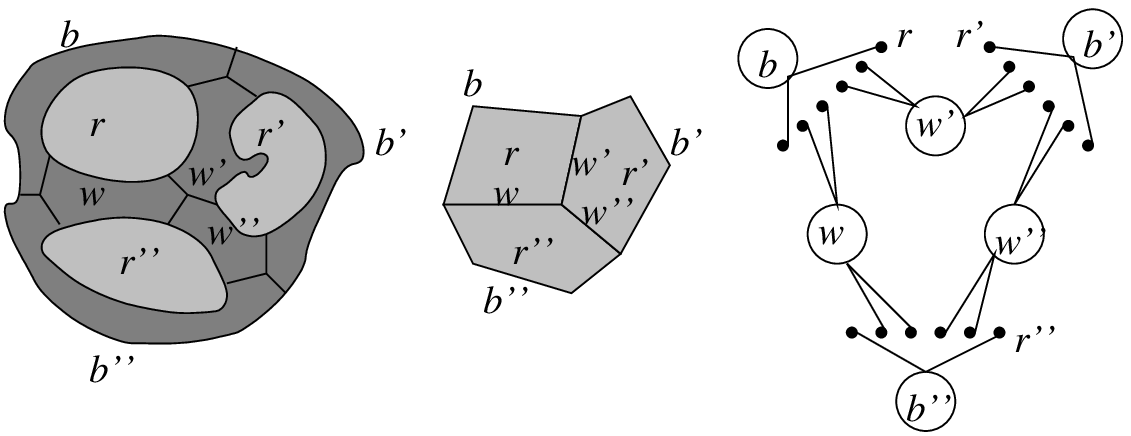}
\end{center}

In a davKd associated to some van Kampen diagram with at most $K$ faces,
as we only consider, since the number of matches is at most $4K$, the
number of internal and boundary translators in $\Gamma$ is at most $4K$
as well.

Note that each translator corresponds, via Proposition~\ref{coarsening},
to a word in the van Kampen diagram which is equal to $e$ in $G$:
translators are nothing else but matches of the davKd.  Indeed,
fulfillability of the davKd implies that for each (say internal)
translator in $\Gamma$, we can find a word $w$ which is equal to $e$ in
$G$, and such that $w=w_1uw_2v$ where $u$ and $v$ are short (of length at
most $E\log\ell$) and that $w_1$ and $w_2$ are the subwords of the
relators of $R$ to which the translator is joined. In the case of random
quotients of $F_m$, we had the relators of $R$ directly connected to each
other, imposing equality of the corresponding subwords; here this
equality happens modulo translators that are equal to $e$ in $G$.

\subsection{Elimination of doublets}

A \emph{doublet} is a vertex of $\Gamma$ that is joined to some
translator by a double edge. This can occur only if in the van Kampen
diagram, two nearly adjacent faces bear the same relator.

Doublets are annoying since the two sides of the translator are not
chosen independently, whereas our argument requires some degree of
independence. We will split the corresponding translators to control the
occurrences of such a situation.

This section is only technical.

\bigskip

Consider a translator in the van Kampen diagram bordered by
two faces bearing the same relator $r$. As a first case, suppose that
these two relators are given the same orientation.

Let $w$ be the translator, $w$ writes $w=u\delta_1 u'\delta_2$ where $u$
and $u'$ are subwords of $r$, and $\delta_{1,2}$ are words of length at
most $2E\log \ell$. The action takes place in $G$. As $u$ and $u'$ need
not be geodesic, they do not necessarily have the same length. Let $u_1$
be the maximum common subword of $u$ and $u'$ (i.e.\ their intersection as subwords of $r$). If
$u_1$ is empty there is no doublet.

There are two cases (up to exchanging $u$ and $u'$): either $u=u_2u_1u_3$
and $u'=u_1$, or $u=u_2u_1$ and $u'=u_1u_3$.

\begin{center}
\includegraphics{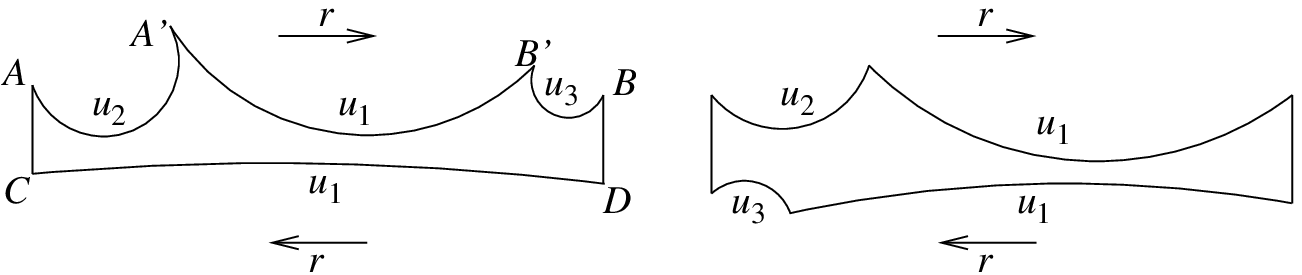}
\end{center}

We will only treat the first case, as the second one is similar.

Redefine $u_1$, $u_2$ and $u_3$ to be geodesic words equal to $u_1$,
$u_2$ and $u_3$ respectively. In any hyperbolic space, any point on a geodesic
joining the two ends of a curve of length $L$ is $(1+\delta\log
L)$-close to that curve (cf.~\cite{BH}, p.~400). So the new geodesic
words are $(1+\delta\log\ell)$-close to the previous words $u_1$, $u_2$,
$u_3$. Hence, up to increasing $E$ a little bit, we can still suppose
that $D$ is fulfillable such that $D'$ is $E\log\ell$-narrow, and that
$u_1$, $u_2$, $u_3$ are geodesic.

Define points $A, A', B, B', C, D$ as in the figure. The word read while
going from $A'$ to $B'$ is the same as that from $D$ to $C$.

By elementary hyperbolic geometry, and given that the two lateral sides
are of length at most $2E\log \ell$, any point on $CD$ is
$(2\delta+2E\log\ell)$-close to some point on $AA'$ or $B'B$, or
$2\delta$-close to some point on $A'B'$.

The idea is to run from $D$ to $C$, and simultaneously from $A'$ to $B'$
at the same speed. When the two trajectories get $E\log\ell$-close to
each other, we cut the translator at this position, and by construction
the resulting two parts do not contain any doublets.

Let $L=\abs{u_1}$ and for $0\leq i \leq L$, let $C_i$ be the point of
$DC$ at distance $i$ from $D$. Now assign to $i$ a number $\phi(i)$
between $0$ and $L$ as follows: $C_i$ is close to some point $C'_i$ of
$AB$, set $\phi(C_i)=0$ if $C'_i\in AA'$, $\phi(C_i)=L$ if $C'_i\in B'B$,
and $\phi(C_i)=\dist(C'_i,A')$ if $C'_i\in A'B'$.

By elementary hyperbolic geometry (approximation of $A'B'DC$ by a tree),
the function $\phi:[0;L]\rightarrow[0;L]$ is decreasing up to $8\delta$
(that is, $i<j$ implies $\phi(i)>\phi(j)-8\delta$). We have $\phi(0)=L$
and $\phi(L)=0$ (up to $8\delta$). Set $i_0$ as the smallest $i$ such
that $\phi(i)<i$. This defines a point $C_{i_0}$ on $DC$ and a point
$C'_{i_0}$ on $AB$.

There are six cases depending on whether $C'_{i_0}$ and $C'_{i_0-1}$
belong to $AA'$, $A'B'$ or $B'B$. In each of these cases we can cut the
diagram in at most three parts, in such a way that no part contains two
copies of some subword of $u_1$ (except perhaps up to small words of
length at most $8\delta$ at the extremities).  The cuts to make are from
$C_{i_0}$ to $C'_{i_0}$ and/or to $C'_{i_0-1}$, and are
illustrated below in each case.

\begin{center}
\includegraphics{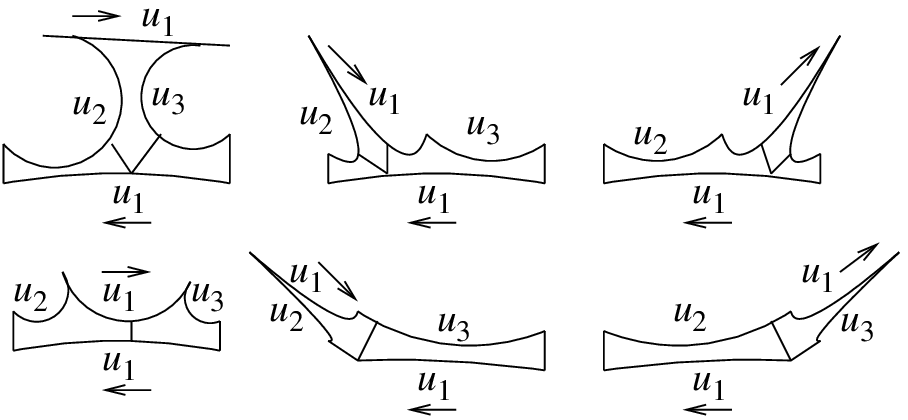}
\end{center}

A translator is a vertex of $\Gamma$ and by ``cutting a
translator'' we mean that we split this vertex into two, and share the
edges according to the figure.

\bigskip

As our second (and more difficult) case, suppose that the translator is
bordered by two faces of the diagram bearing the same relator $r$ of $R$
with opposite orientations. This means that the translator $w$ is equal,
in $G$, to $u\delta_1{u'}\~\delta_2$ where $u$ and $u'$ are subwords of
the relator $r$, and where $\delta_{1,2}$ are words of length at most
$2E\log\ell$.

As above, let $u_1$ be the maximum common subword of $u$ and $u'$ (i.e.\
their intersection as subwords of $r$). There are two cases:
$u=u_2u_1u_3$ and $u'=u_1$, or $u=u_2u_1$ and $u'=u_1u_3$.

\begin{center}
\includegraphics{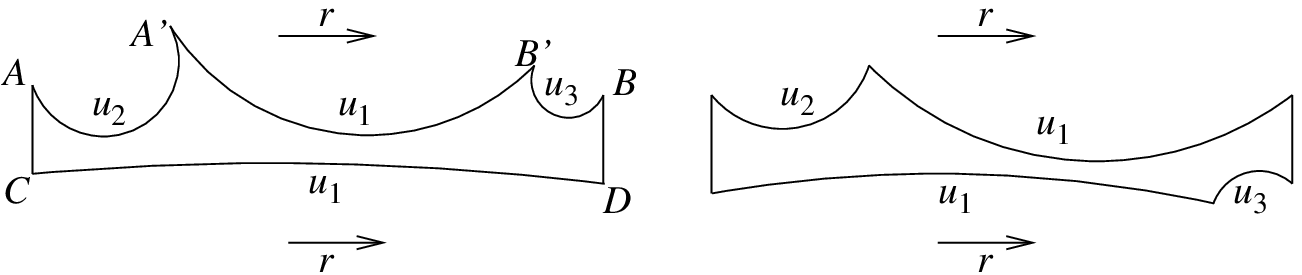}
\end{center}

We will only treat the first case, as the second is similar.

As above, redefine $u_1, u_2$ and $u_3$ to be geodesic.

Define points $A, A', B, B', C, D$ as in the figure. The word read while
going from $A'$ to $B'$ is the same as that from $C$ to $D$.

By elementary hyperbolic geometry, and given that the two lateral sides
are of length at most $2E\log \ell$, any point on $CD$ is
$(2\delta+2E\log\ell)$-close to some point on $AA'$ or $B'B$, or
$2\delta$-close to some point on $A'B'$.

If any point on $CD$ is close to a point on either $AA'$ or $BB'$, we can
simply eliminate the doublets by cutting the figure at the last point of
$CD$ which is close to $AA'$. (As above, by cutting the figure we mean
that we split the vertex of $\Gamma$ representing the translator into
three new vertices and we share its edges according to the figure.) In
this way, we obtain a new graph $\Gamma$ with the considered doublets
removed.

\begin{center}
\includegraphics{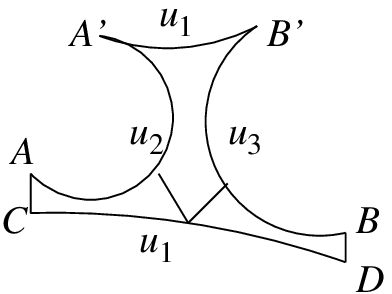}
\end{center}

Otherwise, let $L=\abs{u_1}$ and for $0\leq i \leq L$, let $C_i$ be the
point of $CD$ at distance $i$ from $C$. Now assign to $i$ a number
$\phi(i)$ between $0$ and $L$ as follows: $C_i$ is close to some point
$C'_i$ of $AB$, set $\phi(C_i)=0$ if $C'_i\in AA'$, $\phi(C_i)=L$ if
$C'_i\in B'B$, and $\phi(C_i)=\dist(C'_i,A')$ if $C'_i\in A'B'$.

It follows from elementary hyperbolic geometry (approximation of the
qua\-drilateral $CA'B'D$ by a tree) that $\phi:[0;L]\rightarrow[0;L]$ is an
increasing function up to $8\delta$ (that is, $i<j$ implies
$\phi(i)<\phi(j)+8\delta$). Moreover, let $i$ be the smallest such that
$\phi(i)>0$ and $j$ the largest such that $\phi(j)<L$. Then $\phi$ is, up
to $8\delta$, an isometry of $[i;j]$ to $[\phi(i);\phi(j)]$ (this is
clear on the approximation of $CA'B'D$ by a tree). In other words: the
word $u_1$ is close to a copy of it with some shift $\phi(i)-i$.

Cut the figure in five: cut between $C_i$ and $C'_i$, between $C_i$
and a point of $AA'$ close to it, between $C_j$ and $C'_j$ and between
$C_j$ and a point of $B'B$ close to it (such points exist by definition
of $i$ and $j$).

\begin{center}
\includegraphics{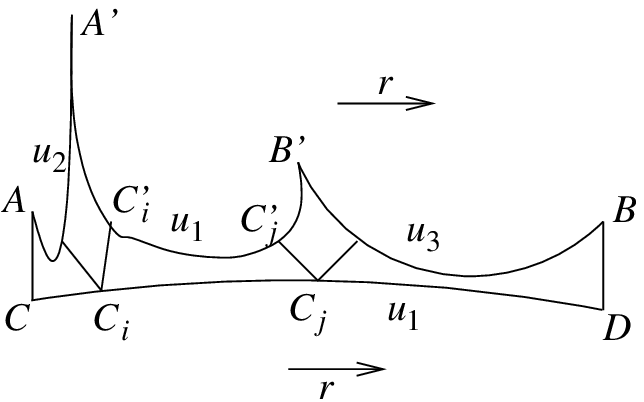}
\end{center}

This way, we get a figure in which only the middle part $C_iC_jC'_jC'_i$
of the figure contains two copies of a given piece of $u_1$. Indeed (from
left to right in the figure) the first part contains letters $0$ to $i$
of the lower copy of $u_1$ and no letter of the upper $u_1$; the second
part contains letters $0$ to $\phi(i)$ of the upper $u_1$ and no letter
of the lower $u_1$; the third part $C_iC_jC'_jC'_i$ contains letters $i$
to $j$ of the lower $u_1$ and letters $\phi(i)$ to $\phi(j)$ of the upper
$u_1$; the fourth and fifth part each contain letters from only one copy
of $u_1$.

First suppose that the intersection of $[i;j]$ and $[\phi(i);\phi(j)]$ is
empty, or that its size is smaller than $\eps_1\abs{u_1}$ (for some small
$\eps_1$ to be fixed later on, depending on $d$ and $G$ but not on
$\ell$). Then, in the new graph $\Gamma$ defined by such cutting of the
translator, at most $\eps_1\abs{u_1}$ of the doublets at play remain.
Simply remove these remaining double edges from the graph $\Gamma$. 

In case the intersection of $[i;j]$ and $[\phi(i);\phi(j)]$ is not
smaller than $\eps_1\abs{u_1}$, let us now deal with the middle piece.

Consider the subdiagram $C_iC_jC'_jC'_i$: it
is bordered by two subwords $u'_1, u''_1$ of $u_1$ of non-empty
intersection. The subword $u'_1$ spans letters $i$ to $j$ of $u_1$,
whereas $u''_1$ spans letters $\phi(i)$ to $\phi(j)$, with
$\phi(j)-\phi(i)=j-i$ up to $8\delta$.

First suppose that the shift $\phi(i)-i$ is bigger than $\eps_2\abs{u_1}$.
Then, chop the figure into sections of size $\eps_2\abs{u_1}$:

\begin{center}
\includegraphics{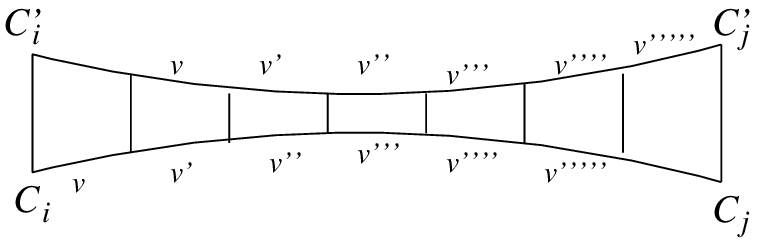}
\end{center}

The word read on one side of a section is equal to the word read on the
other side of the following section, but there are no more doublets. The
original translator has been cut into at most $1/\eps_2$ translators, the
length of each of which is at least $\eps_2\abs{u_1}$.

Second (and last!), suppose that the shift $\phi(i)-i$ is smaller than
$\eps_2\abs{u_1}$. This means that we have an equality $w_1vw_2v\~$ in $G$,
where $v$ is a subword of a random relator $r$, of length at least
$\eps_1\abs{u_1}$, and with $w_1, w_2$ words of length at most
$\eps_2\abs{u_1}$.

As the diagram is strongly reduced, $w_1$ and $w_2$ are non-trivial in
$G$. As the virtual centre of $G$ has been supposed to be trivial, the
probability of this situation is controlled by Axiom~4. Let this
translator as is, but mark it (add some decoration to $\Gamma$) as being
a \emph{commutation translator}.  Furthermore, remove from this
translator all edges that are not double edges, that is, all edges not
corresponding to letters of the $v$ above (there are at most
$2\eps_2\abs{u_1}$ of them).

\bigskip

Follow this process for each translator having doublets. After this, some
doublets have been removed, and some have been marked as being part of a
commutation translator. Note that we suppressed some of the edges of
$\Gamma$, but the proportion of suppressed edges is less than
$\eps_1+2\eps_2$ in each translator.

\subsection{Pause}

Let us sum up the work done so far. Remember the example on
page~\pageref{gammaexample}.

\begin{prop}
For each strongly reduced van Kampen diagram $D$ of the quotient
$G/\langle R\rangle$ such that $\abs{D''}\leq K$ and $\abs{D'}\leq
4\abs{D''}\ell/C$, we have constructed a
graph $\Gamma$ enjoying the following properties:
\end{prop}

\begin{itemize}
\item Vertices of $\Gamma$ are of four types: ordinary vertices, internal
translators, boundary translators, and commutation translators.
\item There are $n\ell$ ordinary vertices of $\Gamma$, grouped in $n$
so-called parts, of $\ell$ vertices each, where $n$ is the number of different relators of $R$ that are present in $D$. Hence each ordinary vertex of $\Gamma$ corresponds to some letter of a relator of $R$. 
\item The edges of $\Gamma$ are between translators and ordinary
vertices.
\item The number of edges at any ordinary vertex is at most equal to the
number of times the corresponding relator of $R$ appears in $D$.
\item For each internal translator $t$, the edges at $t$ are consecutive
vertices of one or two parts of $\Gamma$, representing subwords $u$ and
$v$ of relators of $R$. And there exists a word $w$ such that
$w=\delta_1u\delta_2v$ and $w=e$ in $G$, where $\delta_{1,2}$ have length at most
$2E\log\ell$.
\item For each boundary translator $b$, the edges at $b$ are consecutive
vertices of one part of $\Gamma$, representing a subword $u$ of some relator
of $R$. For each such $b$, there exists a word $w$ such that
$w=\delta_1u\delta_2v$ and $w=e$ in $G$, where $v$ is a subword of the boundary of
$D$, and where $\delta_{1,2}$ have length at most
$2E\log\ell$.
\item For each commutation translator $c$, the edges at $c$ are double
edges to successive vertices of one part of $\Gamma$, representing a
subword $u$ of some relator of $R$. And there exists a word $w$ such that
$w=\delta_1u\delta_2u\~$ and $w=e$ in $G$, where $\delta_{1,2}$ have length at most
$\eps_2\abs{u}$.
\item There are no double edges except those at commutation translators.
\item There are at most $4K/\eps_2$ translators.
\item The total number of edges of $\Gamma$ is at least
$\abs{D''}\ell(1-\eps_1-2\eps_2)$.
\end{itemize}

The numbers $K$ and $\eps_1,\eps_2$ are arbitrary. The number $E$ depends
on $G$ and $K$ but not on $\ell$.

Axioms~2, 3 and 4 are carefully designed to control the probability that,
respectively, a boundary translator, internal translator, and commutation
translator can be filled.

Note that this graph depends only on the davKd associated to the van
Kampen diagram (up to some dividing done for the elimination of doublets;
say we add some decoration to the davKd indicating how this was
done).

Keep all these properties (and notations) in mind for the sequel.

\subsection{Apparent length}

The line of the main argument below is to fulfill the davKd by filling the
translators one by one.

As the same subword of a relator can be joined to a large number of
different translators (if the relator appears several times in the
diagram), during the construction, at some steps it may happen that one
half of a given translator is filled, whereas another part is not. The
solution is to remember in one way or another, for each half-filled
translator, what is the probability that, given its already-filled side,
the word on the other side will fulfill the translator. This leads to the
notion of apparent length, which we define now.

\bigskip

Say we are given an element $x$ of the group, of norm $\norm{x}$. We try to
know if this is a subword of one of our random words under the
probability measure $\mu_\ell$, and to determine the length of this subword.

Given Axiom~2, a good guess for the length of the subword would be
$\norm{x}/\kappa_2$, with the probability of a longer subword decreasing
exponentially.

Given Axiom~3, a good method would be to take another subword $y$ of
length $\abs{y}$ at random under $\mu_\ell$, and test (taking $u=v=e$ in
Axiom~3) the probability that $xy=1$. If $x$ were a subword under
$\mu_\ell$, this probability would be roughly $(2m)^{-\beta(\abs{x}+\abs{y})}$,
hence an evaluation $-\frac{1}{\beta}\log\Pr(xy=e)-\abs{y}$
for the hypothetical length of the subword $x$.

\bigskip

This leads to the notion of apparent length.

We are to apply Axiom~3 to translators, with $u$ and $v$ of size
$2E\log\ell$. For fixed $x\in G$, let $L\geq 0$ and denote by $p_L(xuyv=e)$ the
probability that, if $y$ is a subword of length $L$ under $\mu_\ell$ (in
the sense of Definition~\ref{subword}) there exist words $u$ and $v$ of
length at most $2E\log\ell$ such that $xuyv=e$.

\newcommand{\al}{\mathbb{L}}

\begin{defi}[ (Apparent length at a test-length)]
\label{defal}
The \emph{apparent length} of $x$ at test-length $L$ is
\[
\al_L(x)=-\frac{1}{\beta}\log p_L(xuyv=e)-L
\]
\end{defi}

By definition, if we have a rough evaluation of $p_L$, we get an evaluation of
$\al_L$ up to $o(\ell)$ terms.

We are to apply this definition for $y$ a not too small subword. That is,
we will have $\eps_3\ell/\kappa_1\leq \abs{y}\leq \kappa_1\ell$ with
$\kappa_1$ as in Axiom~1, for some $\eps_3$ to be fixed soon. We will
also use the evaluation from Axiom~2.

\begin{defi}[ (Apparent length)]
\label{defal2}
The \emph{apparent length} of $x$ is
\[
\al(x)=\min \left( \norm{x}/\kappa_2, \min_{\eps_3 \ell/\kappa_1\leq L
\leq \kappa_1 \ell}
\al_L(x)\right)
\]
\end{defi}

\bigskip

Our main tool will now be the following

\begin{prop}
\label{synthax}
For any subword $x$ under $\mu_\ell$, we have
\[
\Pr\left(\al(x)\leq \abs{x}-\ell'\right)\lesssim (2m)^{-\beta \ell'}
\]
uniformly in $\ell'$.
\end{prop}

As usual, in this proposition the sense of ``for any subword under
$\mu_\ell$'' is that of Definition~\ref{subword}.

\begin{dem}
This is simply a rewriting of Axioms~2 and~3, combined to the observation
that the choice of the test-length and of the small words $u$ and $v$
(which are of length $O(\log\ell)$) only introduces a polynomial factor in
$\ell$.
\end{dem}

In our proof, we will also need the fact that the number of possible
apparent lengths for subwords under $\mu_\ell$ grows subexponentially
with $\ell$. So we need at least a rough upper bound on the apparent
length.

By definition, if $x$ appears with probability $p$ as a subword under
$\mu_\ell$, then by symmetry $y$ will by equal to $x\~$ with the same
probability, and thus the probability that $xuyv=e$ is at
least $p^2$ (taking $u=v=e$). Thus $\al_{\abs{x}}(x)\leq
-\frac{2}{\beta}\log p-\abs{x}$. Reversing the equation, this means that
for any subword $x$ under $\mu_\ell$, we have that $\Pr(\al(x)\geq L)\leq
(2m)^{-\beta(L - \abs{x})/2}$.

In particular, taking $L$ large enough ($L\geq 4\ell$ is enough) ensures
that in a set of $(2m)^{d\ell}$ randomly chosen elements under $\mu_\ell$
with $d<\beta$, subwords of length greater than $L$ only occur with
probability exponentially small as $\ell\rightarrow\infty$. Thus, we can
safely assume that all subwords of words of $R$ have apparent length at
most $4\ell$.

\bigskip

In the applications given in this text to plain random words or random
geodesic words, apparent length has more properties, especially a very
nice behavior under multiplication by a random word. In the geodesic word
model, apparent length is simply the usual length. We do not explicitly
need these properties, though they are present in the inspiration of our
arguments, and thus we do not state them.

\subsection{The main argument}
\label{themainargument}

Now we enter the main step of the proof. We will consider a
davKd and evaluate the probability that it is fulfillable. We will see
that either the davKd satisfies some isoperimetric inequality, or this
probability is very small (exponential in $\ell$). It will then be enough
to sum on all davKd's with at most $K$ faces to prove
Proposition~\ref{mainprop}.

In our graph $\Gamma$, the ordinary vertices represent letters of random
relators. Say $\Gamma$ has $n\ell$ ordinary vertices, that is, the faces
of $D''$ bear $n$ different relators of $R$.

We will use the term \emph{letter} to denote one of these vertices.
Enumerate letters in the obvious way from $1$ to $n\ell$, beginning with
the first letter of the first relator. So, a letter is a number between
$1$ and $n\ell$ indicating a position in some relator. Relators are
random words on elements of the generating set $S$ of $G$, so if $i$ is a
letter let $f_i$ be the corresponding element of $S$.

Since the relators are chosen at random, the $f_i$'s are random
variables.

As in the case of random quotients of the free group, the idea is to
construct the graph $\Gamma$ step by step, and evaluate the probability
that at each step, the conditions imposed by the graph are satisfied by
the random set $R$ of relators. We will construct the graph by groups of
successive letters joined to the same translators, and use the notion of
apparent length (see Definition~\ref{defal}) to keep track of the
probabilities involved at each step.

For a letter $i$, write $i\in t$ if $i$ is joined to translator $t$. For
$1\leq a\leq n$, write $i\in a$ to mean that letter $i$ belongs to the
$a$-th part of the graph. So $r_a$ is the product of the $f_i$'s for $i\in a$.

\bigskip

Consider an internal translator $t$. There is a word $w$ associated to
it, which writes $w=u\delta_1v\delta_2$ where $\delta_{1,2}$ are short
and $u$ and $v$ are subwords of the random relators. The subwords $u$ and
$v$ are products of letters, say $u=f_p\ldots f_q$ and $v=f_r\ldots f_s$.
Reserve these notations $w(t)$, $u(t)$, $v(t)$, $p(t)$, $q(t)$, $r(t)$
and $s(t)$. Give similar definitions for boundary translators and
commutation translators.

Call $u$ and $v$ the \emph{sides} of translator $t$. The translator
precisely imposes that there exist short words $\delta_1, \delta_2$ such
that $u\delta_1v\delta_2=e$ in $G$. We will work on the probabilities of
these events.

Some of the translators may have very small sides; yet we are to apply
asymptotic relations (such as the definition of cogrowth) which ask for
arbitrarily long words. As there are at most $4K/\eps_2$ translators,
with at most two sides each, the total length of the sides which are of
length less than $\eps_3\ell$ does not exceed $\eps_3\ell.8K/\eps_2$.
Setting $\eps_3=\eps_2^2/8K$ ensures that the total length of these sides
is less than $\eps_2\ell$.

Call \emph{zero-sided translator} an internal translator both sides of
which have length less than $\eps_3\ell$. Call \emph{two-sided
translator} an internal translator having at least one side of length at
least $\eps_3\ell$ and its smaller side of length at least $\eps_3$ times
the length of its bigger side. Call \emph{one-sided translators} the rest
of internal translators.

Throw away all zero-sided translators from the graph $\Gamma$.  This
throws away a total length of at most $\eps_2\ell$, and do not call
sides any more the small sides of one-sided translators. Now we have two-sided
translators, one-sided translators, commutation
translators and boundary translators, all sides of which have length at
least $\eps_3^2\ell$. So if $\ell$ is large enough (depending on
$\eps_3$) we can apply the probability evaluations of
Axioms~1-4 without trouble.

For a letter $i$, say that translator $t$ is finished at time $i$ if
$i\geq s(t)$. Say that two-sided translator $t$ is half-finished
at time $i$ if $q(t)\leq i<r(t)$.

\bigskip

Add a further decoration to $\Gamma$ (and to the davKd): for each two-sided 
translator $t$, specify an integer $L(t)$ between $0$ and $4\ell$
(remember we can suppose that every subword has apparent length at most
$4\ell$). This will represent the apparent length of the half-word $u(t)$
associated to the diagram when it is half-finished. In the same vein,
specify an integer $L(b)$ between $0$ and $4\ell$ for each boundary
translator $b$, which will represent the apparent length of the word
$u(b)$ when $b$ is finished. We want to show that the boundary length is
big, so we want to show that these apparent lengths of boundary
translators are big. What we will show is the following: if the sum of
the imposed $L(b)$'s for all boundary translators $b$ is too small, the
probability that the diagram is fulfillable is small.

Now say that a random set of relators $r_1,\ldots,r_n$ \emph{fulfills}
the conditions of $\Gamma$ up to letter $i$ if for any internal or
commutation translator $t$ which is finished at time $i$, the corresponding
word $w(t)$ is trivial in $G$; and if, for any half-finished two-sided
translator $t$, the apparent length of the half-word $u(t)$ is $L(t)$; and if, for each finished boundary translator $b$, the apparent
length of $u(b)$ is $L(b)$.

(An apparent length is not necessarily an integer; by
prescribing the apparent length of $u(t)$, we prescribe only the integer
part. As $\ell$ is big the discrepancy is totally negligible and we will
not even write it in what follows.)

Of course, fulfillability of the davKd implies fulfillability of $\Gamma$
up to the last letter for some choice of $r_1,\ldots,r_n \in R$ and for
some choice of the $L(t)$'s. (It is not exactly equivalent as we threw
away some small proportion $\eps_1$ of the edges.)

For a given relator $r$, there may be some translators having a side made
of an initial and final piece of $r$, so that the side straddles the
first letter of $r$. As we will fill in letters one by one starting with
the first ones, we should treat these kind of translators in a different
way. The simplest way to treat this little problem is a further cutting
of the translators that straddle the beginning of a word, using
Proposition~\ref{coarsening}, as is best explained by a figure (the thick
dot represents the beginning of some relator).

\begin{center}
\includegraphics{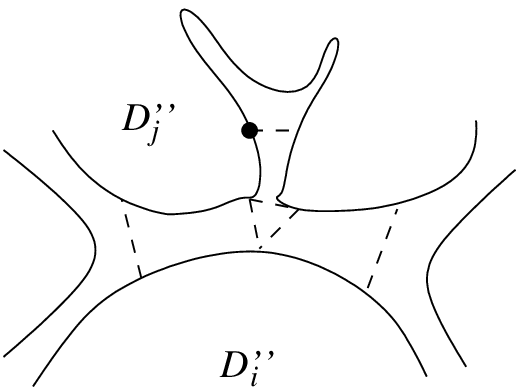}
\end{center}

\bigskip

Up to now there are three free variables in our argument: $K$, the
maximal number of new cells in diagrams we consider; and $\eps_1$ and
$\eps_2$, which are linked to the way we cut translators to eliminate
doublets.

\begin{prop}
\label{isomainprop}
For every density $d<\beta$, for every $K$, there exists
$\eps_1,\eps_2> 0$ such that, if $\ell$ is large enough, then, for any
davKd $\mathcal{D}$, either $\mathcal{D}$ satisfies a
$\frac{\kappa_2}{4}(1-d/\beta)$-isoperimetric
inequality (in the sense of Definition~\ref{isodavKd}),
where $\kappa_2$ is the constant in Axiom~2,
or the probability that $\mathcal{D}$ is fulfillable
is less than $(2m)^{-\ell(\beta-d)/4}$.
\end{prop}

Before proceeding to the proof of this proposition, let us show how it
implies Proposition~\ref{mainprop}, via Proposition~\ref{davKdtovKd}.

If we know that the number of distinct davKd's associated to a van Kampen
diagram satisfying the assumptions of Proposition~\ref{mainprop} is
polynomial in $\ell$, then summing the probability evaluation of
Proposition~\ref{isomainprop} on all such davKd's we can conclude: in
this case, the probability that there exists a davKd violating the
isoperimetric inequality is exponentially small, and so any van Kampen
diagram will satisfy an isoperimetric inequality, since any van Kampen
diagram satisfying the assumptions of Proposition~\ref{mainprop}
and~\ref{lemsmallold} has an associated davKd.
So we will evaluate the number of davKd's with at most $K$ faces.

But by Proposition~\ref{numberofdavKds}, the number of possible davKd's is
polynomial in $\ell$ at fixed $K$. We have to beware we added some extra
decoration to the davKd in between: in the elimination of doublets (we
made at most $K/\eps_2$ more cuttings, which can be kept track of by as
many numbers between $1$ and $K\ell$), and when prescribing an apparent
length for each internal translator (at most $4K/\eps_2$ numbers between
$1$ and $4\ell$). So the number of possibilities remains polynomial in
$\ell$ (all other things being fixed).

This proves that the probability that there exists a davKd violating the
isoperimetric inequality decreases exponentially with $\ell$, hence
Proposition~\ref{mainprop}.

\bigskip

\renewcommand{\Pr}{\mathrm{P}}

\begin{dem}[ of Proposition~\ref{isomainprop}]
Choose some integer $K$.
It is time to fix the parameters $\eps_1$, $\eps_2$. Recall we
set $\eps_3=\eps_2^2/8K$.
Also recall that the sides of translators are of length at least
$\eps_3^2\ell$, so that we will take $\ell$ large enough depending on
$\eps_3$ (that is, depending on $K$ and on the axioms).

With foresight, let
$\eps=\eps_1+3\eps_2+\gamma_4\eps_2/\beta+\eps_3/\kappa_2$ where
$\kappa_2, \gamma_4, \beta$ are the constants appearing in the axioms.
Choose $\eps_1$ and $\eps_2$ small enough so that $\eps\leq
(1-d/\beta)/4$.  These choices depend on $K$, $d$ and $G$ but not on $\ell$
neither on any diagram.

\bigskip

Let $\Pr_i$ be the probability that \emph{some fixed choice} of $n$
relators $r_1,\ldots,r_n \in R$ under our law $\mu_\ell$ fulfills $\Gamma$
up to letter $i$. This does not take into account the choice of $n$
relators among the $(2m)^{d\ell}$ relators of the presentation.
The quantity $\Pr_i$ depends only on the davKd $\mathcal{D}$ and on the law
$\mu_\ell$ of the relators.

Let $1\leq a\leq n$ (recall $n$ is the number of parts of the graph, or
the number of different relators of $R$ appearing in the diagram). Let
$m_a$ be the number of times relator $a$ appears in the diagram. Let
$i_0$ be the first letter of $a$, and $i_f$ the last one.

Let $\Pr^a$ be the probability that \emph{there exists} a choice of relators
$r_1,\ldots,r_a$ in $R$ fulfilling the conditions of $\Gamma$ up to
letter $i_f$ (the last letter of $a$). As there are by definition
$(2m)^{d\ell}$ choices for each relator, we have
\[
\Pr^a/\Pr^{a-1}\leq (2m)^{d\ell} \Pr_{i_f}/\Pr_{i_0-1}
\]
which expresses the fact that when we have fulfilled up to part $a-1$,
to fulfill up to part $a$ is to choose the $a$-th relator in $R$ and to
see if the letters $f_{i_0},\ldots, f_{i_f}$ of this relator fulfill the
conditions imposed on the $a$-th part of the graph by the translators.

Let $A_a$ be the sum of all $L(t)$'s for each two-sided translator $t$ which
is half-finished at time $i_f$, plus the sum of all $L(b)$'s for each boundary
translator $b$ which is finished at time $i_f$. We will study
$A_a-A_{a-1}$.

\begin{lem}
\label{mainlem}
For any davKd $\mathcal{D}$ with at most $K$ faces,
for any $1\leq a \leq n$ we have
\[
A_a-A_{a-1}\geq
m_a\left(\ell(1\!-\eps)+\frac{\log_{2m}\Pr_{i_f}-\log_{2m}\Pr_{i_0-1}}{\beta}\right)
+o(\ell) \quad (\star)
\]
where the constant implied in $o(\ell)$ depends on $K$ but not on the
diagram $\mathcal{D}$.
\end{lem}

Before proving this lemma, let us finish the proof of
Proposition~\ref{isomainprop}.

Recall we saw above that
\[
\Pr^a/\Pr^{a-1}\leq (2m)^{d\ell}\Pr_{i_f}/\Pr_{i_0-1}
\]
where the $(2m)^{d\ell}$ factor accounts for the choice of the relator
$r_a$ in $R$.

Set $d_a=\log_{2m} \Pr^a$ (compare the case of random quotients of
$F_m$). Beware the $d_a$'s are non-positive.
From $(\star)$ we get
\[
A_a-A_{a-1} \geq
m_a\left(\ell(1-\eps)+\frac{d_a-d_{a-1}-d\ell}{\beta}\right)+o(\ell)
\]

Compare this to the equation linking dimension and number of edges on
page~\pageref{edgedim} (and recall that here $A_a$ is not the number of
edges but the apparent length, which varies the opposite way, and that we
want it to be big).

\bigskip

Summing the inequalities above for $a$ from $1$ to $n$ gives
\begin{eqnarray*}
A_n 
&\geq&
\ell (1-\eps) \sum m_a - \frac{d\ell}{\beta}\sum m_a
+\frac1{\beta}\sum m_a(d_a-d_{a-1})+o(\ell)
\\&=&
\abs{\mathcal{D}} \ell\left(1-\eps-\frac{d}{\beta}\right) +\frac1{\beta}\sum
d_a(m_a-m_{a+1})+o(\ell)
\end{eqnarray*}

The number of summands is $n\leq K$, so that the constant in $o(\ell)$ is
controlled by $K$ again.

At the end of the process, all translators are finished, so by definition
$A_n$ is simply the sum of the apparent lengths of all boundary
translators, that is $A_n=\sum_b L(b)$.

Now recall that (if $\mathcal{D}$ is ever fulfillable) a boundary translator $b$ means the existence of an
equality $e=\delta_1u\delta_2 v$ in $G$, with by assumption
$\al(u)=L(b)$, the $\delta$'s of length at most $2E\log\ell$, and $v$
lying on the boundary of the diagram. By the definition of apparent
length (Definition~\ref{defal2} which takes Axiom~2 into account), we have
$\norm{u}\geq\kappa_2 \al(u)=\kappa_2 L(b)$, thus $\norm{v}\geq
\norm{u}-\norm{\delta_1}-\norm{\delta_2}\geq \kappa_2 L(b)+o(\ell)$.
As $v$ lies on the boundary of $\mathcal{D}$ this implies
\[
\abs{\d \mathcal{D}} \geq \kappa_2 A_n + o(\ell)
\]
(once again we can sum the $o(\ell)$'s harmlessly since the number of
translators is bounded by some function of $K$.)

So using the
lower bound for $A_n$ above we get
\[
\abs{\d \mathcal{D}}\geq
\abs{\mathcal{D}}\ell\left(1-\eps-d/\beta\right)\kappa_2
+\frac{\kappa_2}{\beta} \sum d_a(m_a-m_{a+1})\,+o(\ell)
\]

Recall we managed to choose $\eps\leq(1-d/\beta)/4$. Also take $\ell$
large enough so that the $o(\ell)$ term is less than $\ell(1-d/\beta)/4$
(such an $\ell$ depends on $K$).
The inequality above rewrites
\[
\abs{\d \mathcal{D}}\geq \abs{\mathcal{D}}\ell\left(1-d/\beta\right)\kappa_2/2
+\frac{\kappa_2}{\beta} \sum d_a(m_a-m_{a+1})
\]

We are free to choose the order of the enumeration of the parts of the
graph. In particular, we can suppose that the $m_a$'s are non-increasing.

As $\sum m_a=\abs{\mathcal{D}}$, we have $\sum d_a(m_a-m_{a+1})\geq
\abs{\mathcal{D}}\inf
d_a$ (recall the $d_a$'s are non-positive). Thus
\[
\abs{\d \mathcal{D}}\geq \frac{\kappa_2}{2\beta} \abs{\mathcal{D}}\ell\left(\beta-d+2\inf
d_a/\ell\right)
\]

By definition, the probability that the davKd is fulfillable
is less than $(2m)^{d_a}$ for all $a$. This probability is then less than
$(2m)^{\inf d_a}$.

First suppose that $\inf d_a\geq -\ell(\beta-d)/4$.
Then we have the isoperimetric inequality 
\[
\abs{\d \mathcal{D}}\geq
\frac{\kappa_2}{4}\,\ell\abs{\mathcal{D}}(1-d/\beta)
\]
as needed.


Or, second, suppose $\inf d_a<-\ell(\beta-d)/4$.
This means that the probability that the davKd is fulfillable is less
than $(2m)^{-\ell(\beta-d)/4}$.

This proves Proposition~\ref{isomainprop} assuming Lemma~\ref{mainlem}.
\end{dem}

\bigskip

\begin{dem}[ of Lemma~\ref{mainlem}]
The principle of the argument is to look at the evolution of the apparent
length of the translators, where the apparent length of a translator at
some step is the apparent length of the part of this translator which is
filled in at that step. We will show that our axioms imply that when we
add a subword of some
length, the probability that the increase in apparent length is less than
the length of the subword added is exponentially small, such that a
simple equation is satisfied:
\[
\Delta \al \geq \abs{.}+\frac{\Delta \log \Pr}{\beta}
\]
(where $\Delta$ denotes the difference between before and after filling
the subword). This will be the motto of our forthcoming arguments.

But at the end of the process, the word read on any internal translator
is $e$, which is of apparent length $0$, so that the only contribution to
the total apparent length is that of the boundary translators, which we
therefore get an evaluation of.

\bigskip

Now for a rigorous exposition.  The difference between $A_a$ and
$A_{a-1}$ is due to internal translators which are half-finished at time
$i_0$ and are finished at time $i_f$, to internal translators which are
not begun at time $i_0$ and are half-finished at time $i_f$, and to
boundary translators not begun at time $i_0$ but finished at time $i_f$:
that is, all internal or boundary translators joined to a letter between
$i_0$ and $i_f$.

First, consider a two-sided translator $t$ which is not begun at time
$i_0$ and half-finished at time $i_f$. Let $\Delta_t A_a$ be the contribution
of this translator to $A_a-A_{a-1}$, we have $\Delta_t A_a=L(t)$ by
definition. Taking
notations as above, we have an equality $e=u\delta_1v\delta_2$ in $G$. By
assumption, to fulfill the conditions imposed by $\Gamma$ we must have
$\al(u)=L(t)$. The word $u$ is a subword of the part $a$ of $\Gamma$ at
play. But Proposition~\ref{synthax} (that is, Axioms~2 and~3) tells us that, conditionally to
whatever happened up to the choice of $u$, the probability that
$\al(u)=L(t)$ is roughly less than $(2m)^{-\beta(\abs{u}-L(t))}$. Thus,
taking notations as above, with
$p$ the first letter of $u$ and $q$ the last one, we have
\[
\Pr_q/\Pr_{p-1} \lesssim (2m)^{-\beta(\abs{u}-L(t))}
\]
or, taking the log and decomposing $u$ into letters:
\[
\Delta_t A_a\geq \sum_{i\in t, i\in a}
1+\frac{\log_{2m}\Pr_i-\log_{2m}\Pr_{i-1}}{\beta} \ +o(\ell)
\]
where $1$ stands for the length of one letter (!).
Note that a rough evaluation of the probabilities gives an evaluation up
to $o(\ell)$ of the apparent lengths.

This is the rigorous form of our motto above.

Second, consider an internal translator $t$ which is half-finished at time
$i_0$ and finished at time $i_f$. Let $\Delta_t A_a$ be the contribution of
this translator to $A_a-A_{a-1}$, we have $\Delta_t A=-L(t)$.  Taking
notations as above, we have an equality $e=u\delta_1v\delta_2$ in $G$. By
assumption, we have $\al(u)=L(t)$. But the very definition of apparent
length (Definition~\ref{defal}) tells us that given $u$, whatever happened before
the choice of $v$, the probability that there exist such $\delta_{1,2}$
such that $e=u\delta_1v\delta_2$ is at most $(2m)^{-\beta(\al(u)+\abs{v})}$.
Thus
\[
\Pr_s/\Pr_{r-1} \lesssim (2m)^{-\beta(L(t)+\abs{v})}
\]
where $r$ and $s$ are the first and last letter making up $v$. Or, taking the log and decomposing $v$ into letters:
\[
\Delta_t A_a\geq \sum_{i\in t, i\in a}
1+\frac{\log_{2m}\Pr_i-\log_{2m}\Pr_{i-1}}{\beta} \ +o(\ell)
\]
which is exactly the same as above.

Third, consider an internal translator $t$ which is not begun at time
$i_0$ and finished at time $i_f$, that is, $t$ is joined to two subwords of the
part $a$ of the graph at play. As we removed doublets, the subwords $u$
and $v$ are disjoint, and thus we can work in two times and apply the
two cases above, with first $t$ going from not begun state to
half-finished state, then to finished state. The contribution of $t$ to
$A_a-A_{a-1}$ is $0$, and summing the two cases above we get
\[
\Delta_t A_a =0\geq \sum_{i\in t}
1+\frac{\log_{2m}\Pr_i-\log_{2m}\Pr_{i-1}}{\beta} \ +o(\ell)
\]
which is exactly the same as above.

Fourth, consider a commutation translator $t$ which is not begun at time
$i_0$ and is finished at time $i_f$. Write as above that
$e=\delta_1u\delta_2u\~$ in $G$, with $\delta_1$ and $\delta_2$ of length
at most $\eps_2\abs{u}$. By Axiom~4, whatever happened before the choice
of $u$, this event has probability
roughly less than $(2m)^{\gamma_4\eps_2\abs{u}-\beta\abs{u}}$
where $\gamma_4$ is some constant. Take $\eps_4=\gamma_4\eps_2/\beta$, and as
usual denote by $p$ and $q$ the first and last letters making
up $u$. We have shown that
\[
\Pr_q/\Pr_{p-1} \lesssim (2m)^{-\beta\abs{u}(1-\eps_4)}
\]
Take the log, multiply everything by two (since each letter joined to the
commutation diagram $t$ is joined to it by a double edge), so that
\[
\Delta_t A_a=0 \geq \sum_{i \in t} 2(1-\eps_4)
+2\frac{\log_{2m}\Pr_i-\log_{2m}\Pr_{i-1}}{\beta} \ +o(\ell)
\]

Fifth, consider a one-sided translator not begun at time $i_0$ and
finished at time $i_f$. We have an equality $e=u\delta_1v\delta_2$ in
$G$, where $\delta_{1,2}$ have length $O(\log\ell)$ and $\abs{v}\leq
\eps_3\abs{u}$ (by definition of a one-sided translator), so that
$\norm{u}\leq \eps_3\abs{u}+O(\log\ell)$. But by Axiom~2, this has
probability roughly less than $(2m)^{-\beta\abs{u}(1-\eps_3/\kappa_2)}$, so once
again, setting $\eps_5=\eps_3/\kappa_2$:
\[
\Delta_t A_a=0\geq \sum_{i\in t, i\in a}
(1-\eps_5)+\frac{\log_{2m}\Pr_i-\log_{2m}\Pr_{i-1}}{\beta} \ +o(\ell)
\]

Sixth (and last!), consider a boundary
commutator $t$ that is not begun at time $i_0$ and is finished at time
$i_f$. Its situation
is identical to that of an internal translator half-finished at time $i_f$
(first case above), and we get
\[
\Delta_t A_a=L(t)\geq \sum_{i\in t, i\in a}
1+\frac{\log_{2m}\Pr_i-\log_{2m}\Pr_{i-1}}{\beta} \ +o(\ell)
\]

\bigskip

We are now ready to conclude. Sum all the above inequalities for all
translators joined to part $a$: 
\begin{eqnarray*}
A_a-A_{a-1} &=&\sum_{t\text{ translator joined to }a} \hspace{-2em} \Delta_t A_a
\\&\geq &
\sum_{
\begin{subarray}{c}
t\text{ non-commutation translator}
\\i\in t, i\in a
\end{subarray}
}
\hspace{-3em} (1-\eps_5)+\frac{\log_{2m}\Pr_i-\log_{2m}\Pr_{i-1}}{\beta}
\\& & +
\sum_{
\begin{subarray}{c}
t\text{ commutation translator}
\\i\in t, i\in a
\end{subarray}
}
\hspace{-3em} 2(1-\eps_4)+2\frac{\log_{2m}\Pr_i-\log_{2m}\Pr_{i-1}}{\beta}
\\& & +o(\ell)
\end{eqnarray*}

Recall $m_a$ is the number of times the $a$-th relator appears in the
van Kampen diagram. The way we constructed the graph, any vertex representing a
letter of the $a$-th relator is joined to $m_a$ translators (except for a
proportion at most $\eps_1+3\eps_2$ that was removed), counting commutation
translators twice. Thus, in the sum above, each of the $\ell$ letters of $a$ appears
exactly $m_a$ times, and so
\[
A_a-A_{a-1}\geq
m_a\left(\ell(1\!-\!\eps_4\!-\!\eps_5\!-\!\eps_1\!-\!3\eps_2)+\frac{\log_{2m}\Pr_{i_f}-\log_{2m}\Pr_{i_0-1}}{\beta}\right)
+o(\ell)
\]

(Because of the removal of a proportion at most $\eps_1+3\eps_2$ of the
letters, some terms $\log_{2m}\Pr_{i_f}-\log_{2m}\Pr_{i_0-1}$ are
missing in the sum; but as for any $i$, we have $\Pr_i\leq \Pr_{i-1}$,
the difference of log-probabilities $\log_{2m}\Pr_{i}-\log_{2m}\Pr_{i-1}$
is non-positive, and the inequality is true \textit{a fortiori} when we add
these missing terms.)

Note that there is nothing bad hidden in the summation of the $o(\ell)$
terms, since the number of terms in the sum is controlled by the
combinatorics of the diagram (i.e.\ by $K$), and depends in no way on
$\ell$.

Recall we set
$\eps=\eps_1+3\eps_2+\gamma_4\eps_2/\beta+\eps_3/\kappa_2=\eps_1+3\eps_2+\eps_4+\eps_5$,
which is exactly what we get here. So Lemma~\ref{mainlem} is proven.
\end{dem}

All pending proofs are finished;
hence hyperbolicity of the random
quotient when $d<\beta$.

\subsection{Non-elementarity of the quotient}

We now prove that if $d<\beta$, the quotient is infinite and not
quasi-isometric to $\Z$.

\subsubsection{Infiniteness}
\label{infiniteness}

Let $d<\beta$. We will show that the probability that the random
quotient is finite decreases exponentially as $\ell\rightarrow \infty$.

We know from hyperbolicity of the quotient (Proposition~\ref{mainprop2})
that the probability that there exists
a van Kampen diagram of the quotient whose part made of old relators
is reduced and which is strongly reduced with respect to $G$, violating
some isoperimetric inequality, is exponentially close to $0$.

Imagine that $G/\langle R \rangle$ is finite. Then any element of the
quotient is a torsion element. Let $x$ be an element of the quotient,
this means that there exists a van Kampen diagram $D$ bordered by $x^n$ for some
$n$.

Now take for $x$ a random word picked under $\mu_\ell$. We will show that
such a random word is very probably not a torsion element in the
quotient. Instead of applying the previous section's results to the
random quotient of $G$ by $R$, consider the random quotient of $G$ by
$R\cup\{x\}$. Since $x$ is taken at random, $R\cup\{x\}$ is a random set
of words, whose density is only slightly bigger than $d$; this density is
$d'=\frac{1}{\ell}\log_{2m}\left((2m)^{d\ell}+1\right)$ which, if $\ell$
is large enough, is smaller than $\beta$ if $d$ is.

Now, if $G/\langle R \rangle$ is finite then $x$ is of torsion.  Set
$N=\abs{R}=(2m)^{d\ell}$. Consider the following family of diagrams. Let
$D$ be any abstract van Kampen diagram of $G/\langle R \rangle$ of boundary
length $n\ell$ for some $n$. Define the spherical diagram $E$ by gluing
$n$ faces of boundary size $\ell$ on the boundary of $D$ along their
border, and associate to each of the new faces the relator number $N+1$,
so that $D$ is an abstract van Kampen diagram with respect to $R\cup\{x\}$. If
$G/\langle R \rangle$ is finite then $x$ is of torsion, thus at least one
of the diagrams $E$ in this family is fulfillable with respect to $R\cup\{x\}$.

\begin{center}
\includegraphics{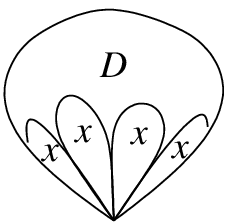}
\end{center}

By Proposition~\ref{strongreduction} we can take the strong reduction of
this diagram. It is non-empty as the faces bearing $x$ cannot be
cancelled (they all have the same orientation).

So there exists a strongly reduced van Kampen diagram of $G/\langle R\cup
\{x\}\rangle$ with boundary length $0$.

But we know by what we already proved (Propositions~\ref{mainprop}
and~\ref{mainprop2}) that, in the random quotient $G/\langle R\cup
\{x\}\rangle$ at density $d'$, the existence of such a diagram is of
probability exponentially close to $0$ as $\ell$ tends to infinity. This
ends the proof.

\subsubsection{Non-quasi$\Z$ness}

We show here that the random quotients for $d<\beta$ are not
quasi-isometric to $\Z$. Of course, we suppose $\beta>0$, which amounts,
in the case we consider (plain, or reduced, or geodesic words), to $G$
itself not being quasi-isometric to $\Z$.

We will reason in a similar manner as above to prove infiniteness. We
will consider a random quotient by a set $R$ of words at density $d$, and
we will add to $R$ two random words $x$ and $y$ picked under $\mu_\ell$,
thus obtaining a new random set of words at a density $d'>d$. As $\ell$
is big, $d'$ is only slightly above $d$, and if $\ell$ is big enough we
still have $d'<\beta$.

Say (from Proposition~\ref{mainprop2}) that any strongly reduced diagram
$D$ of the group $G/\langle R' \rangle$ satisfies an isoperimetric
inequality $\abs{\d D}\geq \alpha \ell \abs{D''}$ for some positive
$\alpha$, notations as above.

Suppose that $G/\langle R\rangle$ is quasi-isometric to $\Z$.

The two random elements $x$ and $y$ are either torsion elements or each of
them generates a subgroup of finite index. The case of torsion is handled
exactly as above in the proof of infiniteness.

Thus, suppose $x$ is not a torsion element. Let $h$ be the index
of the subgroup it generates. Of course $h$ depends on $x$.

For any $n\in \Z$, we can find a $p$ such that $y^n=x^pu$ in $G/\langle
R\rangle$, where $u$ is of length at most $h$. This equality defines a
van Kampen diagram of $G/\langle R\rangle$.

Now glue $n$ faces containing $y$ and $p$ faces containing $x$ to the
boundary of this diagram. This defines a van Kampen diagram of $G/\langle
R'\rangle$, which we can take the strong reduction of. This reduction is
non-empty since faces bearing $x$ and $y$ cannot be cancelled (so in
particular $\abs{D''}\geq n+p$). The boundary of this diagram is $u$.

But $n$ can be taken arbitrarily large, so we can take $n>\abs{u}/\alpha$. Then
the diagram has at least $n$ faces and boundary length $\abs{u}$, which
contradicts our isoperimetric inequality $\abs{\d D}\geq
\alpha\ell\abs{D''}$.

Of course, $u$, $n$ and $p$ depend on the random words $x$ and $y$. But
consider the following family of diagrams: for each $h \in \N$, each
$p\in \N$ and each $n\in \N$ such that $n> h/\alpha$, consider
all abstract van Kampen diagrams of length $h+n\ell+p\ell$, with the numbers on
the faces between $1$ and $N=\abs{R}$. Consider the diagrams obtained
from these by the following process: glue $p$ faces of size $\ell$
bearing number $N+1$ on the boundary, and $n$ faces of size $\ell$
bearing number $N+2$.

The reasoning above shows that if $G/\langle R\rangle$ is quasi-isometric
to $\Z$, then at least one of these abstract van Kampen diagrams is fulfillable
by a strongly reduced van Kampen diagram on the relators of $R'$. But all these
diagrams violate the isoperimetric inequality, hence the conclusion.

\paragraph{Alternate proof.} We give an alternate proof as it uses an
interesting property of the quotients. It works only in the case of a
random quotient by uniformly chosen plain words.

\begin{prop}
If $d>0$, then the abelianized of a random quotient of any group by
uniformly chosen plain random words is (with
probability arbitrarily close to $1$ as $\ell\rightarrow \infty$)
either $\{e\}$ or $\Z/2\Z$.
\end{prop}

(As usual, we find $\Z/2\Z$ when $\ell$ is even and there are no
relations of odd length in the presentation of $G$.)

Of course this is not necessarily true if $d=0$, since in this case the
number of relations added does not tend to infinity.

\begin{dem}

We want to show that a random quotient in density $d>0$ of the
free abelian group $\Z^m$ is trivial or $\Z/2\Z$.

Take a random word of length $\ell$ on $a_1^{\pm 1}, \ldots a_m^{\pm 1}$.
By the central limit theorem (or by an explicit computation on the
multinomial distribution), the number of times generator $a_i$ appears is
roughly $\ell/2m$ up to $\pm \sqrt{\ell}$. 

For the sake of simplicity, say that $\ell$ is a multiple of $2m$. The
probability that a random word $w$ is such that all relators $a_i$
and $a_j\~$ appear exactly $\ell/2m$ times in $w$ is equivalent to
\[
\frac{\sqrt{2m}}{(\pi\ell/m)^{(2m-1)/2}}
\]
by
the central limit theorem with $2m-1$ degrees of freedom or by a direct
computation using Stirling's formula.

This is equivalent as well to the probability that all $a_i$ and $a_j\~$
appear exactly $\ell/2m$ times, except for some $a_{i_0}$ appearing
$1+\ell/2m$ times and some $a_{j_0}$ appearing $\ell/2m-1$ times, this
deviation being negligible.

This probability decreases polynomially in $\ell$. But we choose an
exponential number of random words, namely $(2m)^{d\ell}$. So if $d>0$,
with very high probability we will choose a word $w$ in which all $a_i$ and
$a_j\~$
appear exactly $\ell/2m$ times, except for some $a_{i_0}$ appearing
$1+\ell/2m$ times and some $a_{j_0}$ appearing $\ell/2m-1$ times.

But $w=e$ in the quotient, and $w=e$ in an abelian group is equivalent to
$a_{i_0}a_{j_0}\~=e$ since all other relators appear exactly the same
number of times with exponent $1$ or $-1$ and thus vanish.

As this occurs arbitrarily many times, this will happen for all couples
of $i,j$. So these relators satisfy $a_i=a_j^{\pm 1}$ in the quotient for
all $i,j$. In particular, all relators are equal and moreover we have $a_i=a_i\~$.

Thus the abelianized is either $\{e\}$ or $\Z/2\Z$.
\end{dem}

\begin{cor}
A random quotient of a hyperbolic group by plain random words for $d<\beta$ is not
quasi-isometric to $\Z$.
\end{cor}

\begin{dem}
First take $d>0$. It is well-known (cf.\ \cite{SW}, Theorem~5.12, p.~178)
that a group which is quasi-isometric to $\Z$ has either $\Z$ or the
infinite diedral group $D_\infty$ as a quotient.

If $\Z$ is a quotient of the group, then its abelianized is at least $\Z$,
which contradicts the previous proposition. If $D_\infty$ is a quotient,
note that the abelianized of $D_\infty$ is $D_2=\Z/2\Z\times \Z/2\Z$,
which is still incompatible with the previous proposition. So we are done
if $d>0$.

Now if $d=0$, note that a random quotient with $d>0$ is a quotient of a
random quotient with $d=0$ (isolate the first relators). If the random
group at $d=0$ were quasi-isometric to $\Z$, then all of its quotients
would be either finite or quasi-isometric to $\Z$, which is not the case.
(Note that here we use hyperbolicity of $G$ to be allowed to apply our
main theorem, implying that random quotients are non-trivial for some
$d>0$. It may be that random quotients at $d=0$ of some groups are
quasi-isometric to $\Z$.)
\end{dem}

\bigskip

This ends the proof of Theorem~\ref{techmain}.

\appendix

\newpage

\section{Appendix: The local-global principle, or
Cartan-Hadamard-Gromov theorem}
\label{CHGP}

\bigskip

The Cartan-Hadamard-Gromov theorem allows to go from a local
isoperimetric inequality (concerning small figures in a given space) to
isoperimetry at large scale. It lies at the heart of our argument: to
ensure hyperbolicity of a group, it is enough to check the isoperimetric
inequality for a finite number of diagrams. This finite number depends,
of course, of the quality of the isoperimetric inequality we get on these
small diagrams. In particular, there is an algorithm to detect
hyperbolicity of a given group. We will use the form given by Papasoglu
(see~\cite{Pap}), who has written a completely combinatorial proof. See
also the presentation by Bowditch in~\cite{Bow}.

Let us state the form of the theorem we will use.

Let $X$ be a simplicial complex of dimension $2$ (all faces are
triangles). A \emph{circle drawn in $X$} is a sequence of consecutive
edges such that the endpoint of the last edge is the starting point of
the first one. A \emph{disk drawn in $X$} is a simplicial map from a
triangulated disk to $X$.

\def\tr{_{\mathrm{tr}}}

The \emph{area} $A\tr$ of a disk drawn in $X$ is its number of triangles. The
\emph{length} $L\tr$ of a circle drawn in $X$ is its number of edges. (Both with
multiplicity.) This is, $X$ is considered being made of equilateral
triangles of side $1$ and area $1$.

The \emph{area} of a drawn circle will be the smallest area of a drawn disk
with this circle as boundary, or $\infty$ if no such disk exists. The
\emph{length} of a drawn disk will be the length of its boundary.

\begin{thm}[ (P.~Papasoglu, cf.~\cite{Pap}, after M.~Gromov)]
Let $X$ be a simplicial complex of dimension $2$, simply connected.
Suppose that for some integer $K>0$, any circle $S$ drawn in $X$ whose
area lies between $K^2/2$ and $240K^2$ satisfies
\[
L\tr(S)^2\geq 2\cdot 10^4\,A\tr(S)
\]

Then any circle $S$ drawn in $X$ with $A(S)\geq K^2$ satisfies
\[
L\tr(S)\geq A\tr(S)/K
\]
\end{thm}

This theorem is a particular case of a more general theorem stated by
Gromov in \cite{Gro1}, section 6.8.F, for a length space.  Think of a
manifold. At very small scales, every curve in it satisfies the same
quadratic isoperimetric inequality as in the Euclidean space, with
constant $4\pi$. The theorem means that if, at a slightly larger scale,
the constant in this quadratic isoperimetric inequality becomes better
($2\cdot 10^4$ instead of $4\pi$), then isoperimetry propagates to large
scales, and at these large scales the isoperimetric inequality even becomes
linear. This is analogous to the fact that a control on the
curvature of a manifold (which is a local invariant) allows to deduce
global hyperbolicity properties. This was termed by Gromov
\emph{hyperbolic Cartan-Hadamard theorem} or \emph{local-global principle
for hyperbolic spaces}.

The proof of Papasoglu is based on considering the smallest diagram
violating the inequality to prove, and, by some surgery involving only
cutting it in various ways, to exhibit a smaller diagram violating the
assumptions. As this process only requires to consider subdiagrams of a
given diagram, he proves a somewhat stronger theorem.

\begin{thm}[ (P.~Papasoglu, cf.~\cite{Pap}, after M.~Gromov)]
Let $X$ be a simplicial complex of dimension $2$, simply connected. Let
$P$ be a property of disks in $X$ such that any subdisk of a disk having
$P$ also has $P$.

Suppose that for some integer $K>0$, any disk $D$ drawn in $X$ having
$P$, whose
area lies between $K^2/2$ and $240K^2$ satisfies
\[
L\tr(D)^2\geq 2\cdot 10^4\,A\tr(D)
\]

Then any disk $D$ drawn in $X$, having $P$, with $A(D)\geq K^2$, satisfies
\[
L\tr(D)\geq A\tr(D)/K
\]
\end{thm}

In the previous version, property $P$ was ``having the minimal area for a
given boundary'', hence the change from circles to disks.

\bigskip

We need to extend these theorems to complexes in which not all the faces
are triangular.

Let $X$ be a complex of dimension $2$. Let $f$ be a face of $X$.

The \emph{combinatorial length} $L_c$ of $f$ is defined as the number of
edges of its boundary. The \emph{combinatorial area} $A_c$ of $f$ is
defined as $L_c(f)^2$.

Let $D$ be a disk drawn in $X$. The \emph{combinatorial length} $L_c$ of
$D$ is the length of its boundary. The \emph{combinatorial area} $A_c$ of
$D$ is the sum of the combinatorial areas of its faces.

\begin{prop}
\label{localglobal}
Let $X$ be a complex of dimension $2$, simply connected. Suppose that a
face of $X$ has at most $\ell$ edges.
Let $P$ be a property of disks in $X$ such that any subdisk of a disk having $P$ also has $P$.

Suppose that for some integer $K\geq10^{10}\ell$, any disk $D$ drawn in $X$ having
$P$, whose
area lies between $K^2/4$ and $480K^2$ satisfies
\[
L_c(D)^2\geq 2\cdot 10^{14}\,A_c(D)
\]

Then any disk $D$ drawn in $X$, having $P$, with $A(D)\geq K^2$,
satisfies
\[
L_c(D)\geq A_c(D)/10^4K
\]
\end{prop}


\begin{dem}[ of the proposition]
Of course, we will show this proposition by triangulating $X$ and
applying Papasoglu's theorem.

The naive triangulation (cut a $n$-gon into $n-2$ triangles) does not
work since all triangles do not have the same size.

Triangulate all faces of $X$ in the following way: consider a face of $X$
with $n$ sides as a regular $n$-gon of perimeter $n$ in the Euclidean
plane. Consider a triangulation of the plane by equilateral triangles of
side $1$. (The polygon is drawn here with large $n$, so that it looks
like a circle.)

\begin{center}
\includegraphics{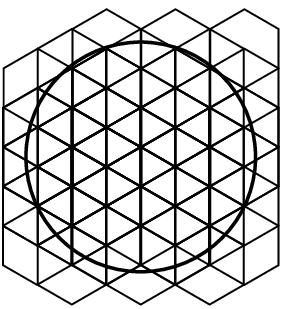}
\end{center}

This is not exactly a triangulation, but with a little work near the
boundary, we can ensure that the polygon is triangulated in such a way
that all triangles have sides between, say, $1/10$ and $10$ and area between
$1/10$ and $10$, so that the distortion between the triangle metric
and the Euclidian metric is a factor at most $10$. Note that the number
of triangles lies between $n^2/100$ and $100n^2$, as the
(Euclidian) area of our $n$-gon is roughly $n^2/4\pi$.

Let $Y$ be the simplicial complex resulting from $X$ by this
triangulation.

Let $L\tr$ and $A\tr$ be the length and area in $Y$ assigning length
$1$ to each edge and area $1$ to each triangle. Let $L_c$ and $A_c$ be
the length and area in $X$ defined above; in $Y$ they can be used for
disks coming from $X$.

Let $L$ and $A$ be the Euclidean length and area in $Y$, that is, each
face of $X$ with $n$ edges is a regular $n$-gon, and the triangles are
given their length and area coming from the triangulation above in the
Euclidean plane.

The discrepancy between $L\tr$, $L$ and $L_c$, and between $A\tr$, $A$
and $A_c$, is at most a factor $100$.

We proceed as follows: We will show that a disk in $Y$ with property $P$,
whose area $A\tr(B)$ lies
between $K^2/2$ and $240K^2$, satisfies $L\tr(B)^2\geq2\cdot 10^4
A\tr(B)$. Then, by the above theorem, any disk $B$ of area $A\tr(B)\geq
K^2$ will satisfy $L\tr(B)\geq A\tr(B)/K$, thus $L_c(B)\geq
A_c(B)/10^4K$ and we will be done.

Let $B$ be a disk in $Y$ with property $P$, whose area $A\tr(B)$ lies
between $K^2/2$ and $240K^2$. We want to show that it satisfies
$L\tr(B)^2\geq2\cdot 10^4 A\tr(B)$.

There are two kinds of disks drawn in $Y$: those who come from a disk
drawn in $X$, and those which there exists faces of $X$ that are
only partially contained in.

For the first kind we are done: by assumption, we have $L_c(B)^2\geq
2\cdot 10^{14} A_c(B)$, which implies $L\tr(B)^2\geq 2\cdot 10^4
A\tr(B)$.

So we want to reduce the problem to this kind of disks.

We will need the following isoperimetric lemmas:

\begin{lem}
\label{lemiso1}
Let $C$ be a regular closed curve in a Euclidean disk $D$. Suppose that $C$ encloses a surface of area at most half the area of
$D$. Then the length of the intersection of $C$ with the boundary of $D$
is at most $32$ times the length of the intersection of $C$ with the
interior of $D$.
\end{lem}

(One would expect an optimal constant $\pi/2$ with optimal $C$ enclosing
a half disk.)

This lemma is shown in \cite{Gro3}, 6.23. The next lemma is a formal
consequence thereof.

\begin{lem}
\label{lemiso2}
Let $C$ be a regular closed curve in a Euclidean disk $D$. Suppose that
$C$ encloses a surface of area at least half the area of $D$. Then the
length of the intersection of $C$ with the interior of $D$ is at least
$1/32$ times the length of $\partial D\setminus C$.
\end{lem}

The next lemma is a consequence of the first one and of the usual
isoperimetric inequality in the Euclidean plane.

\begin{lem}
\label{lemiso3}
Let $C$ be a regular closed curve in a Euclidean disk $D$. Suppose that
$C$ encloses a surface of area at most half the area of $D$. Then the
square of the length of the intersection of $C$ and the interior of $D$
is at least $1/100$ times the area enclosed by $C$.
\end{lem}

\bigskip

Now back to our disk $B$ in $Y$.

Let $D$ be a face of $X$ such that $B$ intersects $D$.

Suppose that $\d B \cap D$ is connected (that is, $B$ intersects $D$ only
once; otherwise, make the following construction for each of the
connected components). Compare the Euclidean area of $B\cap D$ with that
of $D$. If it is more than one half, enlarge $B$ such that it includes
all of $D$.

Follow this process for each face $D$ of $X$ partially intersecting $B$.

Let $B'$ be the disk in $Y$ obtained after this process. By construction,
we have $A(B)\leq A(B')\leq 2A(B)$. By Lemma~\ref{lemiso2}, we have
$L(B')\leq 32L(B)$.

Now, for each face $D$ of $X$ intersecting $B'$, either $D\subset B'$ or
the area of $D\cap B'$ is at most one half the area of $D$.

As a first case, suppose that the cumulated area of all such $D$ which
are included in $B'$ is at least one half of the area of $B'$. Define
$B''$ by amputing $B'$ from all faces $D$ of $X$ which are not totally
included in $B'$. By assumption, we have $A(B')\geq A(B'')\geq A(B')/2$.
And it follows from Lemma~\ref{lemiso1} that $L(B'')\leq 32 L(B')$.

By construction, the disk $B''$ is now a disk made of whole faces of $X$.
As $A(B)/2\leq A(B'')\leq 2A(B)$, we have $K^2/4\leq A(B'')\leq 480K^2$.
We can thus apply the isoperimetric assumption: $L(B'')^2\geq 2\cdot
10^{14} A(B'')$. Since $L(B'')\leq 32^2 L(B)$ and $A(B)\leq 2A(B'')$, we
get that $L(B)^2\geq 2\cdot 10^{10} A(B)$, hence $L\tr(B)\geq 2\cdot 10^4
A\tr(B)$.

As a second case, imagine that the cumulated area of all such $D$ which
are wholly included in $B'$ is less than half the area of $B'$. Let $D_i$
be the faces of $X$ intersecting $B'$ but not wholly contained in $B'$.
Let $a_i=A(D_i\cap B')$. We have $\sum a_i\geq A(B')/2\geq K^2/4$.

Let $m_i=L(\d B'\cap D_i)$. By Lemma~\ref{lemiso3}, we have $m_i^2\geq
a_i/100$.

Since any face of $X$ has at most $\ell$ edges, we have $A_c(D_i\cap
B')\leq \ell^2$, so for any $i$, $a_i\leq 100\ell^2$. Group the indices
$i$ in packs $I$ so that for each $I$, we have $100\ell^2 \leq \sum_{i\in
I}a_i \leq 200\ell^2$. There are at least $K^2/800\ell^2$ packs $I$. Let
$M_I=\sum_{i\in I} m_i$.

We have
\[
M_I=\sum_{i\in I} m_i \geq \sqrt{\sum_{i\in I} m_i^2} \geq
\sqrt{\sum_{i\in I}a_i/100} \geq \ell
\]
and
\[
L(B')^2\geq \left(\sum_i m_i\right)^2=\left(\sum_I M_I\right)^2
\geq \left(\sum_I \ell\right)^2
\]
and as there are at least $K^2/800\ell^2$ packs
\[
L(B')^2\geq K^4/10^6\ell^2\geq A(B') K^2 /10^9\ell^2
\]
as $A(B')\leq 480 K^2$. Now as $L(B')\leq 32L(B)$ and $A(B')\geq A(B)$ we
have
\[
L(B)^2\geq A(B)K^2/10^9\ell^2
\]
or
\[
L\tr(B)^2\geq A\tr(B) K^2/10^{15}\ell^2
\]
and we are done as $K^2\geq 10^{20}\ell^2$.

This ends the proof of the proposition.
\end{dem}

\newpage

\section{Appendix: Conjugacy and isoperimetry in hyperbolic groups}
\label{appiso}

\bigskip

We prove here some of the statements needed in the text about conjugacy
of words and narrowness of diagrams in hyperbolic groups. For general
references on hyperbolic groups and spaces we refer to~\cite{BH},
\cite{CDP} or~\cite{GH}.

Throughout this appendix, $G$ will denote a hyperbolic discrete group
generated by a finite symmetric set $S$, defined by a finite set of
relations $R$ (every discrete hyperbolic group is finitely presented, cf.\ \cite{S}). Let $\delta$ be a hyperbolicity
constant w.r.t.\ $S$.

A \emph{word} will be a word made of letters in $S$. The \emph{length} of a word
$w$ will be its number of letters (regardless of whether it is equal to a
shorter word in the group), denoted by $\abs{w}$.

Equality of words will
always be with respect to the group $G$.

Let $C$ be an isoperimetric constant for $G$, i.e.\ a positive number
such that any simply connected minimal van Kampen diagram $D$ on $G$
satisfies $\abs{\d D}\geq C \abs{D}$. See section~\ref{defs} for
definitions and references about diagrams and isoperimetry.

Let us also suppose that the relations in the presentation $R$ of $G$
have length at most $\lambda$.

\subsection{Conjugate words in $G$}

The goal of this section is to show that if a word $x$ is known to be a
conjugate in $G$ of a short word $y$, then some cyclic permutation of $x$
is conjugate to $y$ by a short word. If $x=uyu\~$, we will say that $x$
is \emph{conjugate to $y$ by $u$}, or that $u$ \emph{conjugates} $x$ and
$y$, or that $u$ is a \emph{conjugating word}. We recall the

\begin{defi*}
A word $w$ is said to be \emph{cyclically geodesic} if it and all of its cyclic
permutations label geodesic words in $G$.
\end{defi*}

The following is well-known (cf.~\cite{BH}, p.~452, where the authors use
``fully reduced'' for ``cyclically geodesic'').

\begin{prop}
\label{frconj}
Let $u$, $v$ be cyclically geodesic words representing conjugate elements of
$G$. Then
\begin{itemize}
\item either $\abs{u}\leq 8\delta+1$ and $\abs{v}\leq 8\delta+1$
\item or else there exist cyclic permutations $u'$ and $v'$ of $u$ and
$v$ which are conjugate by a word of length at most $2\delta+1$.
\end{itemize}
\end{prop}

This immediately extends to:

\begin{prop}
\label{frconj2}
Let $u$, $v$ be cyclically geodesic words representing conjugate elements of
$G$. Then
\begin{itemize}
\item either $\abs{u}\leq 8\delta+1$ and $\abs{v}\leq 8\delta+1$
\item or else there exist a cyclic permutation $v'$ of $v$
which is conjugate to $u$ by a word of length at most $4\delta+1$.
\end{itemize} 
\end{prop}

\begin{dem}
Write $u=u'u''$ and $v=v'v''$ such that the cyclic conjugates $u''u'$ and
$v''v'$ are conjugate by a word $\delta_1$ of length at most $2\delta+1$
as in Proposition~\ref{frconj}. Construct the quadrilateral
$u''u'\delta_1{v'}\~{v''}\~\delta_1\~$. As $u$ and $v$ are cyclically
geodesic, the sides $u''u'$ and $v''v'$ are geodesic, and in this
$\delta$-hyperbolic quadrilateral any point on one side is
$2\delta$-close to some other side. In particular, any point on the side
$u''u'$ is $(2\delta+\abs{\delta_1})$-close to the side $v''v'$.

\begin{center}
\includegraphics{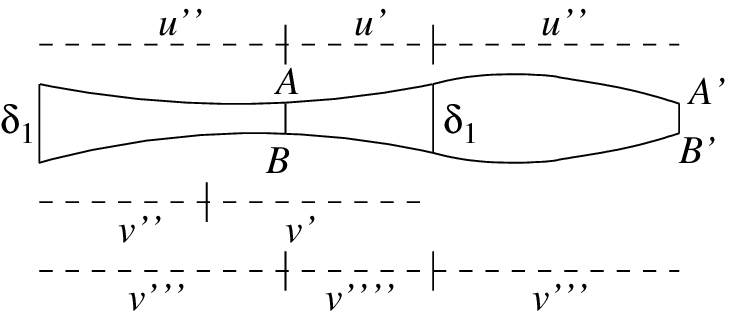}
\end{center}

Let $A$
be the endpoint of $u''$. The point $A$ is
$(2\delta+\abs{\delta_1})$-close to some point $B$ on $v''v'$. Let
$\delta_2$ be a path connecting $A$ to $B$. The point $B$ divides $v''v'$
into two words $v'''$ and $v''''$, and we have
$u=u'u''=\delta_2v''''v'''\delta_2\~$ which ends the proof of the
proposition.
\end{dem}

We will need the following

\begin{prop}
\label{fullyreduced}
Let $w$ be a geodesic word. There exists a cyclically geodesic word $z$
which is conjugate to
$w$ by a word of length at most $(\abs{w}-\abs{z})(\delta+1/2)+4\delta$.
\end{prop}

\begin{dem} 
Set $w_0=w$ and construct a sequence $w_n$ of geodesic words by induction. If $w_n$ is cyclically geodesic, stop. If not, then write $w_n=w_n'w_n''$ such that
$w_n''w_n'$ is not geodesic. Then set $w_{n+1}$ to a geodesic word equal
to $w_n''w_n'$. As length decreases at least by $1$ at each step, the
process stops after a finite number $n$ of steps and $w_n$ is cyclically geodesic.
Note that $n\leq \abs{w}-\abs{w_n}$.

In the Cayley graph of the group, define $W_i$ to be the quasi-geodesic
$(w_0' w_1'\ldots \linebreak[0] w_{i-1}' w_i^k)_{k\in \Z}$ with $w_i'$ as above:

\begin{center}
\includegraphics{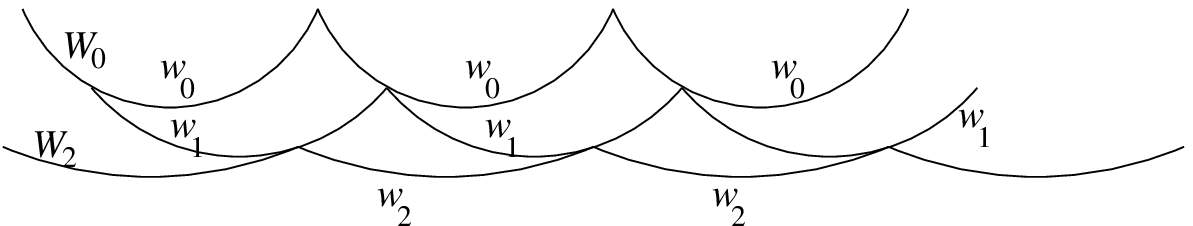}
\end{center}

Consider any of the geodesic triangles made by $w_i$, $w_{i-1}''$,
$w_{i-1}'$. As these are $\delta$-hyperbolic, this means that any point of
$W_i$ is $\delta$-close to the line $W_{i-1}$. Thus, any
point of $W_n$ is $n\delta$-close to $W_0$.

Consider the two endpoints of a copy of $w_n$ lying on $W_n$. These two points
are $n\delta$-close to $W_0$, and since the whole picture is invariant by
translation, this means that we can find a word $u$ of length at most
$n\delta$ such that $u$ conjugates $w_n$ to some cyclic conjugate $w''w'$
of $w$. Now construct the hexagon $w''w'uw_n\~ u\~$.

\begin{center}
\includegraphics{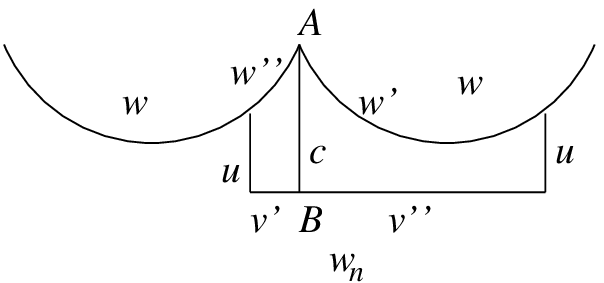}
\end{center}

Let $A$ be the endpoint of $w''$. By elementary $\delta$-hyperbolic
geometry (approximation by a tripod of the triangle consisting of $A$ and
the endpoints of $v$), the distance of $A$ to the side $v$ is at most
$(\abs{w''}+\abs{w'}+2\abs{u}-\abs{w_n})/2+4\delta$. Let $B$ be a point on
side $w_n$ realizing this minimal distance. Let $w_n=v'v''$ such that the
endpoint of $v'$ is $B$. Let $c$ be the word defined by $AB$. Then we
have $w'w''=cv''v'c\~$, so $w$ is conjugate to a cyclic conjugate of $w_n$
by $c$. Taking $z=v''v'$ ends the proof of the proposition.
\end{dem}

Now, in the spirit of Proposition~\ref{frconj}, let $C_c=\max_{x,y} \min
\{\abs{u}, x=uyu\~\}$ where the range of the maximum is the set of all
couples of conjugate words of length at most $8\delta+1$. As this set is
finite we have $C_c<\infty$. Let $C'_c=C_c+4\delta^2+12\delta+2$.

\begin{prop}
\label{gconj}
Let $x$ be a geodesic word and $y$ a conjugate of $x$ of minimal length.
Then some cyclic conjugates of $x$ and $y$ are conjugate by
a word of length at most $C'_c$.
\end{prop}

\begin{dem}
Let $u$ be a conjugating word of minimal length: $x=uyu\~$. This defines
a van Kampen diagram $ABCD$ whose sides are labeled by $u$, $y$, $u\~$ and
$x\~$ in this order.

As $x$, $y$ and $u$ are geodesic words (by minimality assumption), the
$1$-skeleton of this diagram embeds in the Cayley graph of the group, and
we get a hyperbolic quadrilateral $ABCD$ in which every point on any side is
$2\delta$-close to a point on another side.

As a first case, suppose that every point on the side $AB$ is
$2\delta$-close to either $AD$ or $BC$.

\begin{center}
\includegraphics{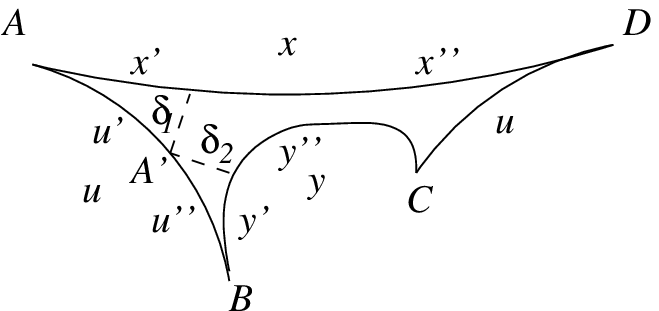}
\end{center}

Let $A'$ be the first point on $AB$ which is $2\delta$-close to $BC$.
Considering the point just before $A'$, we know that $A'$ is
$(2\delta+1)$-close to $AD$.

Then we can write $x=x'x''$, $u=u'u''$ and $y=y'y''$ such that there
exist words $\delta_1$ and $\delta_2$ of length at most $2\delta+1$ such
that $u'=x'\delta_1$ and $u''=\delta_2 {y'}\~$. Then, we have
$x''x'={x'}\~xx'=\delta_1 {u'}\~
uyu\~u'{\delta_1}\~=\delta_1u''y{u''}\~\delta_1\~=\delta_1\delta_2y''y'\delta_2\~\delta_1\~$,
and the cyclic conjugate $x''x'$ of $x$ is conjugate to $y''y'$ by a word
of length at most $4\delta+2$.

By symmetry the same tricks work if $DC$ is close to $DA$ or to $CB$.

Second, if this is not the case, let $A_n$ and $D_n$ be the points on
$AB$ and $DC$ at distance $n$ away from $A$ and $D$, respectively. Let
$n$ be smallest such that either $A_n$ or $D_n$ is not $2\delta$-close to
$AD$ nor to $BC$. By symmetry, let us suppose it is $A_n$ rather than
$D_n$. Let $w$ be a geodesic word joining $A_n$ to $D_n$.

\begin{center}
\includegraphics{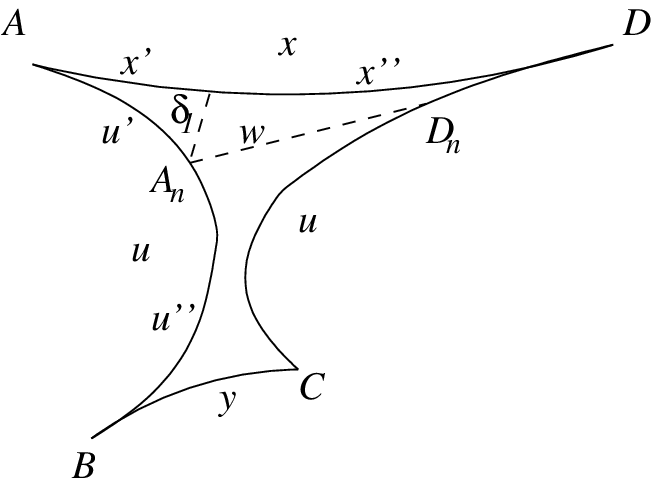}
\end{center}

Let $u'$ be the prefix of $u$ joining $A$ to $A_n$. By definition of $n$
the point $A_n$ is $2\delta+1$-close to $AD$. We have
$u'=x'\delta_1$ where $x'$ is a prefix of $x$, and $\abs{\delta_1}\leq
2\delta+1$. Thus $x''x'$ is conjugate to $w$ by a word of length at most
$2\delta+1$.

Now let us work in $A_nBCD_n$. By definition of $A_n$,
we know there exists a point $A'$ on
$CD_n$ such that $A_n A'\leq 2\delta$. Now we have $A_nD_n\leq
2\delta+A'D_n=2\delta+D_nC-A'C=2\delta+A_nB-A'C \leq 4\delta+A'B-A'C\leq
4\delta+BC$. Thus $\abs{w}\leq 4\delta+\abs{y}$.

By our minimality assumption, $y$ is cyclically geodesic. If $w$ is
cyclically geodesic as well, then we conclude by Proposition~\ref{frconj2}. If not,
use Proposition~\ref{fullyreduced} to find a cyclically geodesic word $z$ which is
conjugate to $w$ by a
word of length at most $(\abs{w}-\abs{z})(\delta+1/2)+4\delta$. By our minimality
assumption on $y$, we have that $\abs{z}\geq \abs{y}$, hence
$\abs{w}-\abs{z}\leq \abs{w}-\abs{y}\leq 4\delta$. Now $z$ and $y$ are both cyclically geodesic and
we conclude by Proposition~\ref{frconj2}.
\end{dem}

\begin{cor}
\label{conj}
Let $x$ be any word and $y$ be a conjugate of $x$ of minimal length. Then
some cyclic conjugates of $x$ and $y$ are conjugate by a word of length at
most $\delta \log_2 \abs{x}+C'_c+1$.
\end{cor}

\begin{dem}
This is because in a hyperbolic space, a geodesic joining the ends of any
curve of length $\ell$ stays at distance at most $1+\delta\log_2\ell$ from this
curve (cf.~\cite{BH}, p.~400). Take a geodesic word $x'$ equal to $x$ and apply
the above proposition; then any cyclic permutation of $x'$ will be 
conjugate to a cyclic permutation of $x$ by a word of length at most
$1+\delta\log_2\abs{x}$.
\end{dem}

\subsection{Cyclic subgroups}

We will also need the following.

\begin{prop}
\label{cyclsubgeod}
There exists a constant $R$ such that, for all hyperbolic $u\in G$, the
Hausdorff distance between the set $(u^n)_{n\in \Z}$ and any geodesic
with the same limit points is at most $\norm{u}+R$.
\end{prop}

\begin{dem}

\begin{lem}
The Hausdorff distance between $(u^n)_{n\in \Z}$ and any geodesic with
the same limit points is finite.
\end{lem}

\begin{dem}[ of the lemma]
From~\cite{GH} (p.~150) we know that $k\mapsto(u^k)_{k\in \Z}$ is a
quasi-geodesic. From ~\cite{GH} (p.~101) we thus know that this
quasi-geodesic lies at finite Hausdorff distance from some geodesic.
From~\cite{GH} (p.~119) we know that any two geodesics with the same limit
points lie at finite Hausdorff distance.
%
%
\end{dem}

Now for the proposition. First, suppose that $u$ is cyclically geodesic.
Let $p$ be a geodesic path joining $e$ to $u$. Let $\Delta$ be the union
of the paths $u^np$, $n\in \Z$. Since $u$ is cyclically geodesic, $\Delta$ is a
$(1,0,\norm{u})$-local quasi-geodesic (notation as in~\cite{GH}). Thus, there
exist constants $R$ and $L$ depending only on $G$ such that, if
$\norm{u}\geq L$, then $\Delta$ lies at Hausdorff distance at most $R$
of some geodesic $\Delta'$ equivalent to it 
(see~\cite{GH}, p.~101), 
hence at Hausdorff distance $16\delta+R$ of any other equivalent geodesic
(\cite{GH}, p.~119).
As there are only a finite number of $u$'s such that
$\norm{u}<L$, and as for each of them the lemma states that $\Delta$ lies
at finite Hausdorff distance from any equivalent geodesic, we are done
when $u$ is cyclically geodesic.

If $u$ is not cyclically geodesic, apply Proposition~\ref{gconj} to get a
cyclically geodesic word $v$ such that $v=xu''u'x\~$ with $u=u'u''$ and
$\abs{x}\leq C'_c$. Apply the above to $(v^k)_{k\in \Z}$: this set stays
at distance at most $R$ of some geodesic $\Delta$. Translate by $u'x\~$.
The set $(u'x\~v^k)_{k\in \Z}$ stays at distance at most $R$ of the
geodesic $u'x\~\Delta$. But since $u^k=u'x\~v^kx{u'}\~$, 
the Hausdorff distance between the sets $(u^k)_{k\in \Z}$ and
$(u'x\~v^k)_{k\in \Z}$ is at most $\norm{x{u'}\~}\leq C'_c+\norm{u}$ and
we are done.
\end{dem}

Since the stabilizer of any point of the boundary is either finite or has
$\Z$ as a finite index subgroup (cf~\cite{GH}, p.~154), we get as an immediate by-product of the
lemma

\begin{cor}
\label{diststab}
Let $\Delta$ be a geodesic in $G$, with limit points $a$ and $b$.
There exists a constant $R(\Delta)$ such that for any $x$ in the
stabilizer of $a$ and $b$, the distance from $x$ to $\Delta$ is at most
$R(\Delta)$.
\end{cor}

\subsection{One-hole diagrams}

We now turn to the study of isoperimetry of van Kampen diagrams with exactly
one hole. Recall that conjugacy of two words $u$ and $v$ is equivalent to
the existence of a one-hole van Kampen diagram bordered by $u$ and $v$.

\begin{prop}
\label{iso1}
There exists a constant $C'>0$ such that for any two conjugate words $u$
and $v$, there exists a one-hole diagram $D$ bordered by $u$ and
$v$, such that $C'\abs{D}\leq \abs{u}+\abs{v}$.
\end{prop}

\begin{dem}
Let us first suppose that $u$ and $v$ are geodesic words. Let $w$ be the
shortest common conjugate of $u$ and $v$.
By Proposition~\ref{gconj}, $u$ and $w$ are conjugate by a word $x$ of
length at most $\abs{u}/2+\abs{w}/2+C'_c$. Thus, there exists a minimal
van Kampen diagram $D$ bordered by $wx\~u\~x$. It follows from the isoperimetry
in $G$ that $\abs{D}\leq (\abs{u}+\abs{w}+2\abs{x})/C$. As $\abs{w}\leq
\abs{u}$ we have $\abs{D}\leq \abs{u}(4+2C'_c)/C$.

Do the same job with $v$ and $w$, to get a diagram $D'$ bordered by
$v\~y\~wy$. Then paste these two diagrams along the $w$'s, getting a
diagram bordered by $v(xy)\~u\~(xy)$. Then transform this diagram into an
annulus by gluing the two $xy$ sides; this leads to a one-hole diagram
bordered by $u$ and $v$. The number of its faces is at most
$(\abs{u}+\abs{v})(4+2C'_c)/C$ and we conclude by setting
$C'=C/(4+2C'_c)$ in case $u$ and $v$ are geodesic.

\begin{center}
\includegraphics{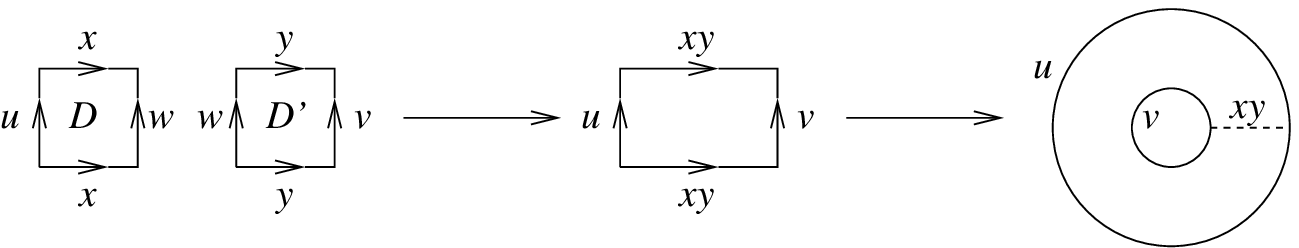}
\end{center}

In case $u$ is not geodesic, let $u'$ be a geodesic word equal to $u$ in
$G$. We know there exists a van Kampen diagram $D_u$ bordered by $u{u'}\~$,
with $\abs{D_u}\leq 2\abs{u}/C$. Let
$D_v$ be a similar diagram for $v$. Let $D$ be as above a one-hole
minimal diagram bordered by $u'$ and $v'$, with $\abs{D}\leq
(\abs{u}+\abs{v})/C'$ with $C'$ as above. Then we can glue $D_u$ and
$D_v$ to $D$ along their common boundaries.

\begin{center}
\includegraphics{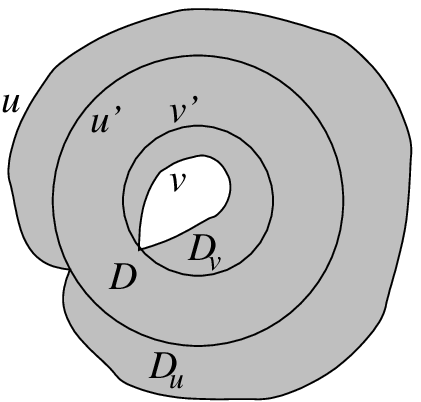}
\end{center}

This leads to a diagram with at most
$(\abs{u}+\abs{v})/C'+2(\abs{u}+\abs{v})/C$ faces, and we conclude by
re-setting $C'$ to $1/(1/C'+2/C)$.
\end{dem}

\subsection{Narrowness of diagrams}

We now prove that diagrams (with or without holes) in a hyperbolic
space are narrow (see section~\ref{defs} for definitions).

Let $\alpha=1/\log(1/(1-C'/\lambda))$ where $C'$ is given by
Proposition~\ref{iso1}. (Recall $\lambda$ is the maximal
length of relators in the presentation of $G$.) Let $\lceil x\rceil$
denote the integer part of $x$ plus one (such that $\lceil \log
\abs{D}\rceil=1$ for $\abs{D}=1$).

\begin{prop}
\label{iso01}
Let $D$ be a minimal diagram with either $0$ or $1$ hole. Then $D$
is $\lceil \alpha\log \abs{D}\rceil$-narrow.
\end{prop}

\begin{dem}
Let $D$ be a minimal van Kampen diagram with $0$ or $1$ hole.
Proposition~\ref{iso1} tells us that
$C' \abs{D}\leq \abs{\d D}$. Let $n$ be the number of faces of $D$ lying
on the boundary. We have $\abs{\d D}\leq \lambda n$. Thus the proportion of
faces of $D$ lying on the boundary is at least $C'/\lambda$.

Let $D'$ be the diagram $D$ with the boundary faces removed. (In case $D'$
is not connected, consider any one of its connected components.) $D'$ has at most one
hole. $D'$ is minimal as a subdiagram of a minimal diagram. Furthermore,
we have $\abs{D'}\leq \abs{D} (1-C'/\lambda)$. By the same argument, the
proportion of boundary faces of $D'$ is at least $C'/\lambda$, and after
removing these faces we get a third diagram $D''$ with at most $\abs{D}
(1-C'/\lambda)^2$ faces. Repeating the argument yields the desired
conclusion as $D$ is exhausted after
$\log\abs{D}/\log(1/(1-C'/\lambda))$ steps.
\end{dem}

\begin{prop}
\label{isohole}
Let $D$ be a minimal $n$-hole diagram. Then $D$ satisfies the
isoperimetric inequality
\[
\abs{\d D}\geq C\abs{D}-n\lambda \left(2+4\lceil \alpha \log
\abs{D}\rceil\right)
\]
\end{prop}

\begin{dem}

\begin{lem}
\label{lemisohole}
Let $D$ be a minimal $n$-hole van Kampen diagram ($n\geq 1$). Either there exists a path in the $1$-skeleton of $D$ joining two holes, with length at most
$\lambda(1+2\lceil\alpha\log \abs{D}\rceil)$, or there exists a path in
the $1$-skeleton of $D$ joining one hole with the exterior boundary, with
length at most $\lambda(1/2+\lceil\alpha\log \abs{D}\rceil)$.
\end{lem}

\begin{dem}[ of the lemma]
We work by induction on $n$. Set $e=\lceil\alpha\log\abs{D}\rceil$.

Observe that a chain of $N$ adjacent faces provides a path of length at
most $N\lambda/2$ in the $1$-skeleton between any two vertices of these
faces.

For $n=1$, the lemma is clear: by the last proposition, the diagram is
$e$-narrow, thus the two components of the boundary are linked by a chain
of at most $2e$ faces, providing a path of length at most $\lambda e$.

Now suppose the lemma is true up to some $n\geq 1$, and let $D$ be a
$(n+1)$-hole van Kampen diagram. For every hole $i$, let $B_i$ be the set of
faces of $D$ lying at distance at most $2e+1$ from the boundary of $i$.

Either, first, there are holes $i\neq j$ such that $B_i$ and $B_j$ have a
common face or edge or vertex. This provides a chain of at most $4e+2$
faces between the boundaries of holes $i$ and $j$, thus a path of length
at most $\lambda(2e+1)$.

Or, second, the $B_i$'s do not meet. Choose any hole $i$.

There can be holes in $B_i$, different from $i$, that can be
filled in $D$. Define $B'_i$ as $B_i$ plus the interiors of these holes in
$D$, in such a manner that all holes of $B'_i$ are holes of $D$.

First, suppose that $B'_i$ does not encircle any hole $j$ of $D$ other
than $i$. As $B_i$ is defined as the ball of
radius $2e+1$ around $i$ in $D$, any face on the exterior boundary of
$B'_i$ is either a face at distance $2e+1$ from hole $i$, or a face on
the boundary of $D$. But as $B'_i$ is a one-hole van Kampen diagram included in $D$,
hence $e$-narrow by Proposition~\ref{iso01}, not all faces of the
exterior boundary of $B'_i$ can be at distance $2e+1$ from $i$. That is, at
least one face of the exterior boundary of $B'_i$ is on the exterior boundary of
$D$, hence a path of length at most $\lambda(e+1/2)$.

Second, imagine that $B'_i$ encircles at least one hole $j\neq i$ of $D$.
Consider the part $D'$ of $D$ comprised between $B'_i$ and $j$, that is,
the connected component of $D\setminus B'_i$ containing $j$. This is a
diagram with at least one hole $j$ (and maybe others), but as it does not
contain $i$ it has at most $n$ holes. As $D$ is minimal, $D'$ is. By the
induction assumption, either two holes in $D'$ are at distance at most
$\lambda(2e+1)$, in which case we are done, or one hole, say $j$, in $D'$
is at distance at most $\lambda(e+1/2)$ of the exterior boundary of $D'$.
But the exterior boundary of $D'$ is part of the boundary of $B'_i$, any
point of which is at distance at most $\lambda(e+1/2)$ of hole $i$. Thus
$i$ and $j$ are linked by a path of length at most $\lambda(2e+1)$, which
ends the proof of the lemma.
\end{dem}

\begin{cor*}[ of Lemma~\ref{lemisohole}]
A minimal $n$-hole diagram can be made simply connected by cutting it
along $n$ curves of cumulated length at most
$n\lambda(2\lceil\alpha\log\abs{D}\rceil\linebreak[0]+1)$.
\end{cor*}

The corollary of the lemma ends the proof of the proposition.
\end{dem}

\begin{cor}
\label{narrowhole}
A minimal $n$-hole diagram $D$ is
$\lceil\alpha\log\abs{D}\rceil+n(4\lceil\alpha\log\abs{D}\rceil+2)$-narrow.
\end{cor}

\begin{dem}
Let $D'$ be a simply connected van Kampen diagram resulting from cutting $D$ along
curves of cumulated length at most
$n\lambda(2\lceil\alpha\log\abs{D}\rceil+1)$ (which run along at most
$n(4\lceil\alpha\log\abs{D}\rceil+2)$ faces as can immediately be seen on
the proof above). Every face in the new diagram is at distance $\lceil
\alpha\log\abs{D}\rceil$ from the boundary of $D'$ by
Proposition~\ref{iso01}. The boundary of $D$
is a subset of the boundary of that of $D'$, but by construction any
face on the boundary of $D'$ is at distance at most
$n(4\lceil\alpha\log\abs{D}\rceil+2)$ from the boundary of $D$.
\end{dem}

\subsection{Coarsenings of diagrams}

If $D$ is a very narrow diagram with holes, then we have an intuitive
feeling of which parts of its boundary ``face'' which. This intuition can
be made clear using the approximation of hyperbolic spaces by trees. We
now give for this intuition a mathematical setting fitted to our needs.

\begin{defi}
\label{defimatching}
Let $w_1,\ldots, w_n$ be $n$ geodesic words in $G$. A
\emph{$(k,\eps)$-matching}
of these words is a set of words $w'_1,\ldots,w'_{k}$, some of which may be
empty, together with a partition $I_1,I_2$ of $\{1,\ldots,k\}$ and a
bijection $\psi$ between $I_1$ and $I_2$, such that:
\begin{itemize}
\item The words $w'_1,\ldots,w'_k$ form a partition of the words
$w_1,\ldots,w_n$.
\item For all $i\in I_1$, there exist words $\delta_1$ and $\delta_2$ of
length at most $\eps$ such that $w'_i=\delta_1 w'_{\psi(i)} \delta_2$ in
$G$ (we will say that $w'_i$ and $w'_{\psi(i)}$ \emph{$\eps$-match}).
\end{itemize}
\end{defi}

This means that we cut the words $w_i$ into at most $k$ subwords such that
each subword ``faces'' another one up to $\eps$. Typically $\eps$ is of
order $\delta$. We have to leave open the possibility that some $w'_i$
are empty since, for example, if one of the $w_i$'s is very short, it
could have to be matched with the empty word.

The following proposition is basically equivalent to the approximation of
finite hyperbolic spaces by trees.

\begin{prop}
Let $w_1,\ldots, w_n$, for $n\geq 2$, be $n$ geodesic words in $G$ such
that $w_1\ldots w_n=e$. There exists a $(4n,4n\delta)$-matching of these words.
\end{prop}

For $n=3$ this closely resembles the definition of thin triangles.

\begin{dem}
Work by induction on $n$. The result is clear for $n=2$.
Suppose that $n+1$ words $w_1,\ldots,w_n,w_{n+1}$ forming a closed piecewise
geodesic path in $G$ are given. Let $x$ be a geodesic word equal
to $w_nw_{n+1}$.  The three geodesic words $x$, $w_n$, $w_{n+1}$ form a
$\delta$-thin triangle. Let $x=x_1x_2$ where the endpoint of $x_1$ lies
at distance at most $2\delta$ from both sides $w_n$ and $w_{n+1}$ of the
triangle.  Now apply the induction assumption to $w_1,\ldots,w_{n-1},x$.
This gives a matching involving a partition of the word $x$ into subwords
$x'_i$, $i\in I$. At most one of the words $x'_i$ straddles the endpoint
of $x_1$.  If some $x'_i$ is included in $x_1$ or $x_2$ and is
$4n\delta$-matched to $w'$ where $w'$ is a subword of the $w_i$'s, then
using thinness of the triangle $x$, $w_n$, $w_{n+1}$ we can
$(4n+2)\delta$-match $w'$ with a subword of $w_n$ or $w_{n+1}$.  If
$x'_i$ straddles the endpoint of $x_1$, and $x'_i$ is $4n\delta$-matched
to $w'$, then we can write $x'_i=x''_i x'''_i$ where the endpoint of
$x''_i$ is that of $x_1$, and also write $w'=w''w'''$ such that $w''$
$(4n+2)\delta$-matches with $x''_i$ and $w'''$ matches with $x'''_i$;
using thinness of the triangle $x$, $w_n$, $w_{n+1}$ we can
$(4n+4)\delta$-match $w''$ and $w'''$ with subwords of $w_n$ and
$w_{n+1}$ respectively. Last, the two parts of $w_n$ and $w_{n+1}$ which
are not $2\delta$-close to $x$ can be matched together.
\end{dem}

Doing the induction more cleverly, one can even obtain a
$(4n,4\delta\log_2 n)$-matching.

\bigskip

We are to apply this construction to diagrams with $n$ holes.
In order to symmetrize the role of holes and of the boundary in the
following proposition, we view a $n$-hole diagram as a $(n+1)$-hole diagram
embedded in the sphere.

\begin{prop}
\label{coarsening}
There exists a constant $B$ (depending on $G$) such that, for any minimal
$(n+1)$-hole diagram $D$ embedded in the sphere, there exists a 
$(8n,Bn\log\abs{D})$-matching of the boundary words of $D$.
This matching is called the \emph{coarsening} of $D$.
\end{prop}

The coarsening of $D$ is best visualized as a planar graph as in the
following picture (black dots mark the points where we
partition the boundary words). The planar graph can be precisely defined
but we do not need it.

\begin{center}
\includegraphics{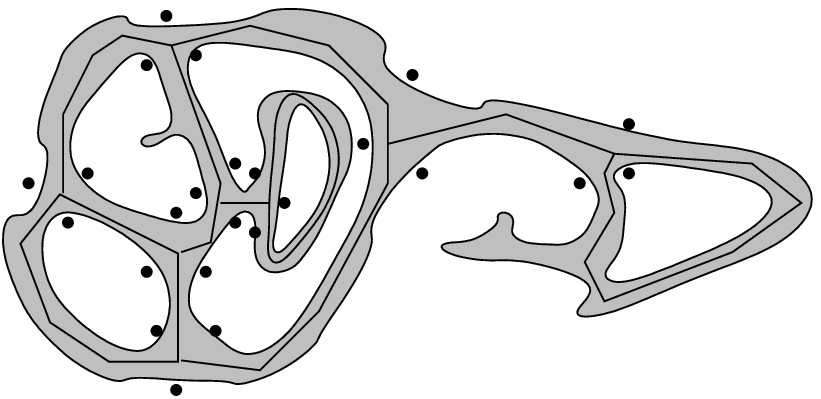}
\end{center}

\begin{dem}
First, using Corollary~\ref{conj}, and the fact already used above that
any geodesic joining the endpoints of a curve of length $\ell$ stays
$(\delta\log\ell+1)$-close to that curve, we can suppose that the
boundary words of $D$ are cyclically geodesic: the length
$\ell$ of any boundary word of $D$ is at most $\lambda\abs{D}$ and so a
$(k,\eps)$-matching for the cyclically reduced words will give a
$(k,\eps+2+C'_c+2\delta\log(\lambda \abs{D}))$-matching of the original words.

Use the corollary of Lemma~\ref{lemisohole} to cut $D$ into a simply
connected diagram $D'$. The boundary word of $D'$ is made of $n'\leq 2n$
pieces $w_1,\ldots,w_{n'}$ partitioning the boundary words of $D$, with
little words $x_1,\ldots,x_{n'}$ of cumulated length at most
$B_1n\log\abs{D}$ in between (for some constant $B_1$ depending on $G$).
Define $y_i$ to be a geodesic word equal to $w_ix_i$.

Now apply the previous proposition to get a $(4n',4n'\delta)$-matching of
the words $y_1$, $y_2$, \ldots, $y_{n'}$. Since the $x_i$'s are of
cumulated length at most $B_1n\log\abs{D}$, this matching defines a
$(4n',4n'\delta+2\delta+B_1n\log\abs{D})$-matching of the $w_i$'s.
\end{dem}

\newpage

\section{Appendix: Cases of harmful torsion}
\label{appcex}

\bigskip

Here we show that the assumption of harmless torsion cannot be removed.
Examples of hyperbolic groups with harmful torsion include such groups as
$(F_m\times \Z/2\Z)\star F_m$ with $m\geq 2$, since the $\Z/2\Z$ factor
has a centralizer which is a free group of rank $m$. More precisely:

\begin{prop}
Theorem~\ref{main} does not hold for the hyperbolic group $(F_4\times
\Z/2\Z)\star F_4$.
\end{prop}

\begin{dem}
Consider the two groups $G_1=(F_m\times \Z/2\Z)\star F_m$ and $G_2=(F_m
\star F_m)\times \Z/2\Z$. In each of these, denote by $u$ a generator for
$\Z/2\Z$ and respectively by $a_1,\ldots,a_m$ and $b_1,\ldots,b_m$ a
standard set of generators for the first and second factor $F_m$. Let
$A_1$ and $A_2$ be the subgroups of $G_1$ and $G_2$, respectively,
generated by the $a_i$'s, and define $B_1$ and $B_2$ similarly.

It is immediate to see that these groups are hyperbolic.

For any group $G$ generated by $k$ elements, let $\lambda(G)$ denote the
spectral radius of the random walk on $G$ (w.r.t.\ the $k$ generators),
and let $\theta=1+\log_{2k}\lambda$ be the gross cogrowth of $G$.

The spectral radius for the free group $F_k$ is $\lambda(F_k)=\sqrt{2k-1}/k$.
By Lemma 4.1 of~\cite{K1}, the spectral radius for
the group $F_k\times \Z/2\Z$ is equal to
$(1+k\lambda(F_k))/(k+1)=(1+\sqrt{2k-1})/(k+1)$.

In particular, the spectral radius of $G_2$ is $(1+\sqrt{4m-1})/(2m+1)$.

Take $m=4$. We have $\theta(G_2)=1+\log_{4m+2} \lambda(G_2)\approx .788$.
In particular, the critical density $d^2_{\mathrm{crit}}$ for random
quotients of $G_2$ by plain random words is about $1-.788=.212$.

Since $G_2$ is a quotient of $G_1$ we have of course $\lambda(G_1)\leq
\lambda(G_2)$. But Theorem~1 of~\cite{K1} states that quotienting a group
by (the normal closure of) a non-amenable subgroup strictly increases the
spectral radius. The kernel of the quotient map $G_1\rightarrow G_2$
contains the two elements $ub_1u^{-1}{b_1}^{-1}$ and
$ub_2u^{-1}{b_2}^{-1}$ which generate a free non-cyclic subgroup in
$\Z/2\Z\star B_1$. Hence the kernel is non-amenable.

Thus, we have $\lambda(G_1)<\lambda(G_2)$, so that if Theorem~\ref{main}
holds for
$G_1$, the critical density $d^1_{\mathrm{crit}}$ for random quotients of
$G_1$ 
satisfies \[
d^1_{\mathrm{crit}}>d^2_{\mathrm{crit}}\approx .212
\]

But we are going to prove that random quotients of $G_1$ are very
probably trivial as soon as $d>d^2_{\mathrm{crit}}$.

\bigskip

Let $R$ be a set of randomly chosen words in
$u,a_1,\ldots,a_m,b_1,\ldots,b_m$, of length $\ell$, at density $d$.
(Note that for the model of random quotients by plain random words, the
law of the relators depends only on the generators and not on the initial
group.) We now study the random quotient $G_1/\langle R\rangle$ and
consider the elements of $R$ as elements of $G_1$.

Let us compute the probability that one of the relators in $R$ belongs to
$C=A_1\times \Z/2\Z\subset G_1$. The number of words of length $\ell$
belonging to $C$ is at least $(2m+2)^{\ell}$, so that the corresponding
density is at least $\log_{4m+2}(2m+2)\approx .797$. So there exists a
density $d_C\leq 1-.797=.203$ such that if $d>d_C$, there will very
probably be some element of $R$ lying in $C$. Note that this is below the
critical density $d^1_{\mathrm{crit}}$ for random quotients of $G_1$
predicted by Theorem~\ref{main} (if it holds). Also, $d_C$ is not $0$
since $G_1/\langle C\rangle$ is not amenable.

By the same argument, for $d>d_C$ it is very probable that $R$ contains
a relator $r$ of the form $r=xc$ where $x$ is one of the generators of
$G_1$ and $c$ is a word of length $\ell-1$ with $c\in C$. As the random
words are sampled uniformly, when $\ell$ is big enough this will occur
for \emph{each} of the relators $x$ of $G_1$.

Let $H$ be the random quotient $G/\langle R \rangle$. By
definition of $C$, $u$ commutes with $c$ in $G_1$, so in $H$ we have
\[
uxu^{-1}x^{-1}=uxcu^{-1}c^{-1}x^{-1}=uru^{-1}r^{-1}=e
\]
since $r=e$ in $H$ by definition.

Thus, in $H$, the generator $u$ commutes with all the generators of
$G_1$. Let $S\subset G_1$ be the set of the commutators of $u$ with these
generators, we have
\[
H=G_1/\langle R\rangle = G_1/\langle R\cup S\rangle = \left(G_1/\langle
S\rangle\right)/\langle R\rangle = G_2/\langle R \rangle
\]

But $G_2/\langle R\rangle$ is a random quotient of $G_2$ (this is because
for random quotients by plain words, the law of $R$ is independent on the
initial group). In particular, if $d>d^2_{\mathrm{crit}}\approx .212$
this group will very probably be trivial, whereas if Theorem~\ref{main} were valid
for $G_1$, the critical value would be
$d^1_{\mathrm{crit}}>d^2_{\mathrm{crit}}$.

This ends the proof.
\end{dem}

\bigskip

So random quotients of $G_1$ behave in a different manner than that of
Theorem~\ref{main}. For densities between $0$ and $d_C<.203$ they behave
``normally'' (in particular, Axiom~4 is satisfied). But for
densities between $d_C$ and $.212$, Axiom~4 is not satisfied, and the
random quotients behave like random quotients of $G_2$, and they vanish
as soon as $d>.212$, whereas the expected critical density in Theorem~\ref{main}
would be higher. (The gap between $.203$ and $.212$ can be made larger
by taking bigger $m$.)

\begin{center}
\begin{picture}(0,0)%
\includegraphics{threephases.pstex}%
\end{picture}%
\setlength{\unitlength}{3947sp}%
\begingroup\makeatletter\ifx\SetFigFont\undefined%
\gdef\SetFigFont#1#2#3#4#5{%
  \reset@font\fontsize{#1}{#2pt}%
  \fontfamily{#3}\fontseries{#4}\fontshape{#5}%
  \selectfont}%
\fi\endgroup%
\begin{picture}(5192,1528)(301,-2519)
\put(2326,-2461){\makebox(0,0)[b]{\smash{\SetFigFont{12}{14.4}{\familydefault}{\mddefault}{\updefault}{\color[rgb]{0,0,0}$ub_1=b_1u$}%
}}}
\put(1126,-2461){\makebox(0,0)[b]{\smash{\SetFigFont{12}{14.4}{\familydefault}{\mddefault}{\updefault}{\color[rgb]{0,0,0}$ub_1\neq b_1 u$}%
}}}
\put(2326,-2236){\makebox(0,0)[b]{\smash{\SetFigFont{12}{14.4}{\familydefault}{\mddefault}{\updefault}{\color[rgb]{0,0,0}Infinite}%
}}}
\put(1126,-2236){\makebox(0,0)[b]{\smash{\SetFigFont{12}{14.4}{\familydefault}{\mddefault}{\updefault}{\color[rgb]{0,0,0}Infinite}%
}}}
\put(4201,-2236){\makebox(0,0)[b]{\smash{\SetFigFont{12}{14.4}{\familydefault}{\mddefault}{\updefault}{\color[rgb]{0,0,0}Trivial}%
}}}
\put(1426,-1861){\makebox(0,0)[b]{\smash{\SetFigFont{12}{14.4}{\familydefault}{\mddefault}{\updefault}{\color[rgb]{0,0,0}Random quotients of $G_1$ are}%
}}}
\put(1726,-1336){\makebox(0,0)[b]{\smash{\SetFigFont{12}{14.4}{\rmdefault}{\mddefault}{\itdefault}{\color[rgb]{0,0,0}$d_C$}%
}}}
\put(2926,-1336){\makebox(0,0)[b]{\smash{\SetFigFont{12}{14.4}{\rmdefault}{\mddefault}{\itdefault}{\color[rgb]{0,0,0}$1-\theta(G_2)$}%
}}}
\put(3751,-1186){\makebox(0,0)[b]{\smash{\SetFigFont{12}{14.4}{\rmdefault}{\mddefault}{\itdefault}{\color[rgb]{0,0,0}$1-\theta(G_1)$}%
}}}
\end{picture}

\end{center}

The two phases are really different: indeed, a difference can be seen in
the ball of radius $2$ in the Cayley graph since, in the random quotient,
the relation $ub_1=b_1u$ (notation as above) holds in the second phase
but not in the first one (since in the first phase, the ``ordinary'' theory
of random quotients holds and in particular, the radius of injectivity
grows with $\ell$).

\bigskip

More than three phases can probably be arranged, using groups such as 
\[\left(\left(\left(F_m\times \Z/2\Z\right)\star F_p\right)\times
\Z/2\Z\right)\star F_q
\]
with different critical densities equal to the
densities of the centralizers of the different torsion elements.

\newpage
\tableofcontents

\end{document}